\documentclass[12pt]{article}
\usepackage{amsfonts}
\usepackage{amssymb}
\usepackage{latexsym}
\usepackage{amsmath}
\usepackage{bibmods}
\usepackage{graphicx,epsfig}
\usepackage{color}
\usepackage[arrow,matrix,curve]{xy}
\usepackage{amssymb}
\usepackage{amscd}
\usepackage{color}

\newtheorem{lemma}{\indent Lemma}[section]
\newtheorem{theorem}[lemma]{\indent Theorem}
\newtheorem{corollary}[lemma]{\indent Corollary}
\newtheorem{proposition}[lemma]{\indent Proposition}

\newtheorem{definition}[lemma]{\indent Definition}
\newtheorem{example}[lemma]{\indent Example}
\newtheorem{remark}[lemma]{\indent Remark}

\newcommand{\bt}{\beta}

\newcommand{\vf}{\varphi}
\newcommand{\ve}{\varepsilon}

\newcommand{\bC}{\mathbb{C}}
\newcommand{\bH}{\mathbb{H}}

\newcommand{\bR}{\mathbb{R}}
\newcommand{\bW}{\mathbb{W}}

\newcommand{\cA}{{\mathcal A}}

\newcommand{\cF}{{\mathcal F}}

\newcommand{\cK}{{\mathcal K}}
\newcommand{\cL}{{\mathcal L}}
\newcommand{\cM}{{\mathcal M}}

\newcommand{\ov}{\overline}
\newcommand{\wt}{\widetilde}
\newcommand{\dR}{\overset{\text{\tiny$\bullet$}}{\mathbb{R}}}

\newcommand{\A}{{\boldsymbol A}}

\newcommand{\K}{{\boldsymbol K}}
\newcommand{\M}{{\boldsymbol M}}
\newcommand{\N}{{\boldsymbol N}}

\newcommand{\V}{{\boldsymbol V}}

\newcommand{\IIm}{\operatorname{Im}}
\newcommand{\RRe}{\operatorname{Re}}
\newcommand{\ind}{\operatorname{ind}}
\newcommand{\Ind}{\operatorname{Ind}}

\newcommand{\sign}{\operatorname{sign}}
\newcommand{\supp}{\operatorname{supp}}
\newcommand{\QED}{\hspace{\fill}$\Box$\medskip\par}

\begin{document}

\begin{center}
{\Large\textbf{MELLIN CONVOLUTION OPERATORS\\[-0.5mm] IN BESSEL POTENTIAL SPACES WITH\\[1.5mm] ADMISSIBLE MEROMORPHIC KERNELS}\footnote{This work was supported by the Shota Rustaveli Georgian National Science
Foundation, Contract No. 13/14.}}
\vspace{5mm}
 \end{center}
\begin{center}
\ \hfill \textit{Dedicated to the memory of Academician Victor Kupradze
on the occasion of his 110-th birthday anniversary}
\vspace{3mm}

\textbf{R. Duduchava}

\vspace{2mm}

I. Javakhishvili Tbilisi State University, Andrea Razmadze
Mathematical Institute,  Tamarashvili str. 6, Tbilisi 0177, Georgia;
\texttt{roland.duduchava@tsu.ge}
 \end{center}


\thispagestyle{empty}

\begin{abstract}
The paper is devoted to Mellin convolution operators with meromorphic kernels in Bessel potential spaces. We encounter such operators while investigating boundary value problems for elliptic equations in planar 2D domains with angular points on the boundary.

Our study is based upon two results. The first concerns commutants of Mellin convolution and Bessel potential operators: Bessel potentials alter essentially after commutation with Mellin convolutions depending on the poles of the kernel (in contrast to commutants with Fourier convolution operatiors.) The second basic ingredient is the results on the Banach algebra $\mathfrak{A}_p$ generated by Mellin convolution and Fourier convolution operators in weighted $\mathbb{L}_p$-spaces obtained by the author in 1970's and 1980's. These results are modified by adding Hankel operators. Examples of Mellin convolution operators are considered.

The first version of the paper was published in  {\em Memoirs on Differential Equations and Mathematical Physics} {\bf 60},  135-177, 2013. The formulations and proofs there contain fatal errors, which are improved in the present preprint. Part of the results, obtained with V. Didenko, are published in the preprint \cite{DiDu15}.
\end{abstract}

\newpage
\tableofcontents

\section*{Introduction}
\addcontentsline{toc}{section}{Introduction}

It is well-known that various boundary value problems for PDE in planar domains with angular points on the boundary, e.g. Lam\'e systems in elasticity (cracks in elastic media, reinforced plates), Maxwell's system and Helmholtz equation in electromagnetic scattering, Cauchy--Riemann systems, Carleman--Vekua systems in generalized analytic function theory etc. can be studied with the help of the Mellin convolution equations of the form
\begin{equation}\label{e0.1}
    \mathbf{A}\varphi(t):=c_0\varphi(t)+\frac{c_1}{\pi i}\int\limits_0^\infty \frac{\varphi(\tau)\,dt}{\tau-t}+
     \int\limits_{0}^\infty\mathcal{K}\Big(\frac t\tau\Big)\varphi(\tau)\,\frac{d\tau}{\tau}=f(t),
\end{equation}
with the kernel $\mathcal{K}$ satisfying the condition
\begin{equation}\label{e0.2}
    \int\limits_0^\infty t^{\beta-1}|\mathcal{K}(t)|\,dt<\infty, \;\; 0<\beta<1,
 \end{equation}
which makes it a bounded operator in the weighted Lebesgue space \linebreak  $\mathbb{L}_p(\mathbb{R}^+,t^\gamma)$, provided $1\leqslant p\leqslant\infty$, $-1<\gamma<p-1$, $\beta:=(1+\gamma)/p$ (cf.~\cite{Du79}).

In particular, integral equations with fixed singularities in the kernel
\begin{multline}\label{e0.3}
\hskip-5mm c_0(t)\varphi(t)+\frac{c_1(t)}{\pi i} \int\limits_0^\infty\frac{\varphi(\tau)\,dt}{\tau-t}
    +\sum_{k=0}^n \frac{c_{k+2}(t)t^{k-r}}{\pi i} \int\limits_0^\infty \frac{\tau^r\varphi(\tau)\,d\tau}{(\tau+t)^{k+1}}=f(t), \;\;
                    0\leqslant t\leqslant1,
\end{multline}
where $0\leqslant r\leqslant k$ are of type \eqref{e0.1} after localization, i.e. after ``freezing'' the coefficients.

The Fredholm theory and the unique solvability of equations \eqref{e0.1} in the weighted Lebesgue spaces were accomplished in \cite{Du79}. This investigation was based on the following observation: if $1<p<\infty$, $-1<\gamma<p-1$, $\beta:=(1+\gamma)/p$, the following mutually invertible exponential transformations
\begin{equation}\label{e0.4}
\begin{gathered}
    Z_\beta:\mathbb{L}_p([0,1],t^\gamma)\longrightarrow\mathbb{L}_p(\mathbb{R}^+), \\
    Z_\beta\varphi(\xi):=e^{-\beta\xi}\varphi(e^{-\xi}),\;\; \xi\in\bR:=(-\infty,\infty), \\
    Z^{-1}_\beta:\mathbb{L}_p(\mathbb{R}^+)\longrightarrow\mathbb{L}_p([0,1],t^\gamma), \\
    Z^{-1}_\beta\psi(t):=t^{-\beta}\psi(-\ln\,t), \;\; t\in\mathbb{R}^+:=(0,\infty),
\end{gathered}
\end{equation}
transform the equation \eqref{e0.1} from the weighted Lebesgue space  $f,\varphi\in\mathbb{L}_p(\mathbb{R}^+,t^\gamma)$ into the Fourier convolution equation $W^0_{\mathcal{A}_\beta}\psi=g$, $\psi=Z_\beta\varphi, g=Z_\beta f\in\mathbb{L}_p(\mathbb{R})$ of the form
\allowdisplaybreaks
\begin{align*}
    W^0_{\mathcal{A}_\beta}\psi(x) & =c_0\psi(x)+\int\limits\limits_{-\infty}^\infty \mathcal{K}_1(x-y)\varphi(y)\,dy, \\
    \mathcal{K}_1(x) & =e^{-\beta x}\Big[\frac{c_1}{1-e^{-x}}+\mathcal{K}(e^{-x})\Big].
 \end{align*}
Note that the symbol of the operator $W^0_{\mathcal{A}_\beta}$, viz. the Fourier transform of the kernel
\begin{align}\label{e0.5}
    \mathcal{A}_\beta(\xi) & :=c_0+\int\limits_{-\infty}^\infty e^{i\xi x}\mathcal{K}_1(x)\,dx\nonumber\\
    & :=c_0-ic_1\cot\pi(\beta-i\xi)+\int\limits_{-\infty}^\infty e^{(i\xi-\beta)x}\mathcal{K}(e^{-x})\,dx, \;\; \xi\in\mathbb{R}
\end{align}
is a piecewise continuous function. Let us recall that the theory of Fourier convolution operators with discontinuous symbols is well developed, cf.   \cite{Du75a, Du75b, Du77, Du78, Th85}. This allows one to investigate various properties of the operators \eqref{e0.1}, \eqref{e0.3}. In particular, Fredholm criteria, index formula and conditions of unique solvability of the equations \eqref{e0.1} and \eqref{e0.3} have been established in \cite{Du79}.

Similar integral operators with fixed singularities in kernel arise in the theory of singular integral equations with the complex conjugation
$$  a(t)\varphi(t)+\frac{b(t)}{\pi i}\int\limits_\Gamma\frac{\varphi(\tau)\,dt}{\tau-t}+
        \frac{e(t)}{\pi i}\overline{\int\limits_\Gamma\frac{\overline{\varphi(\tau)}\,dt}{\tau-t}}=f(t), \;\; t\in\Gamma $$
and in more general R-linear equations
\begin{multline*}
    a(t)\varphi(t)+b(t)\overline{\varphi(t)}+\frac{c(t)}{\pi i}\int\limits_\Gamma\frac{\varphi(\tau)\,dt}{\tau-t}+
        \frac{d(t)}{\pi i}\int\limits_\Gamma\frac{\overline{\varphi(\tau)}\,dt}{\tau-t}+ \\
    +\frac{e(t)}{\pi i}\overline{\int\limits_\Gamma\frac{\varphi(\tau)\,dt}{\tau-t}}+
        \frac{g(t)}{\pi i}\overline{\int\limits_\Gamma \frac{\overline{\varphi(\tau)}\,dt}{\tau-t}}=f(t),\;\; t\in\Gamma,
\end{multline*}
if the contour  $\Gamma$ possesses corner points. Note that a complete theory of such equations is presented in \cite{DL85, DLS95}, whereas approximation methods have been studied in \cite{DiRS95, DiV96}.

Let $t_1,\dots,t_n\in\Gamma$ be the corner points of a piecewise-smooth contour $\Gamma$, and let $\mathbb{L}_p(\Gamma,\rho)$ denote the weighted $\mathbb{L}_p$-space with a power weight $\rho(t):=\prod\limits_{j=1}^n|t-t_j|^{\gamma_j}$. Assume that the parameters $p$ and $\beta_j:=(1+\gamma_j)/p$ satisfy the conditions
\begin{eqnarray*}
    1<p<\infty, \;\; 0<\beta_j<1,\;\; j=1,\dots,n.
\end{eqnarray*}
If the coefficients of the above equations are piecewise-continuous matrix functions, one can construct a function $\mathcal{A}_{\vec\beta}(t,\xi)$, $t\in\Gamma$, $\xi\in\mathbb{R}$, $\vec\beta:=(\beta_1,\dots,\beta_n)$, called the symbol of the equation (of the related operator). It is possible to express various properties of the equation in terms of $\mathcal{A}_{\vec\beta}$:
\begin{itemize}
\vskip+0.15cm
\item The equation is Fredholm in $\mathbb{L}_p(\Gamma,\rho)$ if and only if its symbol is elliptic., i.e. iff $\inf_{(t,\xi)\in\Gamma\times\mathbb{R}}|\;\mathcal{A}_{\vec\beta}(t,\xi)|>0$;

\vskip+0.15cm
\item To an elliptic symbol $\mathcal{A}_{\vec\beta}(t,\xi)$ there corresponds an integer valued index  \linebreak     $\ind\mathcal{A}_{\vec\beta}(t,\xi)$, the winding number, which coincides with the Fredholm index of the corresponding operator modulo a constant multiplier.
\end{itemize}

For more detailed survey of the theory and various applications to the problems of elasticity we refer the reader to \cite{Du75a, Du75b, Du77, Du79, Du82, Du84a, Du84b, Du86, Sc85}.

Similar approach to boundary integral equations on curves with corner points based on Mellin transformation has been exploited by M.~Costabel and E.~Stephan \cite{Cos83, CS84}.

However, one of the main problems in boundary integral equations for elliptic partial differential equations is the absence of appropriate results for Mellin convolution operators in Bessel potential spaces, cf. \cite{Du82, Du84b, Du86} and recent publications on nano-photonics \cite{BCC12a, BCC12b, GB10}. Such results are needed to obtain an equivalent reformulation of boundary value problems into boundary integral equations in Bessel potential spaces. Nevertheless, numerous works on Mellin convolution equations seem to pay almost no attention to the mentioned problem.

The first arising problem is the boundedness results for Mellin convolution operators in Bessel potential spaces. The conditions on kernels known so far are very restrictive. The following boundedness result for the Mellin convolution operator can be proved
 %
\begin{proposition}\label{p1.3}
Let $1<p<\infty$ and let $m=1,2,\ldots$ be an integer. If a function $\mathcal{K}$ satisfies the condition
\begin{equation}\label{e1.3}
    \int\limits_0^1 t^{\frac1p-m-1}|\mathcal{K}(t)|\,dt+\int\limits_1^\infty t^{\frac1p-1}|\mathcal{K}(t)|\,dt<\infty,
\end{equation}
then the Mellin convolution operator (see \eqref{e0.1})
\begin{equation}\label{e1.3a}
    \A=\mathfrak{M}^0_{\cA_{1/p}}:\widetilde{\mathbb{H}}{}^s_p(\mathbb{R}^+)\longrightarrow \mathbb{H}^s_p(\mathbb{R}^+)
\end{equation}
with the symbol $($see \eqref{e0.5}$)$
\begin{equation}\label{e1.4}
    \cA_{1/p}(\xi):=c_0+c_1\coth\pi\Big(\frac ip+\xi\Big)+\int\limits_{0}^{\infty}t^{\frac1p-i\xi} \mathcal{K}(t)\,\frac{dt}t\,,
                \;\; \xi\in\mathbb{R},
\end{equation}
is bounded for any $0\leqslant s\leqslant m$.
\end{proposition}

Note that the condition
\begin{equation}\label{e1.3b}
    K_\beta:=\int\limits_0^\infty t^{\beta-1}|\mathcal{K}(t)|\,dt<\infty
\end{equation}
and the constraints \eqref{e1.1} ensure that the operator
\[  \mathfrak{M}^0_a:\mathbb{L}_p(\mathbb{R}^+,t^\gamma)\longrightarrow \mathbb{L}_p(\mathbb{R}^+,t^\gamma)     \]
is bounded and the norm of the Mellin convolution
\begin{equation}\label{e1.3c}
\mathfrak{M}^0_{a_\beta}\varphi(t):=\int\limits_{0}^\infty\mathcal{K}\Big(\frac{t}
    \tau\Big)\varphi(\tau)\,\frac{d\tau}{\tau}
\end{equation}
admits the estimate $\|\mathfrak{M}^0_{a_\beta}\|\leqslant K_\beta$.

The above-formulated result has very restricted application. For example, the operators
\begin{equation}\label{e1.5}
\begin{aligned}
    \N_\alpha\varphi(t) & :=\frac{\sin\alpha}{\pi}\int\limits_0^\infty \frac{t\,\varphi(\tau)\,d\tau}{t^2+\tau^2-2t\tau\cos\alpha}\,, \\
    \N^*_\alpha\varphi(t) & :=\frac{\sin\alpha}{\pi}\int\limits_0^\infty \frac{\tau\,\psi_j(\tau)\,d\tau}{t^2+\tau^2-2t\tau\cos\alpha}\,, \\
    \M_\alpha\varphi(t) & :=\frac1{2\pi} \int\limits_{\bR^+} \frac{[\tau\,\cos\alpha-t]\varphi(\tau)\,d\tau}{t^2+\tau^2-2t\,\tau\cos\alpha}\,,\;\;
                -\pi<\alpha<\pi,
\end{aligned}
\end{equation}
which we encounter in boundary integral equations for elliptic boundary value problems (see \cite{BDKT13}), as well as the operators
\begin{equation}\label{e1.18x}
    \N_{m,k}\varphi(t):=\frac{t^k}{\pi i} \int\limits_0^\infty\frac{\tau^{m-k}\varphi(\tau)\,d\tau}{(\tau+t)^{m+1}}\,, \;\; k=0,\dots,m,
\end{equation}
represented in \eqref{e0.3}, do not satisfy the conditions \eqref{e1.3}. In particular, $\N_\alpha$ satisfies condition \eqref{e1.3} only for $m=1$ and $\N_{m,k}$ only for $m=k$. Although, as we will see below in Theorem \ref{t1.5}, all operators $\N_\alpha$, $\N^*_\alpha$ and $\N_{m,k}$ are bounded in Bessel potential spaces in the setting \eqref{e1.2} for all $s\in\bR$.

In the present paper we introduce {\em admissible kernels}, which are meromorphic functions on the complex plane $\mathbb{C}$, vanishing at the infinity
\begin{equation}\label{e1.7}
\begin{gathered}
\mathcal{K}(t):=\sum_{j=0}^\ell \frac{d_j}{t-c_j}+\sum_{j=\ell+1}^\infty
    \frac{d_j}{(t-c_j)^{m_j}}\,, \;\; c_j\not=0, \;\; j=0,1,\dots, \\
c_0,\dots,c_\ell\in\mathbb{R},\;\; 0<\alpha_k:=|\arg c_k|\leqslant\pi,\;\;
    k=\ell+1,\ell+2,\dots.
\end{gathered}
\end{equation}
$\mathcal{K}(t)$ have poles at $c_0,c_1,\ldots\in\mathbb{C}\setminus\{0\}$ and complex coefficients $d_j\in\mathbb{C}$. The Mellin convolution operator
\begin{equation}\label{e0.11}
    \mathbf{K}^m_c\varphi(t):=\displaystyle\frac1\pi\int\limits_0^\infty \frac{\tau^{m-1}\varphi(\tau)\,d\tau}{(t-c\,\tau)^m}.
 \end{equation}
with the kernel
 \[
\mathcal{K}^m_c(t):=\frac1{(t-c)^{m}},\qquad c_j\not=0
 \]
(see Definition \ref{d1.4}) turns out to be bounded in the Bessel potential spaces (see Theorem \ref{t1.5}).

In order to study Mellin convolution operators in Bessel potential spaces, we use the ``lifting'' procedure, performed with the help of the Bessel potential operators $\mathbf{\Lambda}^s_+$ and $\mathbf{\Lambda}^{s-r}_-$, which transform the initial operator $\mathfrak{M}_a^0$ into the lifted operator $\mathbf{\Lambda}^{s-r}_-\mathfrak{M}_a^0\mathbf{\Lambda}^{-s}_+$ acting already on a Lebesgue $\mathbb{L}_p$ spaces. However, the lifted operator is neither Mellin nor Fourier convolution and to describe its properties,  one has to study the commutants of Bessel potential operators and Mellin convolutions with meromorphic kernels. It turns out that Bessel potentials alter after commutation with Mellin convolutions and the result depends essentially on poles of the meromorphic kernels. These results allows us to show that the lifted operator $\mathbf{\Lambda}^{s-r}_-\mathfrak{M}_a\mathbf{\Lambda}^{-s}_+$ belongs to the Banach algebra of operators generated by Mellin and Fourier convolution operators with discontinuous symbols. Since such algebras have been studied before \cite{Du87}, one can derive various information (Fredholm properties, index, the unique solvability) about the initial Mellin convolution  equation $\mathfrak{M}_a^0\varphi=g$ in Bessel potential spaces in the settings $\varphi\in\widetilde{\mathbb{H}}^s_p(\mathbb{R}^+)$,  $g\in\widetilde{\mathbb{H}}^{s-r}_p(\mathbb{R}^+)$  and in the settings $\varphi\in\widetilde{\mathbb{H}}^s_p(\mathbb{R}^+)$, $g\in\mathbb{H}^{s-r}_p(\mathbb{R}^+)$.

The results of the present work is already applied in \cite{DTT14} to the investigation of some boundary value problems studied before by Lax--Milgram Lemma in \cite{BCC12a, BCC12b}. Note that the present approach is more flexible and provides better tools for analyzing the solvability of the boundary value problems and the asymptotic behavior of their solutions.

It is worth noting that the obtained results can also be used to study Schr\"{o}dinger operator on combinatorial and quantum graphs. Such a problem has attracted a lot of attention recently, since the operator mentioned above possesses interesting properties and has various applications, in particular, in nano-structures (see \cite{Ku04, Ku05} and
the references there). Another area for application of the present results are Mellin pseudodifferential operators on graphs. This problem has been studied in \cite{RR12}, but in the periodic case only. Moreover, some of the results can be applied in the study of stability of approximation methods for Mellin convolution equations in Bessel potential spaces.

The present paper is organized as follows. In the first section we observe Mellin and Fourier convolution operators with discontinuous symbols acting on Lebesgue spaces. Most of these results are well known and we recall them for convenience. In the second section we define Mellin convolutions with admissible meromorphic kernels and prove their boundedness in Bessel potential spaces. In Section~\ref{s2} is proved the key result on commutants of the Mellin convolution operator (with admissible meromorphic kernel) and a Bessel potential. In Section~\ref{s3} we enhance results on Banach algebra generated by Mellin and Fourier convolution operators by adding explicit definition of the symbol of a Hankel operator, which belong to this algebra. In Sections \ref{s4} the obtained results are applied to describe Fredholm properties and the index of Mellin convolution operators with admissible meromorphic kernels in Bessel potential spaces.

\section{Mellin Convolution and Bessel Potential Operators}
\label{s1}

Let $N$ be a positive integer. If there arises no confusion, we write $\mathfrak{A}$ for both scalar and matrix $N\times N$ algebras with entries from $\mathfrak{A}$. Similarly, the same notation $\mathfrak{B}$ is used for the set of $N$-dimensional vectors with entries from $\mathfrak{B}$. It will be usually clear from the context what kind of space or algebra is considered.

The integral operator  \eqref{e0.1} is called Mellin convolution. More generally, if $a\in\mathbb{L}_\infty(\mathbb{R})$ is an essentially bounded measurable $N\times N$ matrix function, the Mellin convolution operator $\mathfrak{M}_a^{0}$ is defined by
\begin{eqnarray}\label{e3.14a}
\mathfrak{M}_a^0\varphi(t):=\mathcal{M}^{-1}_\beta a\mathcal{M}_\beta \varphi(t)=
    \frac1{2\pi}\int\limits_{-\infty}^{\infty} \!a(\xi)\int\limits_0^\infty\! \Big(\frac t\tau\Big)^{i\xi-\beta}\varphi(\tau)\,\frac{d\tau}\tau\,d\xi,\\
 \varphi\in\mathbb{S}(\mathbb{R}^+),\nonumber
\end{eqnarray}
where $\mathbb{S}(\mathbb{R}^+)$ is the Schwartz space of
fast decaying functions on $\mathbb{R}^+$, whereas
$\mathcal{M}_\beta$ and $\mathcal{M}^{-1}_\beta$ are the Mellin
transform and its inverse, i.e.
\begin{align*}
    \mathcal{M}_\beta\psi(\xi) & :=\int\limits_0^\infty t^{\beta-i\xi}\psi(t)\,\frac{dt}t,\;\; \xi\in\mathbb{R},\\
    \mathcal{M}^{-1}_\beta\varphi(t) & :=\frac1{2\pi}\int\limits_{-\infty}^{\infty}t^{i\xi-\beta} \varphi(\xi)\,d\xi,\;\;t\in\mathbb{R}^+.
\end{align*}
The function $a(\xi)$ is usually referred to as a symbol of the Mellin operator $\mathfrak{M}_a^{0}$. Further, if the corresponding Mellin convolution operator $\mathfrak{M}_a^0$ is bounded on the weighted Lebesgue space $\mathbb{L}_p(\mathbb{R}^+,t^\gamma)$ of $N$-vector functions endowed with the norm
 \[
\big\|\varphi\mid\mathbb{L}_p(\mathbb{R}^+,t^\gamma)\big\|:=\bigg[\int\limits_0^\infty t^\gamma|\varphi(t)|^p\,dt\bigg]^{1/p},
 \]
then the symbol $a(\xi)$ is called a Mellin {\em $\mathbb{L}_{p,\gamma}$  multiplier}.

The transformations
\begin{align*}
&{\bf Z}_\beta:\mathbb{L}_p(\mathbb{R}^+,t^\gamma)\longrightarrow\mathbb{L}_p(\mathbb{R}),
     \quad{\bf Z}_\beta\varphi(\xi):=e^{-\beta t}\varphi(e^{-\xi}),\;\;\xi\in\mathbb{R},\\
&{\bf Z}^{-1}_\beta:\mathbb{L}_p(\mathbb{R})\longrightarrow\mathbb{L}_p(\mathbb{R}^+,
     t^\gamma),\quad{\bf Z}^{-1}_\beta\psi(t):=t^{-\beta}\psi(-\ln\,t),\;\; t\in\mathbb{R}^+,
\end{align*}
arrange an isometrical isomorphism between the corresponding $\mathbb{L}_p$-spaces. Moreover, the relations
\begin{equation}\label{e3.17}
\begin{gathered}
    \mathcal{M}_\beta=\mathcal{F} {\bf Z}_\beta,\quad \mathcal{M}_\beta^{-1}={\bf Z}^{-1}_\beta\mathcal{F}^{-1}, \\
    \mathfrak{M}_a^{0}=\mathcal{M}_\beta^{-1} a\mathcal{M}_\beta=
        {\bf Z}^{-1}_\beta\mathcal{F}^{-1} a\mathcal{F} {\bf Z}_\beta={\bf Z}^{-1}_\beta W^0_a {\bf Z}_\beta,\\
    -1<\gamma<p-1,\qquad \beta:=\frac{1+\gamma}p.\qquad 0,\beta<1,
\end{gathered}
\end{equation}
where $\mathcal{F}$ and  $\mathcal{F}^{-1}$ are the Fourier
transform and its inverse,
$$  \mathcal{F}\varphi(\xi):=\int\limits_{-\infty}^\infty e^{i\xi x} \varphi(x)\,dx,\quad
    \mathcal{F}^{-1}\psi(x):=\frac1{2\pi}\int\limits_{-\infty}^\infty e^{-i\xi x}\psi(\xi)\,d\xi,\;\;x\in\mathbb{R}\,,      $$
show a close connection between Mellin $\mathfrak{M}_a^{0}$
and Fourier
$$
W^0_a\varphi:=\mathcal{F}^{-1}a\mathcal{F}\varphi, \;\; \varphi\in\mathbb{S}(\mathbb{R}),      $$
convolution operators, as well as between the corresponding
transforms. Here $\mathbb{S}(\mathbb{R})$ denotes the Schwartz
class of infinitely smooth functions, decaying fast at the
infinity.

An $N\times N$ matrix function $a(\xi)$, $\xi\in\mathbb{R}$
is called a {\em Fourier $\mathbb{L}_p$-multiplier} if the
operator $W^0_a:\mathbb{L}_p(\mathbb{R}) \longrightarrow\mathbb{L}_p(\mathbb{R})$ is bounded. The set of all $\mathbb{L}_p$-multipliers is denoted by $\mathfrak{M}_p(\mathbb{R})$.

From \eqref{e3.17} immediately follows  the following
 %
\begin{proposition}[See \cite{Du79}]\label{p1}
Let $1<p<\infty$. The class of Mellin $\mathbb{L}_{p,\gamma}$-multipli\-ers coincides with the Banach algebra $\mathfrak{M}_p(\mathbb{R})$ of Fourier $\mathbb{L}_p$ -multipliers for arbitrary $-1<\gamma<p-1$ and is independent of the parameter $\gamma$.
\end{proposition}

Thus, a Mellin convolution operator $\mathfrak{M}_a^0$ in \eqref{e3.14a} is bounded in the weighted Lebesgue space $\mathbb{L}_p(\mathbb{R}^+,t^\gamma)$ if and only if $a\in\mathfrak{M}_p(\mathbb{R})$.

It is known, see, e.g. \cite{Du79}, that the Banach algebra $\mathfrak{M}_p(\mathbb{R})$ contains the algebra $\mathbf{V}_1(\mathbb{R})$ of all functions with bounded variation provided that
\begin{equation}\label{e1.1}
    \beta:=\frac{1+\gamma}p, \;\; 1<p<\infty, \;\;  -1<\gamma<p-1.
\end{equation}

As it was already mentioned, the primary aim of the present paper is to study Mellin convolution operators $\mathfrak{M}^0_a$ acting in Bessel potential spaces,
\begin{equation}\label{e1.2}
\mathfrak{M}^0_a:\widetilde{\mathbb{H}}^s_p(\mathbb{R}^+)\longrightarrow
    \mathbb{H}^s_p(\mathbb{R}^+).
\end{equation}
The symbols of these operators are $N\times N$ matrix functions
$a\in C\mathfrak{M}^0_p(\overline{\mathbb{R}})$, continuous
on the real axis $\mathbb{R}$ with the only one possible jump
at infinity. We commence with the definition of the Besseel
potential spaces and Bessel potentials, arranging isometrical
isomorphisms between these spaces and enabling the lifting
procedure, writing a Fredholm equivalent operator (equation)
in the Lebesgue space $\mathbb{L}_p(\mathbb{R}^+)$ for the
operator $\mathfrak{M}^0_a$ in \eqref{e1.2}.

For $s\in\mathbb{R}$ and $1<p<\infty$, the Bessel potential space,
known also as a fractional Sobolev space, is the
subspace of the Schwartz space $\mathbb{S}'(\mathbb{R})$ of
distributions having the finite norm
$$  \big\|\varphi\mid\mathbb{H}^s_p(\mathbb{R})\big\|:=\bigg[\int\limits_{-\infty}^\infty
            \big|\mathcal{F}^{-1}\big(1+|\xi|^2\big)^{s/2}(\mathcal{F}\varphi)(t)\big|^p\,dt\bigg]^{1/p}<\infty.        $$

For an integer parameter $s=m=1,2,\ldots\;$, the space
$\mathbb{H}^s_p(\mathbb{R})$ coincides with the usual Sobolev
space endowed with an equivalent norm
\begin{equation}\label{e1.23}
    \big\|\varphi\mid\mathbb{W}^m_p(\mathbb{R})\big\|:=\bigg[\sum_{k=0}^m\int\limits_{-\infty}^\infty
                \Big|\frac{d^k\varphi(t)}{dt^k}\Big|^p\,dt\bigg]^{1/p}.
 \end{equation}

If $s<0$, one gets the space of distributions. Moreover,
$\mathbb{H}^{-s}_{p'}(\mathbb{R})$ is the dual to the space
$\mathbb{H}^s_p(\mathbb{R}^+)$, provided $p':=\frac p{p-1}$,
$1<p<\infty$. Note that $\mathbb{H}^s_2(\mathbb{R})$ is a
Hilbert space with the inner product
 \[
\langle\varphi,\psi\rangle_s =\int\limits_{\mathbb{R}} (\mathcal{F}\varphi)(\xi) \overline{(\mathcal{F}\psi)(\xi)}(1+\xi^2)^s\,
     d\xi, \;\;\varphi, \psi \in \mathbb{H}^s(\mathbb{R}).
 \]
By $r_\Sigma$ we denote the operator restricting functions or
distributions defined  on $\mathbb{R}$ to the subset
$\Sigma\subset\mathbb{R}$. Thus
$\mathbb{H}^s_p(\mathbb{R}^+)=r_+(\mathbb{H}^s_p(\mathbb{R}))$,
and the norm in $\mathbb{H}^s_p(\mathbb{R}^+)$ is defined by
\[  \big\|f\mid\mathbb{H}^s_p(\mathbb{R}^+)\big\|=\inf_\ell \big\|\ell f\mid\mathbb{H}^s_p(\mathbb{R})\big\|,      \]
where $\ell f$ stands for any extension of $f$ to a distribution
in $\mathbb{H}^s_p(\mathbb{R})$.

Further, we denote by $\widetilde{\mathbb{H}}^s_p(\mathbb{R}^+)$
the (closed) subspace of $\mathbb{H}^s_p(\mathbb{R})$ which
consists of all distributions supported in the closure of
$\mathbb{R}^+$.

Notice that $\widetilde{\mathbb{H}}^s_p(\mathbb{R}^+)$ is
always continuously embedded in $\mathbb{H}^s_p(\mathbb{R}^+)$,
and if $s\in(1/p-1,1/p)$, these two spaces coincide. Moreover,
$\mathbb{H}^s_p(\mathbb{R}^+)$ may be viewed as the
quotient-space $\mathbb{H}^s_p(\mathbb{R}^+):=\mathbb{H}^s_p(\mathbb{R})/\widetilde{\mathbb{H}}^s_p(\mathbb{R}^-)$, $\mathbb{R}^-:=(-\infty,0)$.

Let $a\in\mathbb{L}_{\infty,loc}(\mathbb{R})$ be a locally bounded $m\times m$ matrix function. The Fourier convolution
operator (FCO) with the symbol $a$ is defined by
\begin{equation}\label{e1.46a}
    W^0_a:=\mathcal{F}^{-1}a\mathcal{F}.
\end{equation}
If the operator
\begin{equation}\label{e1.46}
    W^0_a:\mathbb{H}^s_p(\mathbb{R})\longrightarrow \mathbb{H}^{s-r}_p(\mathbb{R})
\end{equation}
is bounded, we say that $a$ is an $\mathbb{L}_p$-multiplier of order $r$ and use "$\mathbb{L}_p$-multiplier" if the order is $0$. The set of all $\mathbb{L}_p$-multipliers of order $r$ (of order $0$) is denoted by $\mathfrak{M}^r_p(\mathbb{R})$ (by $\mathfrak{M}_p(\mathbb{R})$, respectively).

For an $\mathbb{L}_p$-multiplier of order $r$, $a\in\mathfrak{M}^r_p(\mathbb{R})$, the Fourier convolution operator (FCO) on the semi-axis $\mathbb{R}^+$ is defined by the equality
\begin{equation}\label{e1.6}
    W_a=r_+W^0_a:\widetilde{\mathbb{H}}^s_p(\mathbb{R}^+)\longrightarrow \mathbb{H}^{s-r}_p(\mathbb{R}^+)
\end{equation}
and the Hankel operator by the equality
\begin{equation}\label{e1.47}
H_a=r_+\V W^0_a:\widetilde{\mathbb{H}}^s_p(\mathbb{R}^+)\longrightarrow
    \mathbb{H}^{s-r}_p(\mathbb{R}^+),\;\;\V\psi(t):=\psi(-t),
\end{equation}
where $r_+:=r_{\mathbb{R}^+}:\mathbb{H}^s_p(\mathbb{R})\longrightarrow \mathbb{H}^s_p(\mathbb{R}^+)$ is the restriction operator to the semi-axes
$\mathbb{R}^+$. Note that the generalized H\"ormander's kernel
of a FCO $W_a$ depends on the difference of arguments
$\mathcal{K}_a(t-\tau)$, while the H\"ormander's kernel
of a Hankel operator $r_+\V W^0_a$ depends of the sum of
the arguments $\mathcal{K}_a(t+\tau)$.

We did not use in the definition of the class of multipliers $\mathfrak{M}^r_p(\mathbb{R})$ the parameter $s\in\mathbb{R}$. This is due to the fact that $\mathfrak{M}^r_p(\mathbb{R})$ is independent of $s$: if the operator $W_a$ in \eqref{e1.47} is bounded for some $s\in\mathbb{R}$, it is bounded for all other values of $s$. Another definition of the multiplier class $\mathfrak{M}^r_p(\mathbb{R})$ is written as follows: $a\in\mathfrak{M}^r_p(\mathbb{R})$ if and only if $\lambda^{-r}a\in\mathfrak{M}_p( \overline{\mathbb{R}}) =\mathfrak{M}^0_p(\overline{\mathbb{R}})$, where $\lambda^r(\xi) :=(1+|\xi|^2)^{r/2}$. This assertion is one of the consequences of the following theorem.
 %
\begin{theorem}\label{t1.1}
Let  $1<p<\infty$. Then:
\begin{enumerate}
\vskip+0.15cm
\item For any $r,s\in\mathbb{R}$, $\gamma\in\mathbb{C}$, $\IIm\gamma>0$ the convolution operators $(\Psi$DOs$)$
\begin{eqnarray}\label{e1.24}
\begin{gathered}
\begin{aligned}
\mathbf{\Lambda}_\gamma^r=W_{\lambda^r_\gamma}:\widetilde{\mathbb{H}}^s_p
     (\mathbb{R}^+) & \longrightarrow\widetilde{\mathbb{H}}^{s-r}_p (\mathbb{R}^+),\\
\mathbf{\Lambda}_{-\gamma}^r=r_+W^0_{\lambda^r_{-\gamma}}\ell:\mathbb{H}^s_p
     (\mathbb{R}^+) & \longrightarrow\mathbb{H}^{s-r}_p(\mathbb{R}^+),
\end{aligned} \\
    \lambda^r_{\pm\gamma}(\xi):=(\xi\pm\gamma)^r,\;\;\xi\in\mathbb{R}, \;\; \IIm\gamma>0,
\end{gathered}
\end{eqnarray}
where $\ell:\mathbb{H}^s_p(\mathbb{R}^+)\longrightarrow\mathbb{H}^s_p (\mathbb{R})$ is an extension operator and $r_+$ is the restriction from the axes $\mathbb{R}$ to the semi-axes $\mathbb{R}^+$, arrange isomorphisms of the corresponding spaces. The final result is independent of the choice of an extension $\ell$.

\vskip+0.15cm
\item For any  operator  $\mathbf{A}:\widetilde{\mathbb{H}}^s_p(\mathbb{R}^+) \longrightarrow \mathbb{H}^{s-r}_p(\mathbb{R}^+)$ of order $r$, the following diagram is commutative
\begin{equation}\label{e1.25}
    \begin{matrix}
    \xymatrix{ \widetilde{\mathbb{H}}^s_p(\mathbb{R}^+) \ar[r]^{\mathbf{A}} &\mathbb{H}^{s-r}_p(\mathbb{R}^+) \ar[d]^{\mathbf{\Lambda}_{-\gamma}^{s-r}} \\
    \mathbb{L}_p(\mathbb{R}^+)\ar[r]_{\mathbf{\Lambda}_{-\gamma}^{s-r}\mathbf{A}
    \mathbf{\Lambda}^{-s}_\gamma} \ar[u]^{\mathbf{\Lambda}^{-s}_\gamma} &
    \mathbb{L}_p(\mathbb{R}^+)}\end{matrix}.
\end{equation}
The diagram \eqref{e1.24} provides an equivalent lifting of the
operator $\mathbf{A}$ of order $r$ to the operator
$\mathbf{\Lambda}_{-\gamma}^{s-r}\mathbf{A}\mathbf{\Lambda}^{-s}_\gamma:\mathbb{L}_p(
\mathbb{R}^+)\longrightarrow\mathbb{L}_p(\mathbb{R}^+)$ of order~$0$.
\vskip+0.15cm
\item
For any bounded convolution operator $W_a:\wt{\mathbb{H}}^s_p(\mathbb{R}^+) \longrightarrow \mathbb{H}^{s-r}_p(\mathbb{R}^+)$ of order $r$ and for any pair of complex numbers $\gamma_1, \gamma_2$ such that ${\rm Im}\,\gamma_j>0$, $j=1,2$, the lifted operator
\begin{equation}\label{e1.25a}
 \begin{array}{c}
\mathbf{\Lambda}_{-\gamma_1}^\mu W_a\mathbf{\Lambda}_{\gamma_2}^\nu
     =W_{a_{\mu,\nu}}\;:\;\wt{\mathbb{H}}^{s+\nu}_p(\mathbb{R}^+) \longrightarrow\mathbb{H}^{s-r-\mu}_p(\mathbb{R}^+),\\[2mm]
a_{\mu,\nu}(\xi):=(\xi-\gamma_1)^\mu a(\xi)(\xi+\gamma_2)^\nu
 \end{array}
\end{equation}
is again a Fourier convolution.

In particular,  the lifted operator  $W_{a_0}$ in $\mathbb{L}_p$-spaces, $\mathbf{ \Lambda}_{-\gamma}^{s-r}W_a\mathbf{\Lambda}^{-s}_\gamma:\mathbb{L}_p(\mathbb{R}^+)
\longrightarrow\mathbb{L}_p(\mathbb{R}^+)$ has the symbol
\begin{equation*}
a_{s-r,-s}(\xi)=\lambda^{s-r}_{-\gamma}(\xi)a(\xi)\lambda^{-s}_\gamma(\xi)
    =\Big(\frac{\xi-\gamma}{\xi+\gamma}\Big)^{s-r}\,\frac{a(\xi)}{(\xi+i)^r}\,.
\end{equation*}
\item
The Hilbert transform $S_\mathbb{\mathbb{R}^+}=i\mathbf{K}^1_1=W_{{\rm -\,sign}}$ is a Fourier convolution operator and
\begin{gather}\label{e1.10}
\mathbf{\Lambda}^s_{-\gamma_1}\boldsymbol{K}^1_1\mathbf{\Lambda}^{-s}_{\gamma_2}
     =W_{i\,g^s_{-\gamma_1,\gamma_2}\,{\rm sign}},
\end{gather}
where
\begin{eqnarray}\label{e1.11}
g^s_{-\gamma_1,\gamma_2}(\xi):=\left(\frac{\xi-\gamma_1}{\xi+\gamma_2}\right)^s.
\end{eqnarray}
 \end{enumerate}
 \end{theorem}
{\bf Proof:}  For the proof of items 1-3 we refer the reader  to
\cite[Lemma 5.1]{Du79} and \cite{DS93,Es81}. The item 4 is a consequence of the proved items 2 and 3 (see \cite{Du79,Du13}).
\QED
 %
\begin{remark}\label{r1.2}
The class of Fourier convolution operators is a subclass of pseudodifferential operators ($\Psi$DOs). Moreover, for integer parameters $m=1,2,\ldots$ the Bessel potentials
$\mathbf{\Lambda}_\pm^m=W_{\lambda^m_{\pm\gamma}}$, which are Fourier convolutions of order $m$, are ordinary differential operators of the same order~$m$:
\begin{equation}\label{e1.26}
\mathbf{\Lambda}_{\pm\gamma}^m=W_{\lambda^m_{\pm\gamma}}=\Big(i\,
     \frac{d}{dt}\pm\gamma\Big)^m=\sum_{k=0}^m \binom{m}{k}i^k(\pm\gamma)^{m-k}\,\frac{d^k}{dt^k}\,.
\end{equation}
These potentials map both spaces (cf. \eqref{e1.24})
\begin{equation}\label{e1.27}
\begin{aligned}
    \mathbf{\Lambda}_{\pm\gamma}^m & :\widetilde{\mathbb{H}}^s_p(\mathbb{R}^+)\longrightarrow \widetilde{\mathbb{H}}^{s-r}_p(\mathbb{R}^+), \\
    & :\mathbb{H}^s_p(\mathbb{R}^+)\longrightarrow\mathbb{H}^{s-m}_p(\mathbb{R}^+),
\end{aligned}
\end{equation}
but the mappings are not isomorphisms because the inverses $\mathbf{\Lambda}_{\pm \gamma}^{-m}$ do not map both spaces, only those indicated in \eqref{e1.24}.
\end{remark}
 %
\begin{remark}\label{r1.4}
For any pair of multipliers $a\in\mathfrak{M}^r_p(\mathbb{R})$, $b\in\mathfrak{M}^s_p(\mathbb{R})$
the corresponding convolution operators on the half axes $W^0_a$ and $W_b^0$ have the property $W^0_aW^0_b=W^0_bW^0_a=W^0_{ab}$.

For the corresponding Wiener-Hopf operators on the half axes a similar equality
\begin{equation}\label{e28}
W_aW_b=W_{ab}
\end{equation}
holds if and only if either the function $a(\xi)$ has an analytic extension in the lower half plane, or the function $b(\xi)$ has an analytic extension in the upper half plane (see \cite{Du79}).

Note that, actually \eqref{e1.25a} is a consequence of \eqref{e28}.
\end{remark}

\section[Mellin Convolutions with Admissible Meromorphic Kernels]{Mellin Convolutions with Admissible\\ Meromorphic Kernels}
\label{s2}

Now we consider kernels $\mathcal{K}(t)$, exposed in \eqref{e1.7}, which are meromorphic functions on the complex plane $\mathbb{C}$, vanishing at infinity, having poles at $c_0,c_1,\ldots\in\mathbb{C}\setminus\{0\}$ and complex coefficients $d_j\in\mathbb{C}$.
 %
\begin{definition}\label{d1.4}
We call a kernel $\mathcal{K}(t)$ in \eqref{e1.7} admissible iff:
\begin{itemize}
\vskip+0.15cm
\item[(i)] $\mathcal{K}(t)$ has only a finite number of poles
$c_0,\dots,c_\ell$ which belong to the positive semi-axes,
i.e., $\arg c_0=\cdots=\arg c_\ell=0$;

\vskip+0.15cm
\item[(ii)] The corresponding multiplicities are one
$m_0=\cdots=m_\ell=1$;

\vskip+0.15cm
\item[(iii)] The remainder points $c_{\ell+1},c_{\ell+2},\ldots$ do not
condense to the positive semi-axes and their real parts are bounded uniformly
\begin{equation}\label{e1.8a}
\underset{j\longrightarrow\infty}{\underline{\lim}} c_j\not\in[0,\infty),\quad
        \sup\limits_{j=\ell+1,\ell+2,\ldots}\RRe c_j\leqslant K<\infty.
\end{equation}

\vskip+0.15cm
\item[(iv)]
$\mathcal{K}(t)$ is a kernel of an operator, which is a composition of finite number of operators with admissible kernels.
\end{itemize}
\end{definition}

\begin{example}\label{ex1.4a}
The function
\[  \mathcal{K}(t)=\exp\Big(\frac1{t-c}\Big),\;\; \RRe c<0\;\;\text{or}\;\;\IIm c\not=0        \]
is an example of the admissible kernel which also satisfies the condition of the next Theorem \ref{t1.5}. Other examples of operators with admissible kernels (which also satisfies the condition of the next Theorem \ref{t1.5}) are operators which we encounter in \eqref{e0.3}, in \eqref{e1.5} and in \eqref{e1.6} and, in general, any finite sum in \eqref{e1.7}.
\end{example}

\begin{example}\label{ex1.4b}
The function
\[
\mathcal{K}(t)=\frac{\ln\,t-c_1c_2}{t-c_1c_2}\,, \;\; \IIm c_1\not=0,\;\;
    \IIm c_2\not=0
\]
is another example of the admissible kernel and represnts the composition of operators $c_2\mathbf{K}^1_{c_1}\mathbf{K}^1_{c_2}$ (see \eqref{e1.11}) with admissible kernels which also satisfies the condition of the next Theorem \ref{t1.5}. More trivial examples of operators with admissible kernels (which also satisfies the condition of the next Theorem~\ref{t1.5}) are operators which we encounter in \eqref{e0.3}, in \eqref{e1.5} and in \eqref{e1.6} and, in general, any finite sum in \eqref{e1.7}.
\end{example}
 %
\begin{theorem}\label{t1.5}
Let conditions \eqref{e1.1} hold, $\mathcal{K}(t)$ in \eqref{e1.7} be an admissible kernel and
\begin{equation}\label{e1.8b}
    K_\beta:=\sum_{j=0}^\infty2^{m_j}|d_j|\,|c_j|^{\beta-m_j}<\infty.
\end{equation}
Then the Mellin convolution $\mathfrak{M}^0_{a_\beta}$ in \eqref{e1.3c} with the admissible meromorphic kernel $\mathcal{K}(t)$ in \eqref{e1.7} is bounded in the Lebesgue space $\mathbb{L}_p(\mathbb{R}^+,t^\gamma)$ and its norm is estimated by the constant $\|\mathfrak{M}^0_{a_\beta}\mid\cL(\mathbb{L}_p (\mathbb{R}^+,t^\gamma))\|\leqslant MK_\beta$ with some $M>0$.

We can drop the constant $M$ and replace $2^{m_j}$ by
$2^{\frac{m_j}2}$ in the estimate \eqref{e1.8b} provided
$\RRe c_j<0$ for all $j=0,1,\ldots\;$.
\end{theorem}

{\bf Proof:}
The first $\ell+1$ summands in the definition of the admissible
kernel \eqref{e1.7} correspond to the Cauchy operators
 \[
A_0\varphi(t)=\sum_{j=0}^\ell\frac{d_j}\pi\int\limits_0^\infty
     \frac{\varphi(\tau)\,d\tau}{t-c_j\tau},\;\;c_j>0, \;\; j=0,1,\dots,\ell,
 \]
and their boundedness property in the weighted Lebesgue space
\begin{equation}\label{ex.0}
    A_0:\mathbb{L}_p(\mathbb{R}^+,t^\gamma)\longrightarrow\mathbb{L}_p(\mathbb{R}^+,t^\gamma)
\end{equation}
under constraints \eqref{e1.1} is well known (see \cite{Kh57} and also /\cite{GK79}). Therefore we can ignore the first $\ell$ summands in the expansion of the kernel $\mathcal{K}(t)$ in \eqref{e1.7}. To the boundedness of the operator $\mathfrak{M}^0_{a^\ell_\beta}$ with the remainder kernel
3999\begin{gather*}
    \mathcal{K}^\ell(t):=\sum_{j=\ell+1}^\infty \frac{d_j}{(t-c_j)^{m_j}}\,,\;\; c_j\not=0, \;\; j=0,1,\dots, \\ 0<\alpha_k:=|\arg c_k|\leqslant\pi,\;\; k=\ell+1,\ell+2,\ldots
\end{gather*}
(see \eqref{e1.7}), we apply the estimate \eqref{e1.3b}
\begin{eqnarray}\label{ex.1}
    \big\|\mathfrak{M}^0_{a^\ell_\beta}\mid\cL(\mathbb{L}_p(\mathbb{R}^+,t^\gamma))\big\|
    \leqslant \int\limits_0^\infty t^{\beta-1}|\mathcal{K}^\ell(t)|\,dt\nonumber\\
    \leqslant \sum_{j=\ell+1}^\infty |d_j|
                \int\limits_0^\infty \frac{t^{\beta-1}dt}{|t -c_j|^{m_j}}\,.
\end{eqnarray}

Now note that
\begin{eqnarray*}
    |t-c_j|^{-m_j}=\big(t^2+|c_j|^2-2\RRe c_jt\big)^{-\frac{m_j}2}\leqslant
            \Big(\frac{t^2+|c_j|^2}2\Big)^{-\frac{m_j}2}\\
    \leqslant2^{m_j}(t+|c_j|)^{-m_j} \qquad\text{for all}\;\; t\geqslant2K=2\sup|\RRe c_j|>0.
\end{eqnarray*}
due to the constraints \eqref{e1.8a}. On the other hand,
\[  |t-c_j|^{-m_j}\leqslant M(t+|c_j|)^{-m_j} \;\;\text{for all}\;\; 0\leqslant t\leqslant 2K   \]
and a certain constant $M>0$. Therefore
\begin{equation}\label{ex.2}
    |t-c_j|^{-m_j}\leqslant M2^{m_j}(t+|c_j|)^{-m_j} \;\;\text{for all}\;\; 0<t<\infty.
\end{equation}

Next we recall the formula  from \cite[Formula 3.194.4]{GR94}
\allowdisplaybreaks[0]
\begin{gather}\label{e1.28}
    \int\limits_0^\infty\!\!\frac{t^{\beta-1}\,dt}{(t\!+\!c)^m}\!=\!(-1)^{m-1}\binom{\beta\!-\!1}{m\!-\!1}\,\frac{\pi c\,^{\beta-m}}{\sin\pi\beta},
                \;\; -\pi\!<\!\arg c\!<\!\pi,\;\;\RRe\beta\!<\!1, \\
    \binom{\beta-1}{m-1}:=\frac{(\beta-1)\cdots(\beta-m+1)}{(m-1)!}\,, \;\;
        \binom{\beta-1}{0}:=1 \nonumber
\end{gather}
to calculate the integrals. By inserting the estimate
\eqref{ex.2} into \eqref{ex.1} and applying \eqref{e1.28},
we get
\begin{eqnarray}\label{ex.3}
&&\big\|\mathfrak{M}^0_{a^\ell_\beta}\mid\cL(\mathbb{L}_p(\mathbb{R}^+,t^\gamma))\big\|
    \leqslant \sum_{j=\ell+1}^\infty |d_j|\int\limits_0^\infty \frac{t^{\beta-1}dt}{|t -c_j|^{m_j}}\nonumber\\
&&\hskip15mm\leqslant M_0\sum_{j=\ell+1}^\infty 2^{m_j}|d_j|\int\limits_0^\infty
    \frac{t^{\beta-1}dt}{(t+|c_j|)^{m_j}}\nonumber\\
&&\hskip15mm\leqslant \frac{\pi M_0}{\sin\pi\beta}\sum_{j=\ell+1}^\infty 2^{m_j}|d_j|
    \bigg|\binom{\beta-1}{m_j-1}\bigg|c_j^{\beta-m_j}\nonumber\\
&&\hskip15mm\leqslant M\sum_{j=\ell+1}^\infty2^{m_j}
    |d_j|c_j^{\beta-m_j}=MK_\beta,\quad M:=\frac{\pi M_0}{\sin\pi\beta},
\end{eqnarray}
since (see \eqref{e1.28})
\[  \bigg|\binom{\beta-1}{m_j-1}\bigg|\leqslant 1,      \]
where $K_\beta$ is from \eqref{e1.8b}. The boundedness
\eqref{ex.0} and the estimate \eqref{ex.3} imply the claimed
estimate
\[  \big\|\mathfrak{M}^0_{a_\beta}\mid \cL(\mathbb{L}_p(\mathbb{R}^+,t^\gamma))\big\|\leqslant MK_\beta.        \]

If $\RRe c_j<0$ for all $j=0,1,\ldots$, we have
\begin{multline*}
    \frac1{|t-c_j|^{m_j}}=\big(t^2+|c|^2-2\RRe c_jt\big)^{-\frac{m_j}2}
    \leqslant \big(t^2+|c|^2\big)^{-\frac{m_j}2}\leqslant 2^{\frac{m_j}2}\big(t+|c_j|\big)^{-m_j}
\end{multline*}
valid for all $t>0$ and a constant $M$ does not emerge in the
estimate.
\QED

Let us find the symbol (the Mellin transform of the kernel) of
the operator \eqref{e1.11} for $-\pi<\arg(-c)<\pi$, $m=1,2,\ldots$ (see \eqref{e1.10}, \eqref{e1.11}). For this we apply formula \eqref{e1.28}:
\begin{multline}\label{e1.29}
\cM_\bt\cK^{m}_{c}(\xi)=\int\limits_{0}^{\infty}t^{\beta-i\xi-1} \cK^{m}_{c}(t)\,dt=
   \displaystyle\frac1\pi\int\limits_0^\infty\frac{t^{\beta-i\xi-1}}{(t+(-c))^{m}}\,dt= \\
    =\binom{\beta-i\xi-1}{m-1}\,\frac{(-1)^{m-1}(-c)^{\beta-i\xi-m}}{
    \sin\pi(\beta-i\xi)}, \qquad -\pi<\arg(-c)<\pi,
\end{multline}
where, for $\delta\in\bC$, $(-c)^\delta:=|c|^\delta e^{-i\delta\arg\,c}$. In particular,
\begin{align}
\label{e3.11d}
\cM_\beta\cK^1_{-d}(\xi) & =\frac{d\,^{\beta-i\xi-1}}{\sin\pi(\beta-i\xi)}\,,
     \;\; -\pi<\arg d<\pi, \\
\label{e3.11e}
    \cM_\beta\cK^1_{-1}(\xi) & =\frac1{\sin\pi (\beta-i\xi)}\,, \;\;\xi\in\mathbb{R}.
\end{align}

Now let us find the symbol of the Cauchy singular integral
operator $K^1_1=-iS_{\mathbb{R}^+}$ (see \eqref{e1.12}, \eqref{e1.13}). For this we apply Plemeli formula and formula \eqref{e1.28}:
\allowdisplaybreaks
\begin{align}\label{e3.11f}
    \cM_\beta\cK^1_1(t) & :=\int\limits_0^\infty t^{\beta-i\xi-1}\cK^1_1(t)\,dt=-\displaystyle\frac1\pi\int\limits_0^\infty\frac{t^{\beta
    -i\xi-1}\,dt}{t-1} \nonumber \\
    &\,\, =\lim_{\ve\longrightarrow0} \frac1{2\pi}\int\limits_0^\infty \Big[\frac{t^{\beta-i\xi-1}}{t+e^{i(\pi-\ve)}}+\frac{t^{\beta-i\xi-1}}{t+e^{-i(\pi
    -\ve)}}\Big]\,dt \nonumber \\
    &\,\, =\lim_{\ve\longrightarrow0} \frac{e^{i(\pi-\ve)(\beta-i\xi-1)}+{e^{-i(\pi-\ve)(\beta-i\xi-1)}}}{2
    \sin\pi(\beta-i\xi)}= \cot\pi(\beta-i\xi).
\end{align}

For an admissible kernel with simple (non-multiple) poles
$m_0=m_1=\cdots=1$ and $\arg c_0=\arg c_\ell=0$ and $ 0<\pm\arg c_j<\pi$, $j=\ell+1,\ldots$ we get
\begin{multline}\label{e1.31}
    \cM_\bt\cK(\xi)=\cot\pi(\beta-i\xi)\sum_{j=0}^\ell d_jc_j^{\beta-i\xi-1}- \\
    -\frac1{\sin\pi(\beta-i\xi)}\sum_{j=\ell+1}^\infty d_j\binom{\beta-i\xi-1}{m-1} e^{\mp\pi(\beta-i\xi)i}c^{\beta-i\xi-m}.
\end{multline}
%
\begin{theorem}\label{t1.5a}
If $\mathcal{K}$ is an admissible kernel the corresponding Mellin convolution operator with the kernel $\mathcal{K}$
\begin{equation}\label{eqn21}
\begin{array}{c}
{\bf K}\varphi(t):=\displaystyle\int_0^\infty\mathcal{K}\left(\displaystyle\frac t\tau\right)\varphi(\tau)\displaystyle\frac{d\tau}\tau,\\[3mm]
{\bf K}\;:\;\widetilde{\mathbb{H}}^s_p(\mathbb{R}^+)\longrightarrow
\mathbb{H}^s_p(\mathbb{R}^+),
\end{array}
 \end{equation}
is bounded for all $1<p<\infty$ and $s\in\mathbb{R}$.

The condition on the parameter $p$ can be relaxed to $1\leqslant p\leqslant\infty$, provided the admissible kernel $\mathcal{K}$  in \eqref{e1.7} has no poles on positive semi-axes $\alpha_j=\arg c_j\not=0$ for all $j=0,1,\ldots\;$.
\end{theorem}
{\bf Proof:} Due to the representation \eqref{e1.7}, we have to prove the theorem only for a model kernel
\begin{equation}\label{e1.10a}
\mathcal{K}^m_c(t):=\frac1{\pi(t-c)^m}, \qquad c\not=0, \quad 0\leqslant|\arg c|<\pi,
     \quad m=1,2,\ldots\,.
\end{equation}
The respective Mellin convolution operator $\mathbf{K}^m_c$ (see \eqref{e1.11}) is bounded in $\mathbb{L}_p(\mathbb{R}^+)$ for all $1\leqslant p\leqslant\infty$ for arbitrary $0<|\arg c|<\pi$ (cf. \eqref{e0.2}).

To accomplish the boundedness result of $\mathbf{K}^m_c$ in $\mathbb{L}_p(\mathbb{R}^+)$ we need to consider the case $\arg c=0$ (i.e., $c\in(0,\infty)$) and, therefore, $m=1$ (see Definition \ref{d1.4}). Then the operator $\mathbf{K}^1_c$ coincides with the "dilated" Cauchy singular integral operator with a constant multiplier
\begin{equation}\label{e1.13}
\mathbf{K}^1_c\varphi(t):=\frac1\pi\int\limits_{0}^\infty\frac{\varphi(\tau)d\tau}{t-c\,\tau}
     =-\frac ic\,(S_{\mathbb{R}^+}\varphi)\Big(\frac{t}c\Big),
\end{equation}
where
\begin{equation}\label{e1.12}
    S_{\mathbb{R}^+}\varphi(t):=\frac1{\pi i}\int\limits_0^\infty \frac{\varphi(\tau)\,d\tau}{\tau - t},
\end{equation}
and is bounded in $\mathbb{L}_p(\mathbb{R}^+)$ for all $1<p<\infty$ (cf., e.g., \cite{Du79, GK79}).

Now let $0\leqslant\arg c<2\pi$ and $m=1$. Then, if $\vf\in C^\infty_0(\bR^+)$ is a smooth function with compact support and $k=1,2,\ldots$, integrating by parts we get
\begin{align}\label{e1.14}
    \frac{d^k}{dt^k}\,\mathbf{K}^1_{c}\varphi(t) & =\frac1\pi\int\limits_0^\infty \frac{d^k}{dt^k}\,\frac1{t-c\,\tau}\varphi(\tau)\,d\tau=
        \frac{(-c)^{-k}}\pi\int\limits_0^\infty \frac{d^k}{d\tau^k}\,\frac1{t-c\,\tau}\varphi(\tau)\,d\tau \nonumber \\
    & =\frac{c^{-k}}\pi\int\limits_0^\infty \frac1{t-c\,\tau}\,\frac{d^k\varphi(\tau)}{d\tau^k}\,d\tau=
            c^{-k}\Big(\mathbf{K}^1_c\,\frac{d^k}{dt^k}\varphi\Big)(t).
\end{align}
For $m=2,3,\ldots$ and $0<\arg c<2\pi$ we get similarly
\begin{align*}
    \frac{d}{dt}\,\mathbf{K}^m_c\varphi(t) & =\frac1\pi\int\limits_0^\infty \frac{d}{dt}\,\frac{\tau^{m-1}}{(t-c\,\tau)^m}\,\varphi(\tau)\,d\tau \\
    & =\sum_{j=0}^{m-1}\frac{(-c)^{-1-j}}\pi\int\limits_0^\infty \frac{d}{d\tau}\,\frac{\tau^{m-1-j}}{(t-c\,\tau)^{m-j}}\,\varphi(\tau)\,d\tau \\
    & =-\sum_{j=0}^{m-1}\frac{(-c)^{-1-j}}\pi\int\limits_0^\infty \frac{\tau^{m-1-j}}{(t-c\,\tau)^{m-j}}\,\frac{d}{d\tau}\,\varphi(\tau)\,d\tau \\
    & =-\sum_{j=0}^{m-1}(-c)^{-1-j}\Big(\mathbf{K}^{m-j}_c\,\frac{d}{dt}\,\varphi\Big)(t)
\end{align*}
and, recurrently,
\begin{gather}\label{e1.16}
\frac{d^k}{dt^k}\,\mathbf{K}^m_c\varphi(t)\!=\!(-1)^k\!\sum_{j=0}^{m-1}\!(-c)^{-k-j}
    \gamma^k_j\Big(\mathbf{K}^{m-j}_c\,\frac{d^k}{dt^k}\varphi\Big)(t),\;\;k\!=\!1,2,\dots,\\
\gamma^1_j=j+1,\;\;\gamma_0^k=1,\;\;\gamma^k_j:=\sum_{r=0}^j\gamma^{k-1}_r,
    \;\; j=0,1,\dots,m,\;\; k=1,2,\ldots\,. \nonumber
\end{gather}

Recall now that for an integer $s=n$ the spaces
$\mathbb{H}^n_p(\mathbb{R}^+)$, $\wt{\mathbb{H}}^n_p(\mathbb{R}^+)$  coincide with the Sobolev spaces $\mathbb{W}^n_p(\mathbb{R}^+)$, $\wt{\mathbb{W}}^n_p(\mathbb{R}^+)$, respectively (these spaces are isomorphic and the norms are equivalent) and
$C^\infty_0(\bR^+)$ is a dense subset in $\wt{\mathbb{W}}^n_p(\mathbb{R}^+) =\wt{\mathbb{H}}^n_p(\mathbb{R}^+)$. Then, using the equalities \eqref{e1.14}, \eqref{e1.16} and the boundedness of the operators $\mathbf{K}^{m-j}_c$ (see \eqref{e1.10a}--\eqref{e1.12}), we proceed as follows:
\begin{align}\label{e1.17}
\big\|\mathbf{K}^m_c\varphi\mid\mathbb{H}^n_p(\mathbb{R}^+)\big\|&=\sum_{k=0}^n
    \Big\|\frac{d^k}{dt^k}\,\mathbf{K}^m_c\varphi\mid\mathbb{L}_p(\mathbb{R}^+)\Big\|=
    \nonumber \\
&=\sum_{k=0}^n\sum_{j=0}^{m-1} |c|^{-k-j}\gamma^k_j\Big\|\mathbf{K}^{m-j}_c\, \frac
    {d^k}{dt^k}\,\varphi\mid\mathbb{L}_p (\mathbb{R}^+)\Big\|\leqslant \nonumber \\
&\leqslant M\sum_{k=0}^n \Big\|\frac{d^k}{dt^k}\,\varphi\mid \mathbb{L}_p(
    \mathbb{R}^+)\Big\|=M\big\|\varphi\mid\mathbb{H}^n_p(\mathbb{R}^+)\big\|,
\end{align}
where $M>0$ is a constant, and the boundedness \eqref{eqn21} follows for $s=0,1,2,\ldots\;$. The case of arbitrary $s>0$ follows by the interpolation between the spaces $\mathbb{H}^m_p(\mathbb{R}^+)$ and $\mathbb{H}^0_p(\mathbb{R}^+)=\mathbb{L}_p(\mathbb{R}^+)$, also between the spaces $\wt{\mathbb{H}}^m_p(\mathbb{R}^+)$ and $\wt{\mathbb{H}}^0_p(\mathbb{R}^+)=\mathbb{L}_p(\mathbb{R}^+)$.

For $s<0$ the boundedness \eqref{eqn21} follows by duality: the adjoint operator to $\mathbf{K}^m_c$ is
\begin{equation}\label{e1.18}
    \mathbf{K}^{m,*}_c\varphi(t):=\frac1\pi\int\limits_0^\infty \frac{t^{m-1}\varphi(\tau)\,d\tau}{(\tau-c\,t)^m}=\sum_{j=1}^m \omega_j\mathbf{K}^j_{c^{-1}}\varphi(t),
\end{equation}
for some constant coefficients $\omega_1,\dots,\omega_m$. The
operator $\mathbf{K}^{m,*}_c$ has the admissible kernel and, due to the
proved part of the theorem is bounded in the space setting
$\mathbf{K}^{m,*}_c:\wt{\mathbb{H}}^{-s}_{p'}(\mathbb{R}^+)\longrightarrow
\mathbb{H}^{-s}_{p'}(\mathbb{R}^+)$, $p':=p/(p-1)$, since $-s>0$. The initial operator $\mathbf{K}^m_c:\wt{\mathbb{H}}^s_p(\mathbb{R}^+)\longrightarrow\mathbb{H}^s_p(
\mathbb{R}^+)$ is dual to $\mathbf{K}^{m,*}_c$ and, therefore, is bounded as well. \QED

\begin{corollary}\label{c2.5}
Let $1<p<\infty$ and $s\in\bR$. A Mellin convolution operator $\mathfrak{M}^0_a$ with an admissible kernel described in
Definition $\ref{d1.4}$ $($also see Example $\ref{ex1.4b})$ and Theorem $\ref{t1.5}$ is bounded in Bessel potential spaces
\begin{equation*}
    \mathfrak{M}^0_a:\widetilde{\mathbb{H}}^s_p(\mathbb{R}^+)\longrightarrow \mathbb{H}^s_p(\mathbb{R}^+).
\end{equation*}
\end{corollary}

The bounddness property
\begin{equation*}
    \mathfrak{M}^0_a:\mathbb{H}^s_p(\mathbb{R}^+)\longrightarrow \mathbb{H}^s_p(\mathbb{R}^+).
\end{equation*}
does not hold in general for even a simplest Mellin convolution operator $\K_c$,  except the case when the spaces $\widetilde{\mathbb{H}}^s_p(\mathbb{R}^+$ and $\mathbb{H}^s_p(\mathbb{R}^+)$ can be identified, i.e., except the case
$1/p-1<s<1/p$. Indeed, to check this consider a smooth function with a compact support $\vf\in C^\infty_0(\bR^+$ which is constant on the unit interval: $\vf(t)=1$ for $0<t<1$. Obviously, $\vf\in\mathbb{H}^s_p(\mathbb{R}^+)$ and
$\vf\not\in\widetilde{\mathbb{H}}^s_p(\mathbb{R}^+$ for all $s>1/p$. Then,
 \[
\K_c\vf(t)=\displaystyle\frac1\pi\int\limits_0^\infty \frac{\varphi(\tau)\,d\tau}{t-c\,\tau}
     =\displaystyle\frac1\pi\int\limits_0^1 \frac{d\tau}{t-c\,\tau}+\displaystyle\frac1\pi\int\limits_1^\infty \frac{\varphi(\tau)\,d\tau}{t-c\,\tau}=c^{-1}\ln\tau + \vf_0(t),
 \]
where $\vf_0\in\bH^s_p(\bR^+)\cap C^\infty(\bR^+)$, while the first summand $\ln\,\tau$ does not belong to $\bH^s(\bR^+)$  since all functions in this space are continuous and uniformly bounded for $s>1/p$.

We can prove the following very partial result, which has important practical applications.
%
\begin{theorem}\label{t1.5b}
Let $1<p<\infty$ and $c\in\bC$. Then if $\displaystyle\frac1p-1<s<\displaystyle\frac1p+1$, the operator
\begin{equation}\label{e2.52}
\A_c:=c\K_c-c^{-1}\K_{c^{-1}}\;:\;\mathbb{H}^s_p(\mathbb{R}^+)\longrightarrow\mathbb{H}^s_p(\mathbb{R}^+),
 \end{equation}
is bounded while for $\displaystyle\frac1p-2<s<\displaystyle\frac1p$ the operator
\begin{equation}\label{e2.53}
\A^\#_c:=\K_c - \K_{c^{-1}}\;:\;\mathbb{H}^s_p(\mathbb{R}^+)\longrightarrow\mathbb{H}^s_p(\mathbb{R}^+),
 \end{equation}
is bounded.
\end{theorem}
{\bf Proof:} If $\displaystyle\frac1p-1<s<\displaystyle\frac1p$ the  spaces $\widetilde{\mathbb{H}}^s_p(\mathbb{R}^+$ and $\mathbb{H}^s_p(\mathbb{R}^+)$ can be identified and the boundedness \eqref{e2.52}, \eqref{e2.53}  follow from  Theorem \ref{t1.5a}.

Due to \eqref{e1.27} the following diagram
\begin{equation}\label{e2.54}
    \begin{matrix}
    \xymatrix{\mathbb{H}^s_p(\mathbb{R}^+) \ar[r]^{\A_\alpha} &\mathbb{H}^s_p(\mathbb{R}^+) \ar[d]^{\mathbf{\Lambda}_1^{-1}}\\
    \ar[r]_{\mathbf{\Lambda}_1^1\A_\alpha\mathbf{\Lambda}_1^{-1}}
    \mathbb{H}^{s-1}_p(\mathbb{R}^+) \ar[u]^{\mathbf{\Lambda}^{-1}_1} &
    \mathbb{H}^{s-1}_p(\mathbb{R}^+)}\end{matrix}.
\end{equation}
is commutative. The diagram \eqref{e2.54} provides an equivalent lifting of the operator $\A_\alpha$ from the space $\mathbb{H}^s_p(\mathbb{R}^+) $ to the operator $\A^0_\alpha:=\mathbf{\Lambda}_1^1\A_\alpha\mathbf{\Lambda}_1^{-1}$ in the space $\mathbb{H}^{s-1}_p(\mathbb{R}^+) $. On the other hand, $\mathbf{\Lambda}^1_1=i\partial_t+I$ (see \eqref{e1.26}) and it can be checked easily, using the integration by parts, that $\partial_t\A_\alpha=-\A^\#_\alpha\partial_t$.   Then,
\begin{gather*}
\A^0_\alpha=\mathbf{\Lambda}_1^1\A_\alpha\mathbf{\Lambda}_1^{-1}
     =(i\partial_t+I)\A_\alpha\mathbf{\Lambda}_1^{-1}=(\A_\alpha - \A^\#_\alpha)\mathbf{\Lambda}_1^{-1}
     +\A^\#_\alpha(i\partial_t+I)\mathbf{\Lambda}_1^{-1}\\[3mm]
=(\A_\alpha - \A^\#_\alpha)\mathbf{\Lambda}_1^{-1} + \A^\#_\alpha
\end{gather*}
Since $\displaystyle\frac1p-1<s-1<\displaystyle\frac1p$ and the embedding
 \[
\mathbf{\Lambda}_1^{-1} \mathbb{H}^{s-1}_p(\mathbb{R}^+)= \mathbb{H}^s_p(\mathbb{R}^+)
      \subset \mathbb{H}^{s-1}_p(\mathbb{R}^+)
 \]
is continuous, the operator
 \[
(\A_\alpha - \A^\#_\alpha)\mathbf{\Lambda}_1^{-1} + \A^\#_\alpha\;:\;
      \mathbb{H}^{s-1}_p(\mathbb{R}^+) \longrightarrow \mathbb{H}^{s-1}_p(\mathbb{R}^+)
\]
is bounded and the boundedness result in \eqref{e2.52} follows.

Now let $\displaystyle\frac1p-2<s<\displaystyle\frac1p$. Then
\begin{equation}\label{e2.55}
\frac1{p'}-1=-\frac1p<-s<\frac1{p'}+1=2-1/p, \qquad p':=\frac p{p-1}.
 \end{equation}

The pair of the operator $\K_c$ and $-\overline{c}{}^{-1}\K_{\overline{c}{}^{-1}}$ are adjoint to each other. Therefore,  the operator
\begin{equation}\label{e2.56}
 \A_{\overline{c}}:=\overline{c}\K_{\overline{c}}-\overline{c}{}^{-1}\K_{\overline{c}{}^{-1}}
      \;:\;\mathbb{H}^{-s}_{p'}(\mathbb{R}^+)\longrightarrow\mathbb{H}^{-s}_{p'}(\mathbb{R}^+),
 \end{equation}
is the adjoint to the operator  $\A^\#_c$ in \eqref{e2.53}. Since  the parameters $\{-s,p'\}$ satisfy the condition of the first part of the present theorem (see \eqref{e2.55}), the operator  $ \A_{\overline{c}}$ in \eqref{e2.56} is bounded and jaustifies the boundedness of the adjoint operator $\A^\#_c$ in \eqref{e2.53}.         \QED

The next result is crucial in the present investigation. Note that, the case $\arg\,c=0$ is essentially different and will be considered in Theorem \ref{t3.4} below.
 %
\begin{theorem}\label{t1.6}
Let $0<\arg\,c<2\pi$, $\arg\,c\not=0$ and $0<\arg(-c\,\gamma)<\pi$. Then
\begin{gather}\label{e1.38}   
\mathbf{\Lambda}^s_{-\gamma}\mathbf{K}^1_c\varphi=c^{-s}\mathbf{K}^1_c
     \mathbf{\Lambda}^s_{-c\,\gamma}\varphi,\qquad \varphi\in\widetilde{\mathbb{H}}^r_p(\mathbb{R}^+),
\end{gather}
where $c^{-s}=|c|^{-s}e^{-is\arg\,c}$.
\end{theorem}
{\bf Proof:}  First of all note, that due to the mapping properties of Bessel potential operators (see \eqref{e1.24}) and the mapping properties of a Mellin convolution operator with an admissible kernel both operators
\begin{eqnarray}\label{e54}
\begin{array}{lcr}
\mathbf{\Lambda}^s_{-\gamma}\mathbf{K}^1_c&:&\widetilde{\mathbb{H}}^r_p(\mathbb{R}^+)
     \longrightarrow\mathbb{H}^{r-s}_p(\mathbb{R}^+),\\[1mm]
\mathbf{K}^1_c\mathbf{\Lambda}^s_{-c\,\gamma}&:&\widetilde{\mathbb{H}}^r_p(\mathbb{R}^+)
     \longrightarrow\mathbb{H}^{r-s}_p(\mathbb{R}^+)
\end{array}
 \end{eqnarray}
are correctly defined and bounded for all $s\in\mathbb{R}$, $1<p<\infty$, since  $-\pi<\arg (-\gamma)<0$ and $0<\arg(-c\,\gamma)<\pi$.

Second, let us consider the positive integer values $s=n=1,2,\ldots$. Then, with the help of formulae \eqref{e1.26} and \eqref{e1.16} it follows that:
\begin{align}\label{e51}
\mathbf{\Lambda}^n_{-\gamma}\mathbf{K}^1_c\varphi
&=\Big(i\,\frac{d}{dt}-\gamma\Big)^n\mathbf{K}^1_c\varphi=\sum_{k=0}^n\binom{n}{k}
    i^k(-\gamma)^{n-k} \,\frac{d^k}{dt^k}\,\mathbf{K}^1_c\varphi\nonumber \\
& =\sum_{k=0}^n\binom{n}{k}i^k(-\gamma)^{n-k}c^{-k}\Big(\mathbf{K}^1_c\,
    \frac{d^k}{dt^k}\,\varphi\Big)(t)= \nonumber \\
& =c^{-n}\mathbf{K}^1_c\bigg(\sum_{k=0}^n \binom{n}{k}i^k\left(- c\,\gamma
    \right)^{n-k}\,\frac{d^k}{dt^k}\,\varphi\bigg)(t)= \nonumber \\
& =c^{-n}\mathbf{K}^1_c\mathbf{\Lambda}^n_{- c\,\gamma}\varphi,\qquad\varphi\in
    \widetilde{\mathbb{H}}^r_p(\mathbb{R}^+)
\end{align}
and we have proven formula \eqref{e1.38} for positive integers $s=n=1,2,\ldots$.

For negative $s=-1,-2,\ldots$ formulae \eqref{e1.38} follows if we apply the inverse operator $\mathbf{\Lambda}^{-n}_{-\gamma}$ and $\mathbf{\Lambda}^{-n}_{-c\gamma}$ to the proved operator equality
 \[
\mathbf{\Lambda}^n_{-\gamma}\mathbf{K}^1_c=c^{-n}\mathbf{K}^1_c\mathbf{\Lambda}^n_{-c\,\gamma}
 \]
for positive $n=1,2,\ldots$ from the left and from the right, respectively. We obtain
 \[
\mathbf{K}^1_c\mathbf{\Lambda}^{-n}_{-c\gamma}=c^{-n}\mathbf{\Lambda}^{-n}_{-\gamma}
     \mathbf{K}^1_c\quad {\rm or}\quad\mathbf{\Lambda}^{-n}_{-\gamma}\mathbf{K}^1_c
     =c^n\mathbf{K}^1_c\mathbf{\Lambda}^{-n}_{-c\gamma}
 \]
and \eqref{e1.38} is proved also for a negative $s=-1,-2,\ldots$.

In order to derive formula \eqref{e1.38} for non-integer values of $s$, we can confine ourselves to the case $-2<s<-1$. Indeed, any non-integer value $s\in \mathbb{R}$ can be represented in the form $s=s_0+m$, where $-2<s_0<-1$ and $m$ is an integer. Therefore, if for $s=s_0+m$ the operators in \eqref{e54} are correctly defined and bounded, and if the relations in question are valid for $-2<s_0<-1$, then we can write
\begin{align*}
\mathbf{\Lambda}^s_{-\gamma}{\bf K}^1_c=\mathbf{\Lambda}^{s_0 +m}_{-\gamma}{\bf K}^1_c
     &=c^{-m}\mathbf{\Lambda}^{s_0}_{-\gamma}{\bf K}^1_c\mathbf{\Lambda}^{m}_{-c\,\gamma}
     =c^{-s_0-m}{\bf K}^1_c\mathbf{\Lambda}^{s_0}_{-c\,\gamma} \mathbf{\Lambda}^{m}_{-c\, \gamma}\nonumber\\
&=c^{-s_0-m}{\bf K}^1_c\mathbf{\Lambda}^{s_0+m}_{-c\,\gamma}=c^{-s}{\bf K}^1_c
     \mathbf{\Lambda}^s_{-c\,\gamma}.
\end{align*}

Thus let us assume that $-2<s<-1$ and consider the expression
  \begin{equation}\label{e1.19} 
  \mathbf{\Lambda}^s_{-\gamma}{\bf K}^1_c\varphi(t)\,\,=\frac1{2\pi^2}\,r_+\int\limits_{-\infty}^\infty
    e^{-i\xi t} (\xi-\gamma)^s\int\limits_0^\infty e^{i\xi y}\int\limits_0^\infty
    \frac{\varphi(\tau)}{y-c\tau}\;d\tau\,dy\,d\xi,
 \end{equation}
where $r_+$ is the restriction to $\mathbb{R}^+$. It is clear that the integral in the right-hand-side of \eqref{e1.19} exists. Indeed, if $\varphi\in\mathbb{L}_2$, then ${\bf K}^1_c\varphi\in\mathbb{L}_2\cap C^\infty$ and $\mathbf{\Lambda}^s_{-\gamma}{\bf K}^1_c\varphi\in\mathbb{H}^{-s}\cap C^\infty\subset\mathbb{L}_2\cap C^\infty$.

Now consider the function $e^{-izt}(z-\gamma)^se^{i z y}$, $z\in \mathbb{C}$. Since ${\rm Im}\gamma\not=0$, $s<-1$, then for sufficiently small $\varepsilon>0$ this function is analytic in the
strip between the lines $\mathbb{R}$ and $\mathbb{R}+i\varepsilon$ and vanishes  at the infinity for all finite $t\in\mathbb{R}$ and for all $y>0$. Therefore, the integration over the real line $\mathbb{R}$ in the first integral of \eqref{e1.19} can be replaced by the integration over the line $\mathbb{R}+i\varepsilon$, i.e.
 \begin{equation}\label{e1.21a} 
\hskip-1mm\mathbf{\Lambda}^s_{-\gamma}{\bf
K}^1_c\varphi(t)=\frac1{2\pi^2}\,r_+\hskip-3mm\int\limits_{-\infty}^\infty
e^{-i\xi t +\varepsilon t} (\xi+i\varepsilon-\gamma)^s
    \int\limits_0^\infty e^{i\xi  y -\varepsilon y}\int\limits_0^\infty
    \frac{\varphi(\tau)}{y-c\tau}\;d\tau\,dy\,dx.
\end{equation}
Let us use the density of the set $C^\infty_0(\mathbb{R}^+)$ in $\widetilde{\mathbb{H}}^s_p(\mathbb{R}^+)$. Thus for all finite $t\in\mathbb{R}$ and for all functions $\varphi\in C^\infty_0(\mathbb{R})$ with compact supports the integrand in the corresponding triple integral for \eqref{e1.21a} is absolutely integrable. Therefore, for such functions one can use Fubini-Tonelli theorem and change the order of integration in \eqref{e1.21a}. Thereafter, one returns to the integration over the real line $\mathbb{R}$ and obtains
 \begin{equation}\label{e1.19a}
 \mathbf{\Lambda}^s_{-\gamma}{\bf K}^1_c\varphi(t)\,\,=\frac1{2\pi^2}\,r_+\int\limits_0^\infty\varphi
    (\tau)\int\limits_0^\infty\frac1{y-c\tau}\int\limits_{-\infty}^\infty e^{i\xi(y-t)}(\xi-\gamma)^sd\xi\,dy\,d\tau,
 \end{equation}

In order to study the expression in the right-hand side of \eqref{e1.19a}, one can use a well known formula
\begin{equation*}
\begin{aligned}
\displaystyle
 &\int\limits_{-\infty}^\infty(\beta+ix)^{-\nu}e^{-ipx}\,dx=\begin{cases}
     0 \quad &\text{for} \quad p>0, \\
-\displaystyle\frac{2\pi(-p)^{\nu-1}e^{\beta\,p}}{\Gamma(\nu)}
     \quad&\text{for}\;\;p<0,\end{cases}\\[1.5ex]
& {\rm Re}\,\nu>0,\qquad{\rm Re}\beta>0 ,
\end{aligned}
\end{equation*}
 \cite[Formula 3.382.6]{GR94}. It can be rewritten in a more convenient form--viz.,
\begin{align}\label{eqn27}
&&\hskip-15mm\displaystyle\int\limits_{-\infty}^\infty
e^{i\mu\,\xi}(\xi
     -\gamma)^s\,d\xi=\left\{\begin{array}{ll}0 &\quad\text{ if } \;  \mu<0,\; {\rm Im}\,\gamma>0,\\
     \displaystyle\frac{2\pi\,\mu^{-s-1}e^{-\frac\pi2si+\mu\,\gamma i}}{
     \Gamma(-s)}&\quad\text{ if } \; \mu>0, \; {\rm Im}\,\gamma>0.\end{array}\right.
\end{align}

Applying \eqref{eqn27} to the last integral in \eqref{e1.19a}, one
obtains
\begin{align}\label{eqn28}
\mathbf{\Lambda}^s_{-\gamma}{\bf K}^1_c\varphi(t)&\,\,=\displaystyle\frac{e^{-\frac\pi2si}}{\pi\Gamma(-s)}
    r_+\displaystyle\int\limits_0^\infty\varphi(\tau)\,d\tau\displaystyle\int\limits_t^\infty\displaystyle\frac{e^{i(y-t) \gamma}dy}{(y-t)^{1+s} (y-c\tau)}\nonumber\\
&=\displaystyle\frac{e^{-\frac\pi2si}}{\pi\Gamma(-s)}r_+\displaystyle\int\limits_0^\infty\varphi(\tau)\,d\tau
    \displaystyle\int\limits_0^\infty\displaystyle\frac{y^{-s-1}e^{i\gamma\,y}dy}{y+t-c\tau},
\end{align}
where the integrals exist since  $-s-1>-1$ and $0<\arg \gamma<\pi$ (i.e., ${\rm Im}\,\gamma>0$).

Let us recall the formula
 \begin{equation}\label{eqn29}
  \begin{aligned}
\displaystyle
&\int\limits_0^\infty\frac{x^{\nu-1}e^{-\mu\,x}\,dx}{x+\beta}
     =\beta^{\nu-1}e^{\beta\,\mu}\Gamma(\nu)\Gamma(1-\nu,\beta\mu),\\[1ex]
&{\rm Re}\,\nu>0,\quad{\rm Re}\mu>0,\quad |\arg\,\beta|<\pi
 \end{aligned}
 \end{equation}
(cf. \cite[formula 3.383.10]{GR94}). Due to the conditions $0<\arg\,c<2\pi$, $t>0$, $\tau>0$ we have $|\arg(t-c\tau)|<\pi$ and, therefore, we can apply \eqref{eqn29} to the equality \eqref{eqn28}. Then \eqref{eqn28} acquires the following final form:
\begin{equation}\label{eqn30}
\mathbf{\Lambda}^s_{-\gamma}{\bf K}^1_c\varphi(t)
    =\frac{e^{-\frac\pi2si}}\pi r_+\displaystyle\int\limits_0^\infty\displaystyle
    \frac{e^{-i\gamma(t-c\tau)}\Gamma(1+s,-i\gamma(t-c\tau))\varphi(\tau)
    \,d\tau}{(t-c\tau)^{1+s}}.
\end{equation}

Consider now the reverse composition ${\bf K}^1_c\mathbf{\Lambda}^s_{-c\, \gamma} \varphi(t)$. Changing the order of integration in the corresponding expression (see \eqref{e1.19a} for a similar motivation), one  obtains
\begin{align}\label{eqn31}
{\bf
K}^1_c\mathbf{\Lambda}^s_{-c\,\gamma}\varphi(t)&:=\frac1{2\pi^2}\,r_+
    \int\limits_0^\infty\frac1{t-c\,y}\int\limits_{-\infty}^\infty e^{-i\xi\,y}(\xi-c\,\gamma)^s
    \int\limits_0^\infty e^{i\xi\,\tau}\varphi(\tau)d\tau\,d\xi\,dy\nonumber\\
&=\frac1{2\pi^2}r_+\int\limits_0^\infty\varphi(\tau)\int\limits_0^\infty
    \frac1{t-c\,y}\int\limits_{-\infty}^\infty e^{i\xi(\tau-y)}
    (\xi-c\,\gamma)^sd\xi\,dy\,d\tau.
\end{align}
In order to compute the expression in the right-hand side of \eqref{eqn31}, let us recall formula~3.382.7 from \cite{GR94}:
\begin{equation*}
\begin{aligned}
&\displaystyle\int\limits_{-\infty}^\infty(\beta-ix)^{-\nu}e^{-ipx}\,dx=\begin{cases}
     0\quad &\text{for}\quad p<0, \\
\displaystyle\frac{2\pi\,p^{\nu-1}e^{-\beta\,p}}{\Gamma(\nu)}
\quad&\text{for}
     \quad p>0,\end{cases}\\[1.5ex]
&{\rm Re}\,\nu>0,\quad{\rm Re}\beta>0
\end{aligned}
\end{equation*}
and rewrite it in a form more suitable for our consideration--viz.,
\begin{equation}\label{eqn32}
\begin{aligned}
&\displaystyle\int\limits_{-\infty}^\infty
e^{i\mu\,\xi}(\xi+\omega)^s
     \,d\xi=\left\{\begin{array}{ll}0 &\quad\mu>0,\;{\rm Im}\,\omega>0,\\
     \displaystyle\frac{2\pi\,(-\mu)^{-s-1}e^{\frac\pi2si-\mu\,\omega i}}{
     \Gamma(-s)}&\quad\mu<0,\;{\rm
     Im}\,\omega>0,\end{array}\right.\\[1ex]
& {\rm Re}\,s<0,\quad\mu\in\mathbb{R},\quad\omega,\,s\in\mathbb{C}.
\end{aligned}
\end{equation}
Using \eqref{eqn32}, we represent \eqref{eqn31} in the form
\begin{equation}\label{eqn32a}
 \begin{aligned}
{\bf K}^1_c\mathbf{\Lambda}^s_{-c\,\gamma}\varphi(t)
&=\displaystyle\frac{e^{\frac\pi2si}}{\pi\Gamma(-s)}\,r_+
    \displaystyle\int\limits_0^\infty \varphi(\tau)\,d\tau\displaystyle\int\limits_\tau^\infty
    \displaystyle\frac{e^{-ic\,\gamma(y-\tau)}\,dy}{(y-\tau)^{s+1}(t-c\,y)}\\
&=-\displaystyle\frac{e^{\frac\pi2si}}{\pi
c\Gamma(-s)}r_+\displaystyle
 \int\limits_0^\infty\varphi(\tau)\,d\tau\displaystyle
 \int\limits_0^\infty\displaystyle\frac{y^{-s-1}e^{-ic\gamma\,y}\,dy}{y+\tau-c^{-1}t},
\end{aligned}
 \end{equation}
where the integrals exist since $-s-1>-1$ and $-\pi<\arg(c\,\gamma)<0$  (i.e., ${\rm Im}\,c\,\gamma<0$).

Due to the conditions $0<\arg\,c<2\pi$, $t>0$, $\tau>0$ we have $|\arg(\tau-c^{-1}t)|<\pi$, Therefore, we can apply formula formula \eqref{eqn29} to  \eqref{eqn32a} and get the following representation:
\begin{align}\label{eqn33}
{\bf K}^1_c\mathbf{\Lambda}^s_{-c\,\gamma}\varphi(t)
     =&-\frac{c^{-1}e^{\frac\pi2si}}\pi
    r_+\int\limits_0^\infty\displaystyle\frac{e^{-ic\gamma(c^{-1}t-\tau)}\Gamma(1+s,
    -ic\gamma(c^{-1}t-\tau))\varphi(\tau)\,d\tau}{(\tau-c^{-1}t)^{1+s}}\nonumber\\
=& \frac{c^se^{-\frac\pi2si}}\pi r_+\int\limits_0^\infty
    \displaystyle\frac{e^{-i\gamma(t-c\,\tau)}\Gamma(1+s,-i\gamma(t-c\,\tau))
    \varphi(\tau)\,d\tau}{(t-c\,\tau)^{1+s}}.
\end{align}

If we multiply  \eqref{eqn33} by $c^{-s}$ we get precisely the expression in  \eqref{eqn30} and, therefore, $\mathbf{\Lambda}^s_{-\gamma}{\bf K}^1_c\varphi(t)=c^{-s}{\bf K}^1_c\mathbf{\Lambda}^s_{-c\,\gamma} \varphi(t)$, which proves the claimed equality \eqref{e1.38}  for $-2<s<-1$ and accomplishes the proof.    \hfill $\square$
 %
\begin{corollary}\label{c3.2}
Let $0<\arg\,c<2\pi$ and $0<\arg\gamma<\pi$. Then for arbitrary $\gamma_0\in\mathbb{C}$ such that
$0<\arg\,\gamma_0<\pi$ and $-\pi<\arg(c\,\gamma_0)<0$, one has
\begin{equation}\label{eqn36}
\mathbf{\Lambda}^s_{-\gamma}{\bf
K}^1_c=c^{-s}W_{g_{-\gamma,-\gamma_0}}{\bf K}^1_c
   \mathbf{\Lambda}^s_{-c\,\gamma_0},
\end{equation}
where
\begin{align}\label{eqn37}
g^s_{-\gamma,-\gamma_0}(\xi):=\left(\frac{\xi-\gamma}{\xi-\gamma_0}\right)^s.
\end{align}

If, in addition, $1<p<\infty$ and $1/p-1<r<1/p$ then equality
\eqref{eqn36} can be supplemented as follows:
\begin{align}\label{eqn38}
\mathbf{\Lambda}^s_{-\gamma}{\bf K}^1_c=c^{-s}\left[{\bf
K}^1_cW_{g^s_{-\gamma,-\gamma_0}}
   +\mathbf{T}\right]\mathbf{\Lambda}^s_{-c\,\gamma_0},
\end{align}
where $\mathbf{T}\;:\;\widetilde{\mathbb{H}}{}^r_p(\mathbb{R}^+)\to \mathbb{H}^r_p( \mathbb{R}^+)$ is a compact operator, and if $c$ is a real negative number, then $c^{-s}:=|c|^{-s}e^{-\pi si}$.
\end{corollary}
{\sc Proof.}  It follows from equalities {e28} and \eqref{e1.38} that
\begin{eqnarray*}
\mathbf{\Lambda}^s_{-\gamma}{\bf K}^1_c=\mathbf{\Lambda}^s_{-\gamma}
     \mathbf{\Lambda}^{-s}_{-\gamma_0}\mathbf{\Lambda}^s_{-\gamma_0}{\bf K}^1_c
     =c^{-s}W_{g_{-\gamma,-\gamma_0}}{\bf K}^1_c\mathbf{\Lambda}^s_{-c\,\gamma_0}
\end{eqnarray*}
and \eqref{eqn36} is proved. If $1<p<\infty$ and $1/p-1<r<1/p$, then
the commutator
 \[
\mathbf{T}:=W_{g^s_{-\gamma,-\gamma_0}}{\bf K}^1_c-{\bf
K}^1_cW_{g^s_{-\gamma,-\gamma_0}}
     \;:\;\widetilde{\mathbb{H}}{}^r_p(\mathbb{R}^+)\to\mathbb{H}^r_p(\mathbb{R}^+)
 \]
of Mellin and Fourier convolution operators is correctly defined and
bounded. It is compact for $r=0$ and all $1<p<\infty$ (see
\cite{Du74,Du87}). Due to Krasnoselsky's interpolation theorem (see
\cite{Kr60} and also \cite[Sections 1.10.1 and 1.17.4]{Tr95}), the
operator $\mathbf{T}$ is compact in all $\mathbb{L}_r$-spaces for
$1/p-1<r<1/p$. Therefore, the equality \eqref{eqn36} can be
rewritten as
\begin{eqnarray*}
\mathbf{\Lambda}^s_{-\gamma}{\bf K}^1_c=c^{-s}\left[{\bf
K}^1_cW_{g^s_{-\gamma,-\gamma_0}}
     +\mathbf{T}\right]\mathbf{\Lambda}^s_{-c\,\gamma_0}\, ,
\end{eqnarray*}
and we are done.   \QED
\begin{remark}\label{r2.10}
The assumption $1/p-1<r<1/p$ in \eqref{eqn38} cannot be relaxed.
Indeed, the operator $W_{g^s_{-\gamma,-\gamma_0}}{\bf
K}^1_c=\mathbf{\Lambda}^s_{-\gamma}
\mathbf{\Lambda}^{-s}_{-\gamma_0}{\bf K}^1_c\;:\;\widetilde{
\mathbb{H}}{}^r_p(\mathbb{R}^+)\to\mathbb{H}^r_p(\mathbb{R}^+)$ is
bounded for all $r\in\mathbb{R}$ (see \eqref{e54}). But the
operator ${\bf K}^1_cW_{g^s_{-\gamma,
-\gamma_0}}\;:\;\widetilde{\mathbb{H}}{}^r_p
(\mathbb{R}^+)\to\mathbb{H}^r_p(\mathbb{R}^+)$ is bounded only for
$1/p-1<r<1/p$ because the function $g^s_{-\gamma,-\gamma_0}(\xi)$
has an analytic extension into the lower half-plane but not into the
upper one.
\end{remark}

\section{A Local Principle}\label{s3}

In the present section we expose well known, but slightly modified local principle from \cite{Si65}, which we apply intensively.

Let $\mathfrak{B}_1(\Omega)$and $\mathfrak{B}_2(\Omega)$ be Banach spaces of functions on a domain $\Omega\subset \mathbb{R}^n$ and multiplication by uniformly bounded $C^\infty(\overline\Omega)$-functions are bounded operators in both spaces. If $\Omega=\bR^n$, we consider one point compactification $\overline\Omega:=\dR{}^n$ of $\Omega=\bR^n$.

Let $x\in\overline\Omega$ and consider the class of multiplication operators by functions
 \begin{equation}\label{e2.3.1}
 \begin{array}{r}
\Delta_x:=\Big\{vI\;:\;v\in C^\infty(\Omega),
\;\; v(t)=1 \;\; \text{for}\;\; |t-x|<\ve_1,\quad v(x)\geqslant0\\[3mm]
\text{and}\;\; v(t)=0\;\; \text{for}\;\; |t-x|>\ve_2 \Big\}
 \end{array}
 \end{equation}\index{N}{$\Delta_x$}
where $\ve_2>\ve_1>0$ are not fixed and vary from function to function. $\Delta_x$ is, obviously, a localizing class in the algebra of bounded linear operators $\cL(\mathfrak{B}_1(\Omega),\mathfrak{B}_2(\Omega))$ and $\left\{\Delta_x\right\}_{x\in\overline\Omega}$ is a covering class. Indeed, for a system $\{v_xI\}_{x\in\overline\Omega}$ we consider the related covering
 \[
\overline\Omega=\bigcup_{x\in\overline\Omega}U_x,\quad U_x:=\{y\in\overline\Omega\;:\;v_x(y)=1\}\, .
 \]
The set $\overline\Omega$ is compact and there exists a finite covering system
$\overline\Omega=\bigcup_{j=1}^NU_{x_j}$. The corresponding sum is strictly
positive
 \begin{eqnarray}\label{e3.299}
\inf_{y\in\bR^n}g(y)\geqslant1\quad \mbox{\rm for} \quad
    g(y):=\sum_{j=1}^Nv_{x_j}(y)
\end{eqnarray}
and the multiplication operator $\sum_{j=1}^Nv_{x_j}I=gI$ has the inverse $g^{-1}I$. Thus, the system of localizing classes $\{\Delta_x\}_{x\in\overline\Omega}$ is covering.

Most probably we would have to deal with the quotient space $\cL_0'(\mathfrak{B}_1(\Omega),\mathfrak{B}_2(\Omega))$ $:=\cL(\mathfrak{B}_1(\Omega),\mathfrak{B}_2(\Omega))/
\mathcal{C}(\mathfrak{B}_1(\Omega),\mathfrak{B}_2(\Omega))$ of linear bounded operators with respect to the
compact operators.
 %
\begin{definition}\label{d3.1}
A quotient class $[\A]\in\cL'(\mathfrak{B}_1(\Omega),\mathfrak{B}_2(\Omega))$ is called $\Delta_x$-invertible if there exists a quotient class $[\mathbf{R}_x]\in\cL'(\mathfrak{B}_2(\Omega),\mathfrak{B}_1(\Omega))$ and  $v_x\in \Delta_x$ such that the operator equalities $[\mathbf{R}_x\mathbf{A}v_xI_1]=[v_xI_1]$ and $[v_x\mathbf{A}\mathbf{R}_x]=[v_xI_2]$ holds, where $I_1$ and $I_2$ are the identity operators in the spaces $\mathfrak{B}_1(\Omega)$ and $\mathfrak{B}_2(\Omega)$.
\end{definition}

Consider a pair of operators
\begin{eqnarray}\label{e2.3.3a}
\mathbf{A}_j\;:\; \mathfrak{B}_1(\Omega_j)\to \mathfrak{B}_2(\Omega_j),\qquad j=1,2,
 \end{eqnarray}
in the same pairs of function spaces $\mathfrak{B}_1(\Omega_1)$, $\mathfrak{B}_2(\Omega_1)$ and $\mathfrak{B}_1(\Omega_2)$, $\mathfrak{B}_2(\Omega_2)$ defined on different domains $\Omega_1, \Omega_2\subset \mathbb{R}^n$. For this we assume that for any pair of points $x_1\in\overline\Omega_1$ and $x_2\in\overline\Omega_2$ there exist there exists a local diffeomorphism of neighbourhoods
\begin{eqnarray}\label{e2.3.3b}
\beta\;:\;\omega(x_1)\to\omega(x_2),\qquad x_j\in\omega(x_j)\subset\Omega_j,\quad j=1,2.
 \end{eqnarray}
The operators
 \[
\beta_*\vf(x):=\vf(\beta(x)),\qquad \beta^{-1}_*\psi(y):=\psi(\beta^{-1}_(y))
 \]
are inverses to each-other and map the spaces
 \[
\beta_*\;:\;\mathfrak{B}_j(\omega_2)\to\mathfrak{B}_j(\omega_1),\qquad
       \beta^{-1}_*\;:\;\mathfrak{B}_j(\omega_1)\to\mathfrak{B}_j(\omega_2).
 \]
 %
\begin{definition}\label{d3.2}[Quasi localization]\label{d2.3.2}
Let multiplication by uniformly bounded $C^\infty$ functions on corresponding closed domains $\overline{\Omega_1}$ and $\overline{\Omega_1}$ are bounded operators in all respective spaces $\mathfrak{B}_2(\Omega_1)$ and $\mathfrak{B}_1(\Omega_2)$, $\mathfrak{B}_2(\Omega_2)$.

Two classes from the quotient spaces $[\A_1], [\A_2]\in\cL'(\mathfrak{B}_1,\mathfrak{B}_2)$ (see \eqref{e2.3.3a}) are called locally quasi equivalent at $x_1\in\overline\Omega_1$ and $x_2\in\overline\Omega_2$, if
\begin{eqnarray}\label{e2.3.3}
\begin{array}{c}
\inf\limits_{v^1_{x_1},v^2_{x_1}\in\Delta_{x_1}}\|\!|[v^1_{x_1}][\mathbf{A}_1
     -\beta_*\mathbf{A}_2\beta^{-1}_*]\|\!|=\inf\limits_{v^1_{x_1},v^2_{x_1}\in\Delta_{x_1}}\|\!|[\mathbf{A}_1
     -\beta_*\mathbf{A}_2\beta^{-1}_*][v^2_{x_1}I]\|\!|=0,
\end{array}
 \end{eqnarray}
where the norm in the quotient space $\cL'(\mathfrak{B}_1,\mathfrak{B}_2)=\cL(\mathfrak{B}_1,\mathfrak{B}_2)/ \mathcal{C}(\mathfrak{B}_1,\mathfrak{B}_2)$ coincides with the essential norm
 \[
\|[\A]\|:=\|\!|\A\|\!|:=\inf_{T\in\mathcal{C}(\mathfrak{B}_1,\mathfrak{B}_2)}\|A+T\|.
 \]

Such an equivalence we denote as follows $[\mathbf{A}_1]\overset{\Delta_{x_1}}{\sim} \beta \overset{\Delta_{x_2}}{\sim}[\mathbf{A}_2]$ or also $[\mathbf{A}_1]\overset{x_1}{ \sim}\beta\overset{x_2}{\sim}[\mathbf{A}_2]$.

If $\Omega_1=\Omega_2=\Omega$ and $\beta(x)=x$ is the identity map, the equivalence at the point $x\in\Omega$ is denoted as follows $[\mathbf{A}_1]\overset{\Delta_x}\sim [\mathbf{A}_2]$ or also $[\mathbf{A}_1]\overset{x}\sim[\mathbf{A}_2]$.
\end{definition}
 %
\begin{definition}\label{d3.3}
Let  $\mathbf{A}$, $\mathfrak{B}_1(\Omega)$ and $\mathfrak{B}_2(\Omega)$ be the same as in Definition \ref{d3.2}. An operator $\mathbf{A}\;:\; \mathfrak{B}_1(\Omega)\to \mathfrak{B}_2(\Omega)$ is called of local type if $v_1\A v_2I\;:\; \mathfrak{B}_1(\Omega)\to \mathfrak{B}_2(\Omega)$ is compact for all $v_1,v_2\in C^\infty(\Omega)$, provided $\supp\,v_1\cap\supp\,v_2=\emptyset$ (see \cite{Se66});
Or, equivalently, if $v\A-\A vI\;:\; \mathfrak{B}_1(\Omega)\to \mathfrak{B}_2(\Omega)$ is compact for all  $v\in C^\infty(\Omega)$ (see \cite{Se66}).
\end{definition}
 %
\begin{proposition}[\bf Localization principle]\label{p3.4}
Let  $\mathbf{A}$, $\mathfrak{B}_j(\Omega_k)$, $j,k=1,2$, be a group of four function spaces as in Definition \ref{d2.3.2} and
 \[
\mathbf{A}\;:\; \mathfrak{B}_1(\Omega_1)\to \mathfrak{B}_2(\Omega_1),\qquad
\mathbf{B}_x\;:\; \mathfrak{B}_1(\Omega_2)\to \mathfrak{B}_2(\Omega_2),\qquad x\in\Omega_1,
\]
be operators of local type.

If the quasi equivalence $[\mathbf{A}]\overset{x}{\sim}\beta_x \overset{y}{\sim}[\mathbf{B}_x]$ holds for some diffeomorphism  $\beta_x\;:\;\omega(x)\to\omega(y(x))$, $y(x)\in\Omega_2$, then $[\mathbf{A}]$ is locally invertible at $x\in\Omega_1$ if and only if $[\mathbf{B}_x]$ is locally invertible at $y(x)$.

If the quasi equivalence $[\mathbf{A}]\overset{x}{\sim}\beta_x \overset{y(x)}{\sim}[\mathbf{B}_x]$ holds for all $x\in\overline{\Omega_1}$ and $[\mathbf{B}_x]\in\cL'(\mathfrak{B}_1(\Omega_2),\mathfrak{B}_2(\Omega_2))$ are locally invertible at $y(x)\in\Omega_2$ for all $y(x)\in\overline\Omega$, than the quotient class $[\A]$ is globally invertible (i.e., $\mathbf{A}\;:\; \mathfrak{B}_1(\Omega_1)\to \mathfrak{B}_2(\Omega_1)$ is a Fredholm operator).
\end{proposition}
 %
\begin{remark}\label{r2.3.4}
If in the foregoing Proposition \ref{p3.4} we drop the condition that $\mathbf{A}$ and $\mathbf{B}_x$ are of local type, then from the left (from the right) quasi equivalence and the left invertibility of $[\mathbf{B}_x]\in\cL'(\mathfrak{B}_1(\Omega_2),\mathfrak{B}_2 (\Omega_2))$ at $y(x)\in\Omega_2$ for all $y(x)\in\overline\Omega$, follows  is global invertibility of the quotient class $[\A]$ from the left (from the right), i.e., the existence of the left (of the right) regularizer for the operator $\mathbf{A}\;:\; \mathfrak{B}_1(\Omega_1)\to \mathfrak{B}_2(\Omega_1)$).
\end{remark}

\section{Algebra Generated by Mellin and Fourier Convolution
Operators}\label{s4}

Let $\dR:=\bR\cup\{\infty\}$ denote one point compactification of the real axes $\bR$ and $\overline{\bR}:=\bR\cup\{\pm\infty\}$-the two point compactification of $\bR$. By $C(\dR)$ (by $C(\overline{\bR})$, respectively) we denote the space of continuous functions $g(x)$ on $\bR$ which have the equal limits at the infinity $g(-\infty)=g(+\infty)$ (limits at the infinity can be different $g(-\infty)\not=g(+\infty)$. By $PC(\dR)$is denoted the space of piecewise-continuous functions on $\dR$, having the limits $a(t\pm0)$ at all points $t\in\dR$, including the infinity.

Unlike the operators $W^0_a$ and $\mathfrak{M}^0_a$ (see Section~1), possessing the property
\begin{equation}\label{e2.1}
    W^0_aW^0_b=W^0_{ab},\;\; \mathfrak{M}^0_a\mathfrak{M}^0_b=\mathfrak{M}^0_{ab}\;\; \text{for all}\;\; a,b\in\mathfrak{M}_p(\mathbb{R}),
\end{equation}
the composition of the convolution operators on the semi-axes $W_a$ and $W_b$ cannot be computed by the rules similar to \eqref{e2.1}. Nevertheless, the following propositions hold.
 %
\begin{proposition}[{\cite[Section 2]{Du79}}]\label{p2.1}
Let $1<p<\infty$ and $a,b\in\mathfrak{M}_p(\overline{\mathbb{R}}{}^+)\cap PC(\dR)$
be scalar $\mathbb{L}_p$-multipliers, piecewise-continuous on $\bR$ including infinity. Then the commutant $[W_a,W_b]:=W_aW_b-W_bW_a$ of the operators $W_a$ and $W_b$ is a compact operator in the Lebesgue space $[W_a,W_b]:\mathbb{L}_p(\mathbb{R}^+) \longmapsto \mathbb{L}_p(\mathbb{R}^+)$.

Moreover, if, in addition, the symbols $a(\xi)$ and $b(\xi)$ of the operators $W_a$ and $W_b$ have no common discontinuity points, i.e., if
 \[
\big[a(\xi+0)-a(\xi+0)\big]\big[b(\xi+0)-b(\xi+0)\big]=0\;\; \text{for all}\;\;
     \xi\in\dR,
 \]
then $\mathbf{T}=W_aW_b-W_{ab}$ is a compact operator in $\mathbb{L}_p(\mathbb{R}^+)$.
\end{proposition}

Note that the algebra of $N\times N$ matrix multipliers $\mathfrak{M}_2(\mathbb{R})$  coincides with the algebra of $N\times N$ matrix functions essentially  bounded on $\mathbb{R}$. For $p\not=2$,  the algebra $\mathfrak{M}_p(\mathbb{R})$ is rather complicated. There are multipliers $g\in\mathfrak{M}_p(\mathbb{R})$ which are elliptic, i.e. ${\rm ess}\,\inf|g(x)|>0$, but $1/g\not\in\mathfrak{M}_p(\mathbb{R})$. In connection with this, let us consider the subalgebra $PC\mathfrak{M}_p(\mathbb{R})$ which is the closure of the algebra of piecewise-constant functions on $\mathbb{R}$ in the norm of multipliers $\mathfrak{M}_p(\mathbb{R})$
\[  \big\|a\mid\mathfrak{M}_p(\mathbb{R})\big\|:=\big\|W^0_a\mid\mathbb{L}_p(\mathbb{R})\big\|.       \]
Note that any function $g\in PC\mathfrak{M}_p(\mathbb{R})\subset PC(\mathbb{R})$ has limits $g(x\pm0)$ for all $x\in\overline{\mathbb{R}}$, including the infinity. Let
\[  C\mathfrak{M}_p(\overline{\mathbb{R}}):=C(\overline{\mathbb{R}})\cap PC\mathfrak{M}_p^0(\mathbb{R}),\quad
        C\mathfrak{M}_p^0(\dR):=C(\dR)\cap PC\mathfrak{M}_p(\mathbb{R}),     \]
where functions $g\in C\mathfrak{M}_p(\overline{\mathbb{R}})$ (functions $h\in C(\dR)$) might have jump only at the infinity $g(-\infty)\not=g(+\infty)$ (are continuous at the infinity $h(-\infty)=h(+\infty)$).

$PC\mathfrak{M}_p(\mathbb{R})$ is a Banach algebra and contains all functions of bounded variation as a subset for all $1<p<\infty$ (Stechkin's theorem, see \cite[Section~2]{Du79}). Therefore, $\coth\pi(i\beta+\xi)\in C\mathfrak{M}_p(\overline{\mathbb{R}})$ for all $p\in(1,\infty)$.

\begin{proposition}[{\cite[Section 2]{Du79}}]\label{p2.2}
If $g\in PC\mathfrak{M}_p(\overline{\mathbb{R}})$ is an
$N\times N$ matrix multiplier, then its inverse $g^{-1}\in PC\mathfrak{M}_p(\overline{\mathbb{R}})$ if and only if it is elliptic, i.e. $\det g(x\pm0)\neq 0$ for all $x\in\overline{\mathbb{R}}$. If this is the case, the corresponding Mellin convolution operator $\mathfrak{M}^0_g:\mathbb{L}_p(\mathbb{R}^+) \longmapsto \mathbb{L}_p(\mathbb{R}^+)$ is invertible and
$(\mathfrak{M}^0_g)^{-1}=\mathfrak{M}^0_{g^{-1}}$.

Moreover, any $N\times N$ matrix multiplier $b\in C\mathfrak{M}^0_p(\dR)$ can be approximated by polynomials
\[  r_n(\xi):=\sum_{j=-m}^mc_m\Big(\frac{\xi-i}{\xi+i}\Big)^m,\;\; r_m\in C\mathfrak{M}^0_p(\overline{\mathbb{R}}),      \]
with constant $N\times N$ matrix coefficients, whereas any $N\times N$ matrix multiplier $g\in C\mathfrak{M}^0_p(\overline{\mathbb{R}})$ having a jump
discontinuity at infinity can be approximated by $N\times N$
matrix functions $d\coth\pi(i\beta +\xi)+r_m(\xi)$, $0<\beta<1$.
\end{proposition}

Due to the connection between the Fourier and Mellin convolution
operators (see Introduction, \eqref{e0.4}), the following is a
direct consequence of Proposition~\ref{p2.2}.

\begin{corollary}\label{c2.3}
The Mellin convolution operator
\[  \A=\mathfrak{M}^0_{\cA_\beta}:\mathbb{L}_p(\mathbb{R},
t^\gamma),         \]in \eqref{e0.1} with the symbol
$\mathcal{A}_\beta(\xi)$ in \eqref{e0.5} is invertible if and
only if the symbol is elliptic,
\begin{equation}\label{e2.2}
    \inf_{\xi\in\mathbb{R}}\big|\det \mathcal{A}_\beta(\xi)\big|>0
\end{equation}
and the inverse is then written as
$\mathbf{A}^{-1}=\mathfrak{M}^0_{\cA^{-1}_{1/p}}$.
\end{corollary}

The Hilbert transform on the semi-axis
\begin{equation}\label{e2.3}
    S_{\mathbb{R}^+}\varphi(x):=\frac1{\pi i}\int\limits_0^\infty\frac{\varphi(y)\,dy}{y-x}
\end{equation}
is the Fourier convolution $S_{\mathbb{R}^+}=W_{-{\sign}}$ on the
semi-axis $\mathbb{R}^+$ with the discontinuous symbol
$-\sign\xi$ (see \cite[Lemma 1.35]{Du79}), and it is also the
Mellin convolution
\begin{gather}\label{e2.4}
S_{\bR^+}={\mathfrak{M}_{s_\beta}^0}={\bf Z}_\beta W^0_{s_\beta}{\bf Z}^{-1}_\beta,\\
s_\beta(\xi):=\coth\pi(i\beta+\xi)\!=\!\frac{e^{\pi(i\beta+\xi)}\!
    +\!e^{-\pi(i\beta+\xi)}}{e^{\pi(i\beta+\xi)}\!-\!e^{-\pi(i\beta+\xi)}}
    \!=\!-i\cot\pi(\beta\!-\!i\xi),\; \xi\!\in\!\mathbb{R} \nonumber
\end{gather}
(cf. \eqref{e0.5} and \eqref{e1.4}). Indeed, to verify \eqref{e2.4} rewrite $S_{\mathbb{R}^+}$ in the following form
\[
S_{\mathbb{R}^+}\varphi(x):=\frac1{\pi i}\int\limits_0^\infty \frac{\varphi(y)}{1
    -\frac xy}\,\frac{dy}y=\int\limits_0^\infty K\Big(\frac{x}{y}\Big)\varphi(y)\,\frac{dy}y,        \]
where $K(t):=(1/\pi i)(1-t)^{-1}$. Further, using the formula
$$  \int\limits_0^\infty \frac{t^{z-1}}{1-t}\,dt=\pi\cot\pi z,\;\;\RRe z<1,     $$
cf. \cite[formula 3.241.3]{GR94}, one shows that the Mellin
transform $\mathcal{M}_\beta K(\xi)$ coincides with the
function $s_\beta(\xi)$ from \eqref{e2.4}.

Next Theorem \ref{t3.4} is an enhancement of Theorem \ref{e1.6}.
 %
\begin{theorem}\label{t3.4}
Let $1<p<\infty$ and $s\in\bR$. For arbitrary $\gamma_j\in\bC$, $\IIm\,\gamma_j>=0$ (j=1,2) the Hilbert transform
\begin{gather}\label{e1.38x}
\K^1_1=-iS_\bR^+=-iW_{-\sign}=W_{i\sign}\;:\;\wt{\mathbb{H}}^s_p(\mathbb{R}^+) \longrightarrow\mathbb{H}^s_p(\mathbb{R}^+)
\end{gather}
(see \eqref{e1.13}, \eqref{e1.12} and \eqref{e2.3}; the case  $c=1$, $\arg\,c=0$, Theorem \ref{t1.6}). $\K^1_1$ is a Fourier convolution operator and
\begin{gather}\label{e76}
\mathbf{\Lambda}^s_{-\gamma_1}\mathbf{K}^1_1\mathbf{\Lambda}^{-s}_{\gamma_2}
     =W_{i\sign g^s_{-\gamma_1,\gamma_2}}\;:\;\mathbb{L}_p(\mathbb{R}^+) \longrightarrow\mathbb{L}_p(\mathbb{R}^+),
\end{gather}
where $g^s_{-\gamma_1,\gamma_2}(\xi)$ is defined in \eqref{e1.11}.
\end{theorem}
{\bf Proof:} Formula \eqref{e76} follows from \eqref{e1.25a} and \eqref{e1.38x}. \QED

We need certain results concerning the compactness of Mellin and Fourier convolutions in
$\mathbb{L}_p$-spaces. These results are scattered in literature. For the convenience of the reader, we reformulate them here as Propositions \ref{p3.5}--\ref{p3.9}. For more details, the reader can consult \cite{Co69, Du79, Du87}.
 %
\begin{proposition}[{\cite[Proposition 1.6]{Du87}}]\label{p3.5}
Let $1<p<\infty$, $a\in C(\dR{}^+)$, $b\in C\mathfrak{M}^0_p(\dR)$ and $a(0)=b(\infty)=0$. Then the operators $a\mathfrak{M}^0_b,\mathfrak{M}^0_b\,aI:\mathbb{L}_p(\mathbb{R}^+)
\longrightarrow\mathbb{L}_p(\mathbb{R}^+)$ are compact.
\end{proposition}
 %
\begin{proposition}[{\cite[Lemma 7.1]{Du79} and \cite[Proposition 1.2]{Du87}}] \label{p3.6} Let\\ $1<p<\infty$, $a\in C(\dR{}^+)$, $b\in C\mathfrak{M}^0_p(\dR)$ and $a(\infty)= b(\infty)=0$. Then the operators $aW_b,W_b\,aI:\mathbb{L}_p(\mathbb{R}^+)\longrightarrow \mathbb{L}_p(\mathbb{R}^+)$ are compact.
\end{proposition}
 %
\begin{proposition}[{\cite[Lemma 2.5, Lemma 2.6]{Du87} and \cite{Co69}}]\label{p3.7}
Assume that $1<p<\infty$. Then
\begin{enumerate}
\vskip+0.15cm
\item If $g\in C\mathfrak{M}^0_p(\dR)$ and $g(\infty)=0$, the Hankel operator $H_g:\mathbb{L}_p(\mathbb{R}^+)\longrightarrow\mathbb{L}_p(\mathbb{R}^+)$ is compact;

\vskip+0.15cm
\item If the functions $a\in C(\dR)$, $b\in C\mathfrak{M}^0_p(\overline{\mathbb{R}})$, $c\in C(\overline{\mathbb{R}}{}^+)$ and satisfy at least one of the conditions
\begin{itemize}
\vskip+0.15cm
\item[(i)] $c(0)=b(+\infty)= 0$ and $a(\xi)=0$ for all $\xi>0$,

\vskip+0.15cm
\item[(ii)] $c(0)=b(-\infty)= 0$ and $a(\xi)=0$ for all $\xi<0$,
\end{itemize}
then the operators $cW_a\mathfrak{M}^0_b,\,c\mathfrak{M}^0_bW_a,\,W_a\mathfrak{M}^0_b\,cI,\,\mathfrak{M}^0_bW_a\,cI: \mathbb{L}_p(\mathbb{R}^+) \linebreak      \longrightarrow\mathbb{L}_p(\mathbb{R}^+)$ are compact.
\end{enumerate}
\end{proposition}
{\bf Proof:}
Let us comment only on item 1 in Proposition \ref{p3.7}, which
is not proved in \cite{Du87}, although is well known. The kernel
$k(x+y)$ of the operator $H_a$ is approximated by the Laguerre
polynomials $k_m(x+y)=e^{-x-y}p_m(x+y)$, $m=1,2,\ldots\;$, where
$p_m(x+y)$ are polynomials of order $m$ so that the corresponding
Hankel operators converge in norm $\|H_a-H_{a_m}|\,|\cL(\mathbb{L}_p(\mathbb{R}^+))\|\!\longrightarrow0$, where $a_m=\cF k_m$ are the Fourier transforms of the Laguerre polynomials (see, e.g. \cite{GF74}). Since
\[  |k_m(x+y)|=\big|e^{-x-y}p_m(x+y)\big|\leqslant C_m e^{-x}e^{-y}x^my^m, \;\; m=1,2,\ldots\,,        \]
for some constant $C_m$, the condition on the kernel
\[  \int\limits_0^\infty\bigg[\int\limits_0^\infty |k_m(x+y)|^{p'}\,dy\bigg]^{p/p'}\,dx<\infty, \;\; p':=\frac{p}{p-1} \,,      \]
holds and ensures the compactness of the operator $H_{a_m}:\mathbb{L}_p(\mathbb{R}^+)\longrightarrow\mathbb{L}_p(\mathbb{R}^+)$. Then the limit operator $H_a=\lim\limits_{m\longrightarrow\infty}H_{a_m}$ is compact as well.
\QED
 %
\begin{proposition}[{\cite[Lemma 7.1]{Du79} and \cite[Proposition 1.2]{Du87}}]\label{p3.8} Let\\ $1<p<\infty$, $a\in C(\dR{}^+)$, $b\in C\mathfrak{M}^0_p(\dR)$ and $a(\infty)= b(\infty)=0$. Then the operators $aW_b,W_b\,aI:\mathbb{L}_p(\mathbb{R}^+)\longrightarrow \mathbb{L}_p(\mathbb{R}^+)$ are compact.
\end{proposition}
 %
\begin{proposition}[{\cite[Lemma 7.4]{Du79} and \cite[Lemma 1.2]{Du87}}]\label{p3.9}
Let $1<p<\infty$ and let $a$ and $b$ satisfy at least one of
the conditions
\begin{itemize}
\vskip+0.15cm
\item[(i)] $a\in C(\overline{\mathbb{R}}{}^+)$, $b\in\mathfrak{M}^0_p(\mathbb{R})\cap PC(\overline{\mathbb{R}})$,

\vskip+0.15cm
\item[(ii)] $a\in PC(\overline{\mathbb{R}}{}^+)$, $b\in C\mathfrak{M}^0_p(\overline{\mathbb{R}})$.
\end{itemize}
Then the commutants $[aI,W_b]$ and $[aI,\mathfrak{M}^0_b]$ are
compact operators in the space $\mathbb{L}_p(\mathbb{R}^+)$.
\end{proposition}
 %
\begin{remark}\label{r3.10}
Note that, if both, a symbol $b$ and a function $a$, have jumps at finite points, the commutants $[aI,W_b]$ and $[aI,\mathfrak{M}^0_b]$ are not compact. Only jumps of a symbol at the infinity does not matters.
\end{remark}
 %
\begin{proposition}[\cite{Du87}]\label{p3.11}
The Banach algebra, generated by the Cauchy singular integral operator $S_\mathbb{R^+}$ and the identity operator $I$ on the semi-axis $\mathbb{R}^+$, contains Fourier convolution operators with symbols having discontinuity of the jump type only at zero and at the infinity and Mellin convolution operators with continuous symbols on $\dR$ (including the uinfinity).

Moreover, the Banach algebra $\mathfrak{F}_p(\mathbb{R}^+)$ generated by the Cauchy singular integral operators with ``shifts''
$$  S^c_{\mathbb{R}^+}\varphi(x):=\frac1{\pi i}\int\limits_0^\infty \frac{e^{-ic(x-y)}\varphi(y)\,dy}{y-x}=W_{-\sign(\xi-c)}\varphi(x)\;\;
        \text{for all}\;\; c\in\mathbb{R}       $$
and by the identity operator $I$ on the semi-axis $\mathbb{R}^+$ over the field of $N\times N$ complex valued matrices coincides with the Banach algebra generated by Fourier convolution operators with piecewise-constant $N\times N$ matrix symbols contains all Fourier convolution $W_a$ and hankel $H_b$ operators with $N\times N$ matrix symbols $($multipliers$)$ $a,b\in PC\mathfrak{M}_p(\overline{\mathbb{R}})$.
\end{proposition}

Let us consider the Banach algebra $\mathfrak{A}_p(\mathbb{R}^+)$ generated by Mellin convolution and Fourier convolution operators, i.e. by the operators
\begin{equation}\label{e2.5}
    \mathbf{A}:=\sum_{j=1}^m\mathfrak{M}^0_{a_j}W_{b_j},
\end{equation}
and there compositions, in the Lebesgue space $\mathbb{L}_p(\mathbb{R}^+)$. Here $\mathfrak{M}^0_{a_j}$ are Mellin convolution operators with continuous $N\times N$ matrix symbols $a_j\in C\mathfrak{M}_p(\dR)$, $W_{b_j}$ are Fourier convolution operators with $N\times N$ matrix symbols $b_j \in C\mathfrak{M}_p(\overline{\mathbb{R}} \setminus\{0\}) :=C\mathfrak{M}_p(\overline{\mathbb{R}}^-\cup\overline{\mathbb{R}}^+)$. The algebra of $N\times N$ matrix $\mathbb{L}_p$-multipliers $C\mathfrak{M}_p(\overline{ \mathbb{R}}\setminus\{0\})$ consists of those piecewise-continuous $N\times N$ matrix multipliers $b\in\mathfrak{M}_p(\mathbb{R})\cap PC(\overline{\mathbb{R}})$ which are continuous on the semi-axis $\mathbb{R}^-$ and $\mathbb{R}^+$ but might have finite jump discontinuities at $0$ and at the infinity.

This and more general algebras were studied in \cite{Du87} and also in earlier works \cite{Du74, Du86, Th85}.
 %
\begin{remark}\label{r4.12}
If in \eqref{e2.5} we admit more general symbols $a_j\in C\mathfrak{M}_p(\ov{\bR})$ which have different limits at the infinity $a_j(-\infty)\not=a_j(+\infty)$, this will not be a generalization.

Indeed, if $a_j\in C\mathfrak{M}_p(\ov{\bR})$ has different limits at the infinity $a_j(-\infty)\not=a_j(+\infty)$ we can represent
\begin{eqnarray*}
a_j(\xi)=a^0_j(\xi) + a_j(-\infty)\frac{1-\coth\pi\left(\frac ip + \xi\right)}2 + a_j(+\infty)\frac{1+\coth\pi\left(\frac ip +
     \xi\right)}2,\\
a^0_j(\pm\infty)=0
\end{eqnarray*}
and the corresponding Mellin operator is written as follows
\begin{eqnarray*}
\mathfrak{M}^0_{a_j}&=&\mathfrak{M}^0_{a^0_j} +  \frac{ a_j(-\infty)}2\left[I-S_{\bR^+}\right]
     +  \frac{ a_j(-\infty)}2\left[I+S_{\bR^+}\right]\\
&=&\mathfrak{M}^0_{a^0_j} +  \frac{ a_j(-\infty)}2\left[I-W_{-\sign}\right] +  \frac{ a_j(-\infty)}2\left[I+W_{-\sign}\right]
\end{eqnarray*}
(see \eqref{e0.4} and \eqref{e1.12}).  Therefore, the discontinuity at the infinity of the symbols of Mellin convolution operators is taken over in Fourier convolution operators and we can even assume in \eqref{e2.5} that $a^0_j(\pm\infty)=0$ for all $j=1,\ldots,m$.
\end{remark}

In order to keep the exposition self-contained and to improve formulations from \cite{Du87}, the results concerning the Banach algebra generated by the operators \eqref{e2.5} are presented here with some modification and the proofs.

Note that the algebra $\mathfrak{A}_p(\mathbb{R}^+)$ is actually a subalgebra of the Banach algebra $\mathfrak{F}_p(\mathbb{R}^+)$ generated by the Fourier convolution operators $W_a$ with piecewise-constant symbols $a(\xi)$ in the space $\mathbb{L}_p(\mathbb{R}^+)$ (cf. Proposition \ref{p3.9}). Let $\mathfrak{S}(\mathbb{L}_p(\mathbb{R}^+))$ denote the ideal of all compact operators in $\mathbb{L}_p(\mathbb{R}^+)$. Since the quotient algebra $\mathfrak{F}_p(\mathbb{R}^+)/\mathfrak{S}(\mathbb{L}_p(\mathbb{R}^+))$ is commutative in the scalar case $N=1$, the following is true.
 %
\begin{corollary}\label{c3.12}
The quotient algebra $\mathfrak{A}_p(\mathbb{R}^+)/\mathfrak{S}(\mathbb{L}_p (\mathbb{R}^+))$ is commutative in the scalar case $N=1$.
\end{corollary}

To describe the symbol of the operator \eqref{e2.5}, consider
the infinite clockwise oriented ``rectangle'' $\mathfrak{R}:=\Gamma_1\cup\Gamma_2^-\cup\Gamma_2^+ \cup\Gamma_3$, where (cf. Figure~1)
\[  \Gamma_1:=\overline{\mathbb{R}}\times\{+\infty\},\;\;\Gamma^\pm_2:=\{\pm\infty\}\times\overline{\mathbb{R}}^+,\;\;
            \Gamma_3:=\overline{\mathbb{R}}\times\{0\}.
\]

\vskip20mm 
\setlength{\unitlength}{0.4mm} \hskip-2mm
\begin{picture}(300,120)
\put(0,-20){\epsfig{file=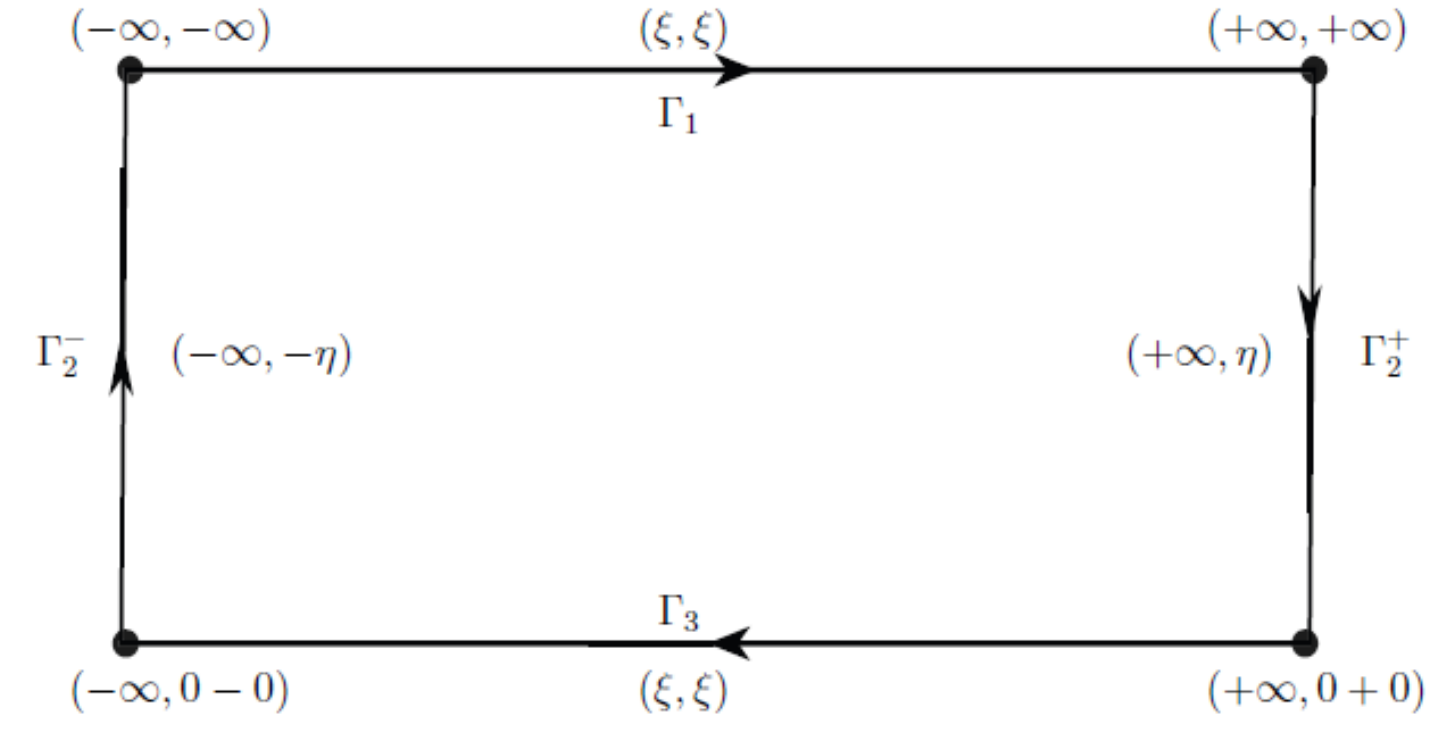,height=55mm, width=130mm}}
\put(180,-35){\makebox(0,0)[bc]{The domain $\mathfrak{R}$ of definition of the symbol $\mathcal{A}_p(\xi, \eta)$.}}
\end{picture}
\vskip20mm
The symbol $\mathcal{A}_p(\omega)$ of the operator $\A$ in \eqref{e2.5} is a function on the set $\mathfrak{R}$, viz.
\begin{equation}\label{e2.7}
    \mathcal{A}_p(\omega):=\begin{cases}
        \displaystyle \sum_{j=1}^m a_j(\xi)(b_j)_p(\infty,\xi), & \omega=(\xi,\infty)\in\overline{\Gamma_1}, \\
        \displaystyle \sum_{j=1}^m a_j(\infty)b_j(\eta), & \omega=(+\infty,\eta)\in\Gamma^+_2, \\
        \displaystyle \sum_{j=1}^m a_j(\infty)b_j(-\eta), & \omega=(-\infty,\eta)\in\Gamma^-_2, \\
        \displaystyle \sum_{j=1}^m a_j(\infty)(b_j)_p(0,\xi), & \omega=(\xi,0)\in\overline{\Gamma_3}.
                \end{cases}
\end{equation}
The symbol $\mathcal{A}_p(\omega)$, when $\omega=(\xi,\infty)$ ranges through the infinite interval $\Gamma_1$ (see Fig. 1) fills  the gap between the values
 \[
\displaystyle \sum_{j=1}^m a_j(\xi)b_j(-\infty)\qquad \text{and}\qquad
     \displaystyle\sum_{j=1}^m a_j(\xi)b_j(+\infty),
 \]
and when $\omega=(\xi,0)$ ranges through the infinite interval $\Gamma_3$ (see Fig. 1) it fills  the gap between the values
 \[
\displaystyle \sum_{j=1}^m a_j(\xi)b_j(0-0)\qquad \text{and}\qquad
     \displaystyle\sum_{j=1}^m a_j(\xi)b_j(0+0).
 \]
The connecting function $g_p(\infty,\xi)$ in \eqref{e2.7} for a piecewise continuous function $g\in PC(\overline{\mathbb{R}})$ is defined as follows
\begin{eqnarray}\label{e2.8}
g_p(x,\xi)&\hskip-3mm:=&\hskip-3mm\frac12\,\big[g(x+0)+g(x-0)\big]-\frac i2\,
     \big[g(x+0)-g(x-0)\big]\cot\pi\Big(\frac1p-i\xi\Big)\nonumber\\
&\hskip-3mm=&\hskip-3mm e^{i\pi\frac{g^+_x+g^-_x}2}\frac{\cos\pi\Big(\frac1p+\frac{g^+_x
     -g^-_x}2-i\xi\Big)}{\sin\pi\Big(\frac1p-i\xi\Big)},\qquad \xi\in\mathbb{R},\\
&&\hskip-10mm g_x^\pm:=\frac1{\pi i}\ln\,g(x\pm0),\quad \text{Re}\,g_x^\pm
     =\frac1\pi\arg\,g(x\pm0),\quad x\in\dR:=\bR\cup\{\infty\}. \nonumber
\end{eqnarray}
The function $g_p(\infty,\xi)$ fills up the discontinuity (the jump) of $g(\xi)$ at $\infty$ between $g(-\infty)$ and $g(+\infty)$ with an oriented arc of the circle such that from every point of the arc the oriented interval $[g(-\infty),g(+\infty)]$ is seen under the angle $\pi/p$. Moreover, the oriented arc lies above the oriented interval if $1/2<1/p<1$ (i.e., if $1<p<2$) and the oriented arc is under the oriented interval if $0<1/p<1/2$ (i.e., if $2<p<\infty$).  For $p=2$ the oriented arc coincides with the oriented interval (see Figure 2).

\vskip15mm 
\setlength{\unitlength}{0.4mm} \hskip-2mm
\begin{picture}(300,140)
\put(22,100){\epsfig{file=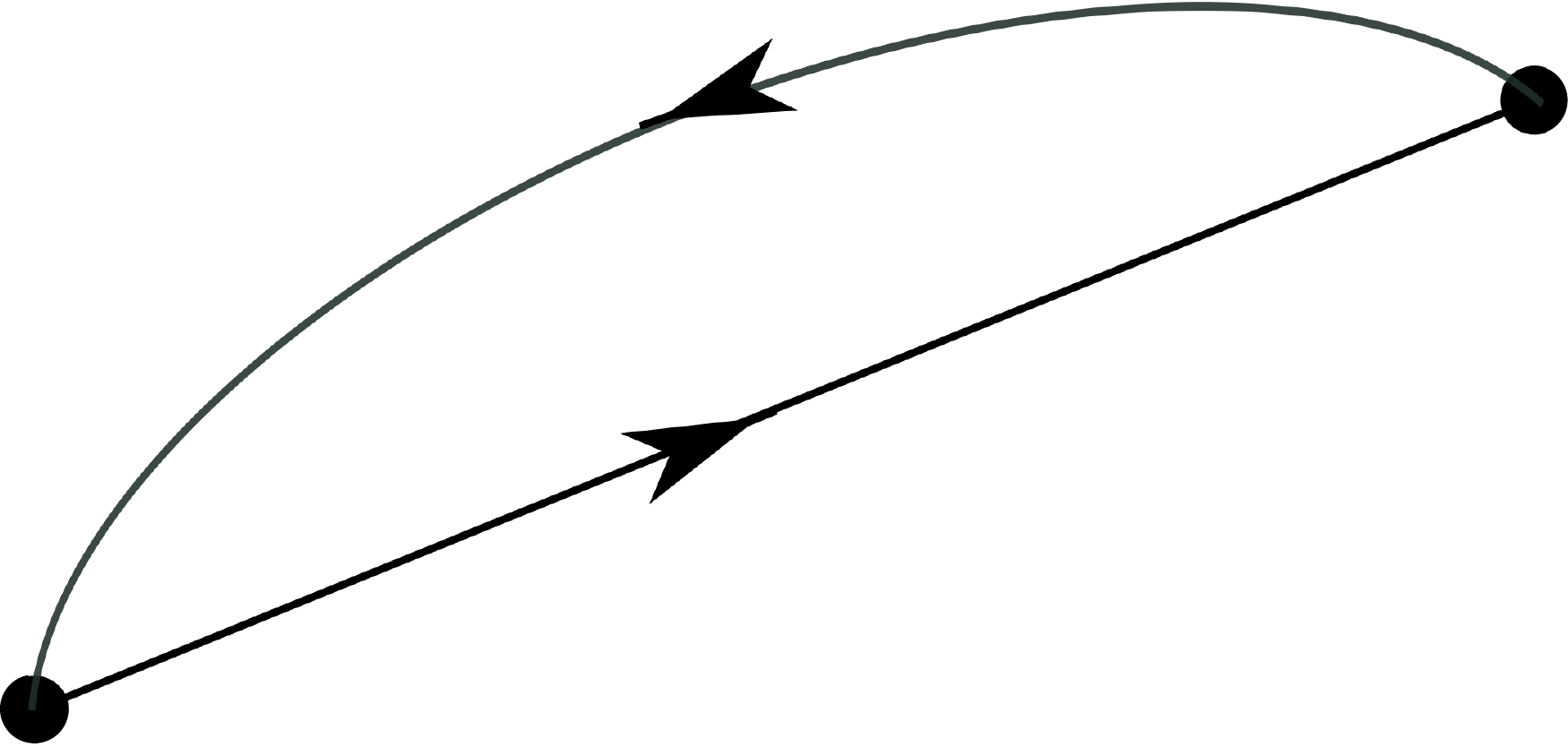,height=20mm, width=35mm }}
\put(155,100){\epsfig{file=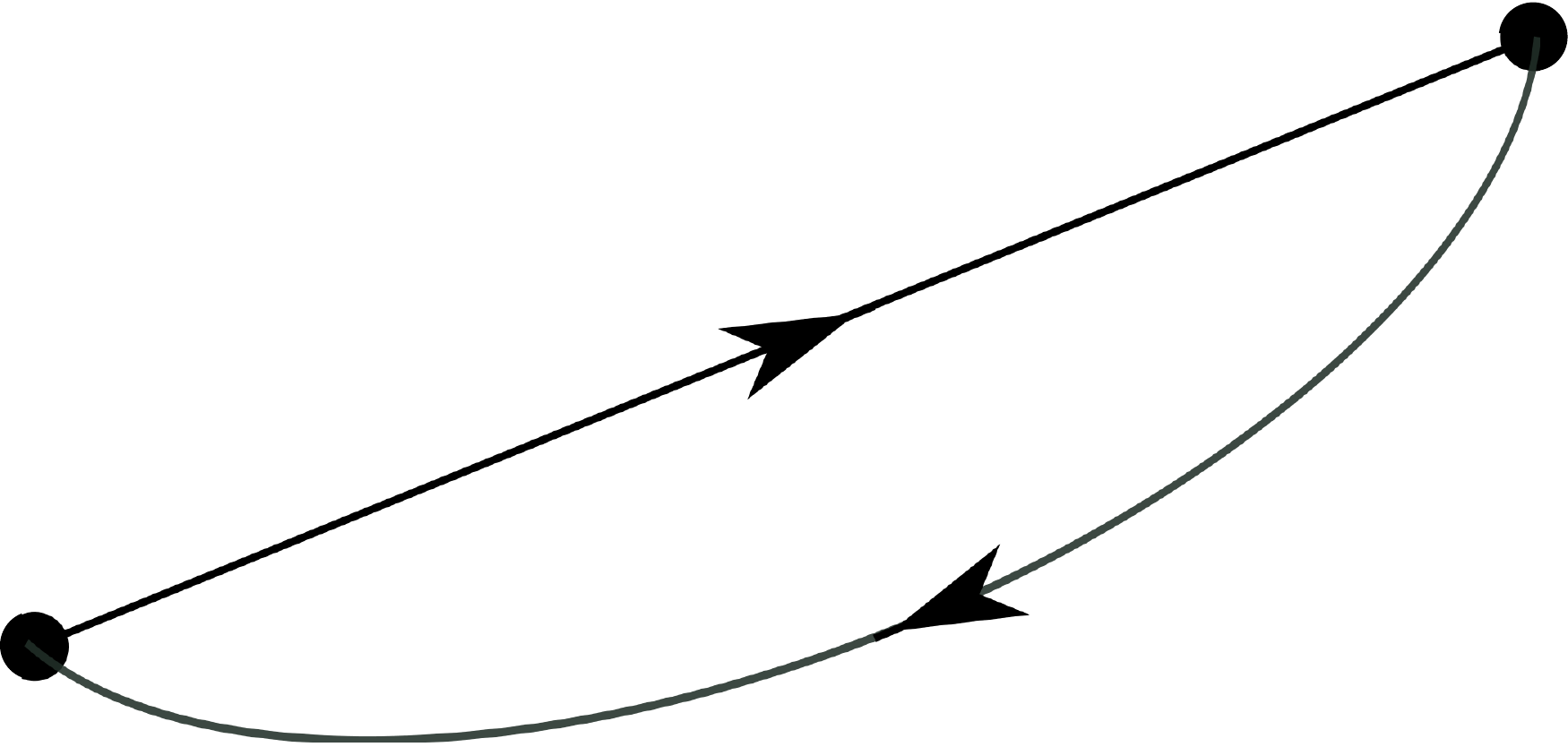,height=20mm, width=35mm }}
\put(265,100){\epsfig{file=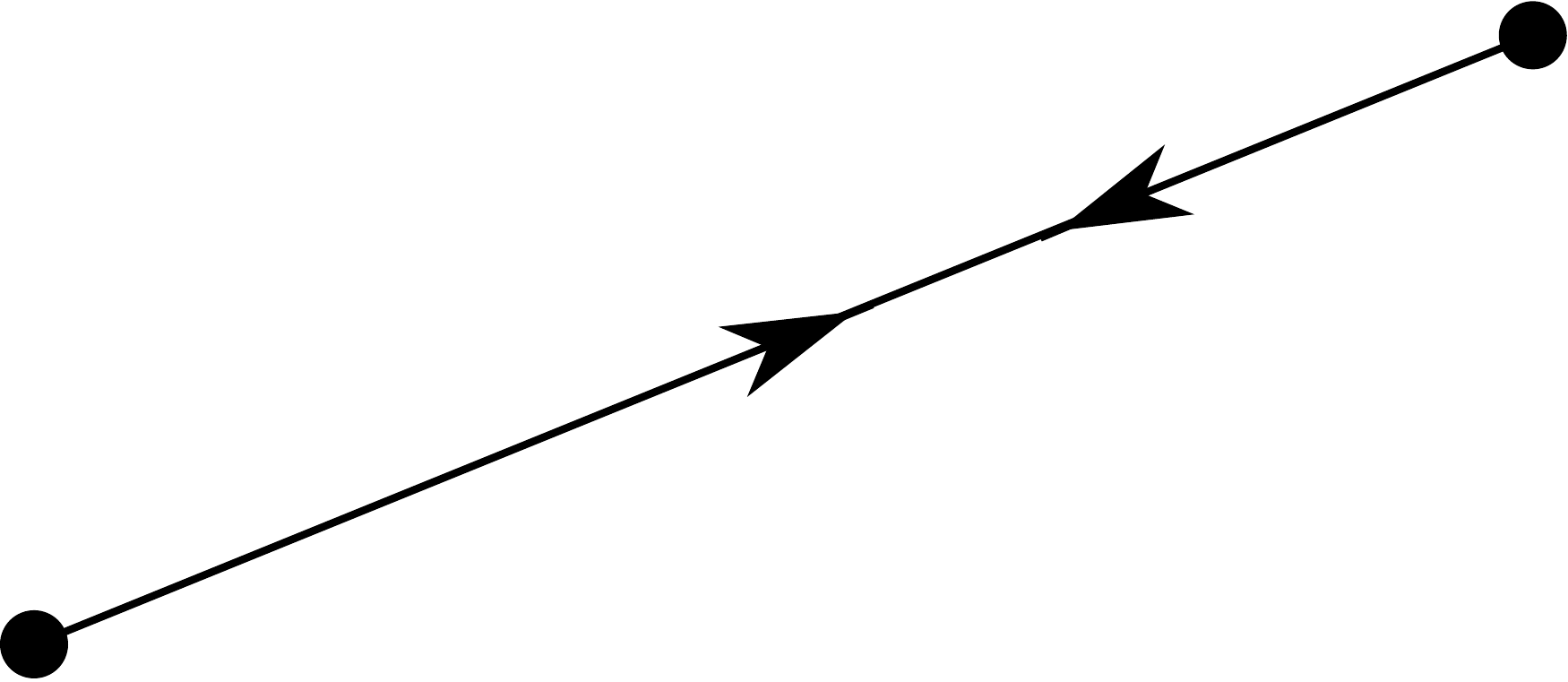,height=20mm, width=35mm }}
\put(07,100){\makebox(0,0)[bc]{$g(-\infty)$}}
\put(137,100){\makebox(0,0)[bc]{$g(-\infty)$}}
\put(247,100){\makebox(0,0)[bc]{$g(-\infty)$}}
\put(125,145){\makebox(0,0)[bc]{$g(+\infty)$}}
\put(220,145){\makebox(0,0)[bc]{$g(+\infty)$}}
\put(330,145){\makebox(0,0)[bc]{$g(+\infty)$}}
\put(45,85){\makebox(0,0)[bc]{$1<p<2$}}
\put(190,85){\makebox(0,0)[bc]{$2<p<\infty$}}
\put(300,85){\makebox(0,0)[bc]{$p=2$}}
\put(180,70){\makebox(0,0)[bc]{Figure 2: Arc condition}}
\end{picture}

\vskip-25mm

A similar geometric interpretation is valid for the function $g_p(t,\xi)$, which connects the point $g(t-0)$ with $g(t+0)$ at the point $t$ where $g(\xi)$ has a jump discontinuity.

To make the symbol $\mathcal{A}_p(\omega)$ continuous, we endow
the rectangle $\mathfrak{R}$ with a special topology. Thus let
us define the distance on the curves $\Gamma_1$, $\Gamma^\pm_2$, $\Gamma_3$ and on $\overline{\mathbb{R}}$ by
\[  \rho(x,y):=\Big|\arg\frac{x-i}{x+i}-\arg\frac{y-i}{y+i}\Big|\;\;\mbox{for arbitrary}\;\; x,y\in\overline{\mathbb{R}}.       \]
In this topology, the length $|\mathfrak{R}|$ of $\mathfrak{R}$ is $6\pi$, and the symbol $\mathcal{A}_p(\omega)$ is continuous everywhere on $\mathfrak{R}$. The image of the function $\det\mathcal{A}_p(\omega)$, $\omega\in \mathfrak{R}$ ($\det\mathcal{B}_p(\omega)$) is a closed curve in the complex plane. It follows from the continuity of  the symbol at the angular points of the rectangle $\mathfrak{R}$ where the one-sided limits coincide. Thus
\begin{align*}
    \mathcal{A}_p(\pm\infty,\infty) & =\sum_{j=1}^m[a_j(\infty)b_j(\mp\infty), \\
    \mathcal{A}_p(\pm\infty,0) & =\sum_{j=1}^m[a_j(\infty)b_j(0\mp0).
\end{align*}
Hence, if the symbol of the corresponding operator is elliptic, i.e. if
\begin{equation}\label{e2.9}
    \inf_{\omega\in\mathfrak{R}} \big|\det\mathcal{A}_p(\omega)\big|>0,
\end{equation}
the increment of the argument $(1/2\pi)\arg \mathcal{A}_p(\omega)$ when $\omega$ ranges through $\mathfrak{R}$ in the positive direction is an integer, is called the winding number or the index and it is denoted by $\ind\det\mathcal{A}_p$.
 %
\begin{theorem}\label{t2.10}
Let $1<p<\infty$ and let $\mathbf{A}$ be defined by \eqref{e2.5}. The operator $\mathbf{A}:\mathbb{L}_p(\mathbb{R}^+)\longrightarrow\mathbb{L}_p(\mathbb{R}^+)$ is Fredholm if and only if its symbol $\mathcal{A}_p(\omega)$ is elliptic. If $\mathbf{A}$ is Fredholm, the index of the operator has the value
\begin{equation}\label{e2.10}
    \Ind\mathbf{A}=-\ind\det\mathcal{A}_p.
\end{equation}
\end{theorem}
{\bf Proof:} Note that our study is based on a localization technique. For more details concerning this approach we refer the reader to \cite{Du79, Du84a, DiSi08, GK79, Si65}.

Let us apply the Gohberg--Krupnik local principle to the operator $\mathbf{A}$ in \eqref{e2.8}, ``freezing'' the symbol of $\mathbf{A}$ at a point $x\in\overline{\mathbb{R}}:=\mathbb{R}\cup \{-\infty\}\cup \{+\infty\}$. For $x\in\mathbb{R}$ and $ \ell\in \mathbb{N}$, $\ell\geq 1$, let $C^\ell_x(\overline{\mathbb{R}})$ denote the set of all $\ell$-times differentiable non-negative functions which are supported in a neighborhood of $x\in\mathbb{R}$ and are identically one everywhere in a smaller neighborhood of $x$. For $x\in\{-\infty\}\cup \{+\infty\}\cup \{\infty\}$, the functions from the corresponding classes $C^\ell_{+\infty}(\overline{\mathbb{R}})$ and $C^\ell_{-\infty}(\overline{\mathbb{R}})$ vanish on semi-infinite intervals $[-\infty,c)$ and $(-c,\infty]$, respectively, for certain $c>0$ and are identically one in smaller neighborhoods. It is easily seen that the system of localizing classes $\{C^\ell_x(\overline{\mathbb{R}})\}_{x\in \overline{\mathbb{R}}}$ is covering in the algebras $C(\overline{\mathbb{R}})$, $\mathfrak{M}_p(\overline{\mathbb{R}})$, respectively (cf. \cite{Du79, Du84a, DiSi08, GK79}).

Let us now consider a system of localizing classes $\{\mathfrak{L}_{\omega,x}\}_{(\omega,x) \in\mathfrak{R}\times\overline{\mathbb{R}}{}^+}$  in the quotient algebra $\mathfrak{A}_p(\mathbb{R}^+)/\mathfrak{S}(\mathbb{L}_p( \mathbb{R}^+))$. These localizing classes depend on two variables, viz. on
$\omega\in\mathfrak{R}$ and $x\in\overline{\mathbb{R}^+}$. In particular, the class $\mathfrak{L}_{\omega,x}$ contains the operator $\Lambda_{\omega,x}$,
\begin{equation}
    \Lambda_{\omega,x}:=\begin{cases}
            \big[h_0\mathfrak{M}^0_{v_\xi}W_{g_\infty}\big]=\big[h_0\mathfrak{M}^0_{v_\xi}\big] \\
            \qquad\qquad \text{if}\;\;\omega=(\xi,\infty)\in\Gamma_1,\;\; x=0; \\
            \big[h_x\mathfrak{M}^0_{v_{\pm\infty}}W_{g_\infty}\big]=\big[h_x\mathfrak{M}^0_{v_{\pm\infty}}W_{g_{\mp\infty}}\big] \\
            \qquad\qquad \text{if}\;\;\omega=(\pm\infty,\infty)\in\Gamma^\pm_2\cap\Gamma_1,\;\;x\in\mathbb{R}^+; \\
            \big[h_\infty\mathfrak{M}^0_{v_{\pm\infty}}W_{g_\eta}\big]=\big[h_\infty\mathfrak{M}^0_{v_{\pm\infty}}W_{g_{\mp\,\eta}}\big] \\
            \qquad\qquad \text{if}\;\;\omega=(\pm\infty,\eta)\in\Gamma^\pm_2,\;\; x=\infty; \\
            \big[h_\infty\mathfrak{M}^0_{v_\infty}W_{g_0}\big]=\big[\mathfrak{M}^0_{v_\xi}W_{g_0}\big] \\
            \qquad\qquad \text{if}\;\;\omega=(\xi,0)\in\overline\Gamma_3,\;\; x=\infty,
        \end{cases} \label{e2.11}
\end{equation}
where $h_x\in C_x^1(\overline{\mathbb{R}}{}^+)$, $v_\xi\in C^1_\xi(\overline{\mathbb{R}}{}^+)$, $g_\eta\in C^1_\eta(\overline{\mathbb{R}}{}^+)$, and $[\mathbf{A}]\in  \linebreak      \mathfrak{A}_p(\mathbb{R}^+)/\mathfrak{S}(\mathbb{L}_p(\mathbb{R}^+))$ denotes the coset containing the operator
$\mathbf{A}\in\mathfrak{A}_p(\mathbb{R}^+)$.

To verify the equalities in \eqref{e2.11}, one has to show that
the difference between the operators in the square brackets is
compact.

Consider the first equality in \eqref{e2.11}: The operator
 \[
h_0W_{g_\infty}-h_0I=h_0W_{(g_\infty-1)} =h_0W_{g_0}
 \]
is compact, since both functions $h_0$ and $1-g_\infty=g_0$ have compact supports, so Proposition \ref{p3.5} applies.

To check the second equality in \eqref{e2.11}, let us note that
$h_x(0)=0$,  $v_{\pm\infty}(\mp\infty)  \linebreak  =0$ and $g_{\pm\infty}(\xi)=0$ for all $\mp\,\xi>0$. From the fourth part of Proposition~\ref{p3.7} we derive that for any $x\in\mathbb{R}^+$ the operator $h_x\mathfrak{M}^0_{v_{\pm\infty}}W_{g_{\pm\infty}}$ is compact. This leads to the claimed equality since
\[  \big[h_x\mathfrak{M}^0_{v_{\pm\infty}}W_{g_\infty}\big]=\big[h_x\mathfrak{M}^0_{v_{\pm\infty}}\{W_{g_{-\infty}}+W_{g_{+\infty}}\}\big]=
            \big[h_x\mathfrak{M}^0_{v_{\pm\infty}}W_{g_{\mp\infty}}\big].       \]
The third identity in \eqref{e2.11} can be verified analogously.
As far as the fourth identity in \eqref{e2.11} is concerned,
one can replace $h_\infty$ by $1$ because the difference
$h_\infty W_{g_0}-W_{g_0} =(1-h_\infty)W_{g_0}=h_0W_{g_0}$ is
compact due to Proposition \ref{p3.5}.

Consider now other properties of the system $\{\mathfrak{L}_{\omega,x}\}_{(\omega,x)\in\mathfrak{R}\times \overline{\mathbb{R}}{}^+}$. Propositions \ref{p3.5}--\ref{p3.8} imply that
$$  \big[h_x\mathfrak{M}^0_{v_\xi}W_{g_\infty}\big]=0\;\;\text{for all} \;\;
            (\xi,\eta,x)\in\overline{\mathbb{R}}\times\overline{\mathbb{R}}\times
                    \overline{\mathbb{R}}{}^+\setminus\mathfrak{R}\times \overline{\mathbb{R}}{}^+.     $$
Therefore, the system of localizing classes $\{\mathfrak{L}_{\omega,x}\}_{(\omega,x)\in\mathfrak{R}\times\overline{\mathbb{R}}{}^+}$ is covering: for a given system $\{\Lambda_{\omega,x}\}_{(\omega,x)\in\mathfrak{R}\times\overline{\mathbb{R}}{}^+}$ of localizing operators one can select a finite number of points $(\omega_1,x_1)=(\xi_1,\eta_1,x_1),\dots,(\omega_s,x_s)=(\xi_s,\eta_s,x_s)\in\mathfrak{R}$ and add appropriately chosen terms $[h_{x_{s+j}}\mathfrak{M}^0_{v_{\xi_{s+j}}}W_{g_{s+j}}]=0$ with     \linebreak     $(\xi_{s+j},\eta_{s+j},x_{s+j}))\in\overline{\mathbb{R}}\times\overline{\mathbb{R}}\times \overline{\mathbb{R}}{}^+\setminus(\mathfrak{R}\times\overline{\mathbb{R}}{}^+)$, $j=1,2,\dots,r$ so, that the equality
\begin{equation}\label{e2.12}
    \sum_{j=1}^r\sum_{k=1}^s \big[c_{x_j}\mathfrak{M}^0_{a_{\xi_j}}W_{b_{\eta_k}}\big]=\big[c\mathfrak{M}^0_aW_b\big]
\end{equation}
holds and the functions $c\in C(\overline{\mathbb{R}}{}^+)$, $a\in C\mathfrak{M}_p(\overline{\mathbb{R}})$, $b\in C\mathfrak{M}_p(\overline{\mathbb{R}})$ are all elliptic. This implies the invertibility of the coset $[c\mathfrak{M}^0_aW_b]$ in the
quotient algebra $\mathfrak{A}_p(\mathbb{R}^+)/\mathfrak{S}(\mathbb{L}_p( \mathbb{R}^+))$ and the inverse coset is $[c\mathfrak{M}^0_aW_b]^{-1}=[c^{-1}\mathfrak{M}^0_{a^{-1}}W_{b^{-1}}]$.

Note that the choice of a finite number of terms in
\eqref{e2.12} is possible due to Borel--Lebesgue lemma and the
compactness of the sets $\overline{\mathbb{R}}$ and $\overline{\mathbb{R}}{}^+$ (two point and one point compactification of $\mathbb{R}$ and of $\mathbb{R}^+$, respectively).

Moreover, localization in the quotient algebra
$\mathfrak{A}_p(\mathbb{R}^+)/\mathfrak{S}(\mathbb{L}_p(\mathbb{R}^+))$ leads to the following local representatives of the cosets containing Mellin and Fourier convolution operators with symbols $a,b\in C\mathfrak{M}_p(\overline{\mathbb{R}})$:
\begin{subequations}
\begin{align}
\label{c2.12a}
    [\mathfrak{M}^0_a] & \overset{\mathfrak{M}^0_{v_{\xi_0}}}{\sim}[\mathfrak{M}^0_{a(\xi_0)}]=[a(\xi_0)I]
            \;\;\text{if}\;\; \xi_0\in\overline{\mathbb{R}}, \\
\label{c2.12b}
    [\mathfrak{M}^0_a] & \overset{v_{x_0}I}{\sim}[\mathfrak{M}^0_{a^\infty}] \;\;\text{if}\;\; \xi_0\in\overline{\mathbb{R}{}^+},\;\; x_0\not=0,\\
\label{c2.12c}
    [\mathfrak{M}^0_a] & \overset{v_0I}{\sim}[\mathfrak{M}^0_a] \;\;\text{if}\;\; \xi_0=0,\\
\label{c2.12d}
    [W_b] & \overset{W_{b_{\eta_0}}}{\sim}[W_{b(\eta_0)}]=[b(\eta_0)I]\;\;\text{if}\;\; \eta_0\in\mathbb{R}\setminus\{0\},\\
\label{c2.12e}
    [W_b] & \overset{W_{b_0}}{\sim}[W_{b^0}]=[\mathfrak{M}^0_{b_p(0,\cdot)}]\;\;\text{if}\;\; \eta=0,\\
\label{c2.12f}
    [W_b] & \overset{W_{g_\infty}}{\sim}[W_{b^\infty(\infty,\cdot)}]=[\mathfrak{M}^0_{b_p(\infty,\cdot)}]\;\;\text{if}\;\; \eta_0=\pm\infty,\\
\label{c2.12g}
    [W_b] & \overset{{v_{x_0}I}}{\sim}[W_{b^\infty}]=[\mathfrak{M}^0_{b_p(\infty,\cdot)}]\;\;\text{if}\;\; x_0\in\mathbb{R}^+,\\
\label{c2.12h}
    [W_b] & \overset{{v_{\infty}I}}{\sim}[W_{b}]\;\;\text{if}\;\; x_0=\infty,
\end{align}
\end{subequations}
where
\begin{equation}\label{e2.13}
\begin{aligned}
    g^\infty(\xi) & :=\frac12\,\big[g(+\infty)+g(-\infty)\big]+\frac12\,\big[g(+\infty)-g(-\infty)\big]\sign\xi= \\
    &\qquad\qquad  =g(-\infty)\chi_-(\xi)+ g(+\infty)\chi_+(\xi), \\
    g^0(\xi) & :=\frac12\,\big[g(0+0)+g(0-0)\big]+\frac12\,\big[g(0+0)+g(0-0)\big]\sign\xi= \\
    &\qquad\qquad =g(0-0)\chi_-(\xi)+g(0+0)\chi_+(\xi),
\end{aligned}
\end{equation}
and $\chi_\pm(\xi):=(1/2)(1\pm\sign\xi)$. Note that in the equivalency relations \eqref{c2.12e}--\eqref{c2.12g} we used the identities, cf. \eqref{e2.3} and \eqref{e2.8},
\begin{align*}
    W_{g^\infty} & =\frac12\,\big[g(-\infty)-g(+\infty)\big]-\frac12\,\big[g(-\infty)-g(+\infty)\big]S_{\mathbb{R}^+}=
                \mathfrak{M}_{g_p(\infty,\cdot)}, \\
    W_{g^0} & =\frac12\,\big[g(0+0)+g(0-0)\big]-\frac12\,\big[g(0+0)-g(0-0)\big]S_{\mathbb{R}^+}=\mathfrak{M}_{g_p(0,\cdot)},
\end{align*}
which means that the Fourier convolution operators with
homogeneous of order $0$ symbols  $g^\infty(\xi)$ and $g^0(\xi)$
are, simultaneously, Mellin convolutions with the symbols
$g_p(\infty,\xi)$, $g_p(0,\xi)$.

Using the equivalence relations \eqref{c2.12a}--\eqref{c2.12h}
and the compactness of the corresponding operators, cf.
Propositions \ref{p3.5}--\ref{p3.7}, one finds easily the
following local representatives of the operator (coset)
$\mathbf{A}\in\mathfrak{A}_p(\mathbb{R}^+)/\mathfrak{S} \mathbb{L}_p(\mathbb{R}^+)$ (see \eqref{e2.8} for the operator $\mathbf{A}$):
\begin{subequations}
\begin{align}
    [\mathbf{A}] & \overset{\Lambda_{(\xi_0,\infty),0}}{\sim}\Big[\sum_{j=1}^m \mathfrak{M}^0_{a_j(\xi_0)}W_{(b_j)^\infty}\Big]= \nonumber \\
    &\qquad =\Big[\sum_{j=1}^m \mathfrak{M}^0_{a_j(\xi_0)(b_j)_p(\infty,\cdot)}\Big]\overset{\Lambda_{(\xi_0,\infty),0}}{\sim}
                    \Big[\sum_{j=1}^m\mathfrak{M}^0_{a_j(\xi_0)(b_j)_p(\infty,\xi_0)}\Big]= \nonumber \\
    &\qquad =\big[\mathcal{A}_p(\xi_0,\infty)I\big]\;\;\text{if}\;\;\omega=(\xi_0,\infty)\in\Gamma_1,\;\; x_0=0, \label{e2.14a} \\
    [\mathbf{A}] & \overset{\Lambda_{(\pm\infty,\infty),x_0}}{\sim}\Big[\sum_{j=1}^m\mathfrak{M}^0_{a_j(\pm\infty)}W_{(b_j)^\infty}\Big]=
            \Big[\sum_{j=1}^m\mathfrak{M}^0_{a_j(\pm\infty)(b_j)_p(\infty,\cdot)}\Big]= \nonumber \\
    &\qquad =\big[\mathfrak{M}^0_{\mathcal{A}_p(\pm\infty,\cdot)}\big]\overset{\Lambda_{(\pm\infty,\infty),x_0}}{\sim}
        \big[\mathcal{A}_p(\pm\infty,\infty)I\big] \label{e2.14b}\\
    &\qquad\qquad\qquad\qquad \text{if}\;\;\omega=(\pm\infty,\infty)\in\overline{\Gamma^\pm_2}\cap \overline{\Gamma_1},\;\;0<x_0<\infty; \nonumber \\
    [\mathbf{A}] & \overset{\Lambda_{(\pm\infty,\mp\,\eta_0),\infty}}{\sim}
        \Big[\sum_{j=1}^m\mathfrak{M}^0_{a_j(\pm\infty)}W_{b_j(\mp\,\eta_0)}\Big]=\Big[\sum_{j=1}^m a_j(\pm\infty)b_j(\mp\,\eta_0)I\Big]= \nonumber \\
    &\qquad =\!\big[\mathcal{A}_p(\pm\infty,\mp\,\eta_0)I\big] \;\;\text{if}\;\;\eta_0\!>\!\!0,\;
        \omega\!=\!(\pm\infty,\mp\,\eta_0)\!\in\!\!\Gamma^\pm_2, \; x_0\!\!=\!\infty; \label{e2.14c} \\
    [\mathbf{A}] & \overset{\Lambda_{(\xi_0,0),\infty}}{\sim}\Big[\sum_{j=1}^m\mathfrak{M}^0_{a_j}W_{b^0_j}\Big]= \nonumber \\
    &\qquad =\Big[\sum_{j=1}^m a_j(\xi_0)\mathfrak{M}_{(b_j)_p(0,\cdot)}\Big]\overset{\Lambda_{(\xi_0,0),\infty}}{\sim}
            \Big[\sum_{j=1}^m a_j(\xi_0)(b_j)_p(0,\xi_0)\Big]= \nonumber\\
    &\qquad =\big[\mathcal{A}_p(\xi_0,0)I\big]\;\;\text{if}\;\;\omega=(\xi_0,0)\in\overline{\Gamma}_3, \;\; x_0=\infty; \label{e2.14d} \\
    [\mathbf{A}] & \overset{\Lambda_{(\pm\infty,\eta),\infty}}{\sim}\Big[
    \sum_{j=1}^m\mathfrak{M}^0_{a_j(\pm\infty)}W_{b_j(0)}\Big]=
        \Big[\sum_{j=1}^m a_j(\pm\infty)b_j(0)I\Big]= \nonumber \\
    &\qquad =\big[\mathcal{A}_p(\pm\infty,0)I\big]\;\;\text{if}\;\;\omega=(\pm\infty,0)\in\overline{\Gamma}_3,\;\; x_0=\infty. \label{e2.14e}
\end{align}
\end{subequations}
It is remarkable that the local representatives \eqref{e2.14a}--\eqref{e2.14e} are just the quotient classes of multiplication operators by constant $N\times N$ matrices
$\left[\mathcal{A}_p(\xi_0,\eta_0)I\right]$. If $\det\mathcal{A}_p(\xi_0,\eta_0)=0$, these representatives are not invertible, both locally and globally. On the other hand, they are globally invertible if $\det\mathcal{A}_p(\xi_0,\eta_0)\not=0$. Thus, the conditions of the local invertibility for all points $\omega_0=(\xi_0,\eta_0)\in\mathfrak{R}$ and the global invertibility of the operators under consideration coincide with the ellipticity condition for the symbol $\inf\limits_{(\xi_0,\eta_0)\in\mathfrak{R}}\det\mathcal{A}_p(\xi_0,\eta_0)\not=0$.

The index $\Ind\mathbf{A}$ is a continuous integer-valued
multiplicative function $\Ind\mathbf{A}\mathbf{B}=\Ind\mathbf{A} + \Ind\mathbf{B}$ defined on the group of Fredholm operators of $\mathfrak{A}_p(\mathbb{R}^+)$. On the other hand, the index function $\ind\det\mathcal{A}_p$ defined on $L_p$-symbols $\mathcal{A}_p$ possesses the same property $\ind\det\mathcal{A}_p\mathcal{B}_p=\ind\det\mathcal{A}_p+\ind\det\mathcal{B}_p$, see explanations after \eqref{e2.9}. Moreover, the set of operators \eqref{e2.8} is dense in the algebra $\mathfrak{A}_p(\mathbb{R}^+)$ and the corresponding set of their symbols is dense in the algebra $C(\mathfrak{R})$ of all continuous functions on $\mathfrak{R}$. For $p=2$ these algebras even coincide. Therefore, there is an algebraic homeomorphism between the quotient algebra $\mathfrak{A}_p(\mathbb{R}^+)/\mathfrak{S}(\mathbb{L}_p(\mathbb{R}^+))$ and the algebra of their symbols which is a dense subalgebra of $C(\mathfrak{R})$. Hence, two various index functions can be only connected by the relation $\Ind\mathbf{A}=M_0\ind\det\mathcal{A}_p$ with an integer constant $M_0$ independent of $\mathbf{A}\in \mathfrak{A}_p(\mathbb{R}^+)/\mathfrak{S}(\mathbb{L}_p(\mathbb{R}^+))$. Since for any Fourier convolution operator $\mathbf{A}=W_a$ the index formula is $\Ind\mathbf{A}=-\ind\det\mathcal{A}_p$ \cite{Du74, Du75a, Du79}, the constant $M_0=-1$, and the index formula \eqref{e2.10} is proved.
\QED

While proving the foregoing Theorem \ref{t2.10} with the help of local principle, as a byproduct we have proved the following
 %
\begin{corollary}\label{c3.14}
Let $1<p<\infty$, $s\in\mathbb{R}$ and $\mathbf{A}$ be defined by \eqref{e2.5}. The operator
 \begin{equation}\label{e19b}
\mathbf{A}:\widetilde{\mathbb{H}}{}^s_p(\mathbb{R}^+)\longrightarrow
     \mathbb{H}^s_p (\mathbb{R}^+)
 \end{equation}
is locally invertible at $0\in\mathbb{R}^+$ if and only if its symbol $\mathcal{A}^s_p(\omega)$, defined in \eqref{e2.7}--\eqref{e2.8}, is elliptic only on $\Gamma_1$:
 \[
\inf_{\omega\in\Gamma_1}\left|\det\,\mathcal{A}^s_p(\omega)\right|
     =\inf_{\xi\in\mathbb{R}}\left|\det\,\mathcal{A}^s_p(\xi,\infty)\right|>0.
 \]
\end{corollary}
 %
\begin{remark}\label{r3.14}
Let us emphasize that the formula \eqref{e2.10} does not contradict the invertibility of ``pure Mellin convolution'' operators $\mathfrak{M}_a^0:\mathbb{L}_p (\mathbb{R}^+)\longrightarrow\mathbb{L}_p(\mathbb{R}^+)$ with an elliptic matrix symbol $a\in C\mathfrak{M}_p^0(\mathbb{R})$, $\inf\limits_{\xi\in \mathbb{R}}| a(\xi)|>0$, stated in Proposition \ref{p1.3}, even if $\ind a\not=0$.

In fact, computing the symbol of $\mathfrak{M}_a^0$ by formula \eqref{e2.7}, one obtains
\[  (\mathfrak{M}_a^0)_p(\omega):=\begin{cases}
                a(\xi), & \omega=(\xi,\infty)\in\overline{\Gamma_1}, \\
                a(+\infty), & \omega=(+\infty,\eta)\in\Gamma^+_2, \\
                a(-\infty), & \omega=(-\infty,\eta)\in\Gamma^-_2, \\
                a(\xi), & \omega=(\xi,0)\in\overline{\Gamma_3}.
                        \end{cases}     \]
Noting that on the sets $\Gamma_1$ and $\Gamma_3$ the variable
$\omega$ runs in opposite direction, the increment of the
argument $[\arg \det(\mathfrak{M}_a^0)_p(\omega)]_\mathfrak{R}=0$
is zero, implying $\Ind\mathfrak{M}^0_a=0$.

In contrast to the above, the pure Fourier convolution operators \linebreak
$W_b:\mathbb{L}_p(\mathbb{R}^+)\longrightarrow\mathbb{L}_p(\mathbb{R}^+)$ with elliptic matrix symbol $b\in C\mathfrak{M}_p^0(\mathbb{R})$, \linebreak $\inf\limits_{\xi\in\mathbb{R}}|b_p(\xi,\eta)|>0$ can possess non-zero indices. Since
\[  b_p(\omega):=\begin{cases}
                b_p(\infty,\xi), & \omega=(\xi,\infty)\in\overline{\Gamma_1},\\
                b(-\eta), &\omega=(+\infty,\eta)\in\Gamma^+_2,\\
                b(\eta), &\omega=(-\infty,\eta)\in\Gamma^-_2,\\
                b(0), &\omega=(\xi,0)\in\overline{\Gamma_3},
                        \end{cases}             \]
one arrives at the well-known formula
$$  \Ind W_b=-\ind b_p.       $$
Moreover, in the case where the symbol $b(-\infty)=b(+\infty)$
is continuous, one has $b_p(\xi,\eta)=b(\xi)$. Thus the
ellipticity of the corresponding operator leads to the formula
$$  \ind b_p=\ind\det b.        $$
\end{remark}

If $\mathcal{A}_p(\omega)$ is the symbol of an operator $\mathbf{A}$ of \eqref{e2.5}, the set $\mathcal{R}(\mathcal{A}_p):=\{\mathcal{A}_p(\omega)\in\mathbb{C}:\omega\in \mathfrak{R}\}$ coincides with the essential spectrum of $\mathbf{A}$. Recall that the essential spectrum $\sigma_{ess}(\mathbf{A})$ of a bounded operator $\mathbf{A}$ is the set of all $\lambda\in \mathbb{C}$ such that the operator $\mathbf{A}-\lambda I$ is not Fredholm in $\mathbb{L}_p(\mathbb{R}^+)$ or, equivalently, the coset $[\mathbf{A}-\lambda I]$ is not invertible in the quotient algebra $\mathfrak{A}_p(\mathbb{R}^+)/ \mathfrak{S}(\mathbb{L}_p(\mathbb{R}^+))$. Then, due to Banach theorem, the essential norm $\|\!|\mathbf{A}\|\!|$ of the operator $\mathbf{A}$ can be estimated as follows
\begin{equation}\label{e2.15}
\sup_{\omega\in\omega}|\mathcal{A}_p(\omega)|\leqslant\|\!|\mathbf{A}\|\!|
     :=\inf_{\mathbf{T}\in\mathfrak{S}(\mathbb{L}_p(\mathbb{R}^+))}
     \big\|(\mathbf{A}+\mathbf{T})\mid\mathcal{L}(\mathbb{L}_p(\mathbb{R}^+))\big\|.
\end{equation}
The inequality \eqref{e2.15} enables one to extend continuously
the symbol map \eqref{e2.7}
\begin{equation}\label{e2.16}
[\mathbf{A}]\longrightarrow\mathcal{A}_p(\omega),\quad [\mathbf{A}]\in\mathfrak{A}_p(
     \mathbb{R}^+)/\mathfrak{S}(\mathbb{L}_p(\mathbb{R}^+))
\end{equation}
on the whole Banach algebra $\mathfrak{A}_p(\mathbb{R}^+)$. Now, applying Theorem \ref{t2.10} and atandard methods, cf. \cite[Theorem 3.2]{Du87}, one can derive the following result.
 %
\begin{corollary}\label{c3.15}
Let $1<p<\infty$ and $\mathbf{A}\in\mathfrak{A}_p(\mathbb{R}^+)$. The operator $\mathbf{A}:\mathbb{L}_p(\mathbb{R}^+)\longrightarrow\mathbb{L}_p(\mathbb{R}^+)$ is Fredholm if and only if it's symbol $\mathcal{A}_p(\omega)$ is elliptic. If $\mathbf{A}$ is Fredholm, then
$$  \Ind\mathbf{A}=-\ind\mathcal{A}_p.      $$
\end{corollary}
 %
\begin{corollary}\label{c3.13}
The set of maximal ideals of the commutative Banach quotient
algebra $\mathfrak{A}_p(\mathbb{R}^+)/\mathfrak{S}(\mathbb{L}_p(\mathbb{R}^+))$ generated by scalar $N=1$ operators in \eqref{e2.5}, is homeomorphic to $\mathfrak{R}$, and the symbol map in \eqref{e2.7}, \eqref{e2.16} is a Gelfand homeomorphism of the corresponding Banach algebras.
\end{corollary}
{\bf Proof:}
The proof is based on Theorem \ref{t2.10} and Corollary \ref{c3.15} and if similar to  \cite[Theorem 3.1]{Du87}. The details of the proof is left to the reader.
\QED
%
\begin{remark}\label{r3.17}
All the above results are valid in a more general setting viz.,
for the Banach algebra $\mathfrak{PA}^{N\times N}_{p,\alpha}(\mathbb{R}^+)$ generated in the weighted Lebesgue space of $N$-vector-functions $\mathbb{L}^N_p(\mathbb{R}^+, x^\alpha)$ by the operators
\begin{equation}\label{e2.17}
    \mathbf{A}:=\sum_{j=1}^m \Big[d^1_j\mathfrak{M}^0_{a^1_j}W_{b^1_j}+d^2_j\mathfrak{M}^0_{a^2_j}H_{c^1_j}+d^3_jW^0_{b^2_j}H_{c^2_j}\Big]
\end{equation}
when coefficients $d^1_j, d^2_j, d^3_j\in PC^{N\times N}(\overline{\mathbb{R}})$ are piecewise-continuous $N\times N$ matrix functions, symbols of Mellin convolution operators $\mathfrak{M}^0_{a^1_j}$, $\mathfrak{M}^0_{a^2_j}$, Winer--Hopf (Fourier convolution) operators $W_{b^1_j}$, $W_{b^2_j}$ and Hankel operators $H_{c^1_j}$, $H_{c^2_j}$ are $N\times N$ piecewise-continuous matrix  $\mathbb{L}_p$-multipliers $a^k_j,b^k_j,c^k_j\in PC^{N\times N}\mathfrak{M}_p( \mathbb{R})$.

The spectral set $\Sigma(\mathfrak{PA}^{N\times N}_{p,\alpha}(\mathbb{R}^+))$ of such Banach algebra (viz., the set where the symbols are defined, e.g. $\mathfrak{R}$ for the Banach algebra $\mathfrak{A}^{N\times N}_p(\mathbb{R}^+)$ investigated above) is more sophisticated and described in the papers \cite{Du77, Du78, Du87, Th85}. Let $\mathfrak{CA}_{p,\alpha}(\mathbb{R}^+)\mathfrak{S}(\mathbb{L}_p(\mathbb{R}^+))$ be the sub-algebra of $\mathfrak{PA}_{p,\alpha}(\mathbb{R}^+)=\mathfrak{PA}^{1 \times 1}_{p,\alpha}(\mathbb{R}^+)$ generated by scalar operators \eqref{e2.17} with continuous coefficients $c_j,h_j\in C(\overline{\mathbb{R}})$ and scalar piecewise-continuous $\mathbb{L}_p$-multipliers) $a_j,b_j,d_j,g_j\in PC\mathfrak{M}_p(\mathbb{R})$. The quotient-algebra \\ $\mathfrak{CA}_{p,\alpha}(\mathbb{R}^+)\mathfrak{S}(\mathbb{L}_p(\mathbb{R}^+))$ with respect to the ideal of all compact operators is a commutative algebra and the spectral set $\Sigma(\mathfrak{PA}_{p,\alpha} (\mathbb{R}^+))$ is homeomorphic to the set of maximal ideals.

We drop further details about the Banach algebra  $\mathfrak{PA}^{N\times N}_{p,\alpha}(\mathbb{R}^+)$, because the result formulated above are sufficient for the purpose of this and subsequent papers dealing with the BVPs in domains with corners at the boundary.
\end{remark}

\section{Mellin Convolution Operators in Bessel Potential Spaces. The Boundednes and Lifting}
\label{s5}

As it was already mentioned, the primary aim of the present paper is to study Mellin convolution operators $\mathfrak{M}^0_a$ acting in Bessel potential spaces,
\begin{equation}\label{e5.1}
    \mathfrak{M}^0_a:\widetilde{\mathbb{H}}^s_p(\mathbb{R}^+)\longrightarrow\mathbb{H}^s_p(\mathbb{R}^+).
\end{equation}
The symbols of these operators are $N\times N$ matrix functions $a\in C\mathfrak{M}^0_p (\overline{\mathbb{R}}),$ continuous on the real axis $\mathbb{R}$ with the only possible jump at infinity.
 %
\begin{theorem}\label{t4.1}
Let $s\in\mathbb{R}$ and $1<p<\infty$.

If conditions of Theorem \ref{t1.6} hold, the Mellin convolution operator between Bessel potential spaces
\begin{eqnarray}\label{e65}
\mathbf{K}^1_c\;:\;\widetilde{\mathbb{H}}{}^r_p(\mathbb{R}^+)\to\mathbb{H}^r_p(\mathbb{R}^+)
\end{eqnarray}
is lifted to the equivalent operator
\begin{eqnarray}\label{e66}
\mathbf{\Lambda}^s_{-\gamma}\mathbf{K}^1_c\mathbf{\Lambda}^{-s}_{\gamma}
     =c^{-s}\mathbf{K}^1_cW_{g^s_{-c\gamma,\gamma}}
     \;:\;\mathbb{L}_p(\mathbb{R}^+)\to\mathbb{L}_p(\mathbb{R}^+),
\end{eqnarray}
where $c^{-s}=|c|^{-s}e^{-is\arg\,c}$ and the function $g^s_{-c\gamma,\gamma}$ is defined in \eqref{e1.11}. 

If conditions of Corollary \ref{c3.2} hold, the Mellin convolution operator between Bessel potential spaces \eqref{e65} is lifted to the equivalent operator
\begin{eqnarray}\label{e67}
\mathbf{\Lambda}^s_{-\gamma}\mathbf{K}^1_c\mathbf{\Lambda}^{-s}_{\gamma}
     =c^{-s}W_{g^s_{-\gamma,-\gamma_0}}\mathbf{K}^1_cW_{g^s_{-c\gamma_0,\gamma}}=c^{-s}
     \mathbf{K}^1_cW_{g^s_{-\gamma,-\gamma_0}g^s_{-c\gamma_0,\gamma}}+\mathbf{T}\nonumber\\
\;:\;\mathbb{L}_p(\mathbb{R}^+)\to\mathbb{L}_p(\mathbb{R}^+),
\end{eqnarray}
where $\mathbf{T}\;:\;\mathbb{L}_p(\mathbb{R}^+) \to\mathbb{L}_p(\mathbb{R}^+)$ is a compact operator.
\end{theorem}
{\bf Proof:} The equivalent operator after lifting is
 \[
\mathbf{\Lambda}^s_{-\gamma}\mathbf{K}^1_c\mathbf{\Lambda}^{-s}_{\gamma}
     \;:\;\mathbb{L}_p(\mathbb{R}^+)\to\mathbb{L}_p(\mathbb{R}^+)
 \]
(see Theorem \ref{t1.1}). To proceed we need two formulae
 \begin{eqnarray}\label{e68}
\mathbf{\Lambda}^s_{-c\gamma}\mathbf{\Lambda}^{-s}_{\gamma}=W_{g^s_{-c\gamma,\gamma}},
     \qquad W_{g^s_{-\gamma,-\gamma_0}}W_{g^s_{-c\gamma_0,\gamma}}
     =W_{g^s_{-\gamma,-\gamma_0}g^s_{-c\gamma_0,\gamma}}.
 \end{eqnarray}
The first one holds because $0<\arg\,\gamma<\pi$ (see \eqref{e1.25a}) and the second one holds because $g^s_{-\gamma,-\gamma_0}(\xi)$ has a smooth, uniformly bounded analytic extension in the complex lower half plane (see \eqref{e28}).

If conditions of Theorem \ref{t1.6} hold, we apply formula \eqref{e54}, the first formula in \eqref{e68} and derive the equality in \eqref{e66}:
 \[
\mathbf{\Lambda}^s_{-\gamma}\mathbf{K}^1_c\mathbf{\Lambda}^{-s}_{\gamma}
     =c^{-s}\mathbf{K}^1_c\mathbf{\Lambda}^s_{-c\gamma}\mathbf{\Lambda}^{-s}_{\gamma}
     =c^{-s}\mathbf{K}^1_cW_{g^s_{-c\gamma,\gamma}}.
 \]

If conditions of Corollary \ref{c3.2} hold, we apply formulae \eqref{eqn36}, \eqref{eqn38}, both formula in \eqref{e68} and derive the equality in \eqref{e67}:
\begin{eqnarray*}
\mathbf{\Lambda}^s_{-\gamma}\mathbf{K}^1_c\mathbf{\Lambda}^{-s}_{\gamma}
     =c^{-s}W_{g^s_{-\gamma,-\gamma_0}}\mathbf{K}^1_c\mathbf{\Lambda}^s_{-c\gamma}
     \mathbf{\Lambda}^{-s}_{\gamma}=c^{-s}W_{g^s_{-\gamma,-\gamma_0}}\mathbf{K}^1_c
     W_{g^s_{-c\gamma_0,\gamma}}\nonumber\\[2mm]
=c^{-s}\mathbf{K}^1_cW_{g^s_{-\gamma,-\gamma_0}}W_{g^s_{-c\gamma_0,\gamma}}+\mathbf{T}.
\end{eqnarray*}
\vskip-7mm
\QED
\vskip5mm
 %
 \begin{remark}
The case of operator $\mathbf{K}^1_1$ is not covered by the foregoing Theorem \ref{t4.1}, where $\arg\,c\not=0$. This case is essentially different as underlined in Theorem \ref{t3.4} because $\mathbf{K}^1_1$ is a Hilbert transform $\mathbf{K}^1_1=-\pi iS_{\mathbb{R}^+}=\pi iW_{\sign}$ and $\mathbf{K}^1_1$ between Bessel potential spaces \eqref{e65} is lifted to the equivalent Fourier convolution operator
\begin{eqnarray}\label{e67a}
\mathbf{\Lambda}^s_{-\gamma}\mathbf{K}^1_1\mathbf{\Lambda}^{-s}_{\gamma}
     =W_{\pi ig^s_{-\gamma,\gamma}\sign}
     \;:\;\mathbb{L}_p(\mathbb{R}^+)\to\mathbb{L}_p(\mathbb{R}^+),
\end{eqnarray}
as it follows from Theorem \ref{t3.4}.
 \end{remark}
 %
\begin{theorem}\label{t4.2}
Let $c_j, d_j\in\mathbb{C}$, $-\pi\leqslant\arg\,c_j<\pi$ $\arg\,c_j\not=0$, for $j=1,\ldots,n$, $0<\arg \gamma<\pi$, $-\pi<\arg(c_j\gamma)<0$ for $j=1,\ldots,m$  and $0<\arg(c_j\gamma)<\pi$ for $j=m+1,\ldots,n$.

The Mellin convolution operator between Bessel potential spaces
\begin{eqnarray}\label{e69}
\A=\sum_{j=1}^nd_j\mathbf{K}^1_{c_j}\;:\;\widetilde{\mathbb{H}}{}^r_p(\mathbb{R}^+)
     \to\mathbb{H}^r_p(\mathbb{R}^+)
\end{eqnarray}
is lifted to the equivalent operator
\begin{subequations}
\begin{eqnarray}\label{e69a}
&&\hskip-15mm\mathbf{\Lambda}^s_{-\gamma}\A\mathbf{\Lambda}^{-s}_{\gamma}
     =\sum_{j=0}^md_jc^{-s}_j\mathbf{K}^1_{c_j}W_{g^s_{-c_j\gamma,-\gamma}}
     +\sum_{j=m+1}^nd_jc^{-s}_jW_{g^s_{-\gamma,-\gamma_j}}
     \mathbf{K}^1_{c_j}W_{g^s_{-c_j\gamma_j,\gamma}}\\
\label{e70}
&&\hskip-8mm=\sum_{j=0}^md_jc^{-s}_j\mathbf{K}^1_{c_j}W_{g^s_{-c_j\gamma,\gamma}}
     +\sum_{j=m+1}^nd_jc^{-s}_j\mathbf{K}^1_{c_j}W_{g^s_{-\gamma,-\gamma_j}
     g^s_{-c_j\gamma_j,\gamma}}+\mathbf{T},
\end{eqnarray}
\end{subequations}
in the $\mathbb{L}_p(\mathbb{R}^+)$ space, where $c^{-s}=|c|^{-s}e^{-is\arg\,c}$ and $\gamma_j$ are such that $0<\arg \gamma_j<\pi$, $-\pi<\arg(c_j\,\gamma_j)<0$ for $j=m+1,\ldots,n$. $\mathbf{T}\;:\;\mathbb{L}_p(\mathbb{R}^+) \to\mathbb{L}_p(\mathbb{R}^+)$ is a compact operator.
\end{theorem}
{\bf Proof:} The proof is a direct consequence of Theorem \ref{t4.1}. \QED
 %
\begin{theorem}\label{t4.3}
Let $s\in\mathbb{R}$ and $1<p<\infty$.

If conditions of Theorem \ref{t1.6} hold, the Mellin convolution operator between Bessel potential spaces
\begin{eqnarray}\label{e71}
\mathbf{K}^2_c\;:\;\widetilde{\mathbb{H}}{}^r_p(\mathbb{R}^+)\to\mathbb{H}^r_p(\mathbb{R}^+)
\end{eqnarray}
is lifted to the equivalent operator
\begin{eqnarray}\label{e72}
\mathbf{\Lambda}^s_{-\gamma}\mathbf{K}^2_c\mathbf{\Lambda}^{-s}_{\gamma}
     =c^{-s}\left[\mathbf{K}^2_c-sc^{-1}\mathbf{K}^1_c\right]W_{g^s_{-c\gamma,\gamma}}
     +s\,\gamma\,c^{-s}\mathbf{K}^1_cW_{g^{s-1}_{-c\gamma,\gamma}}\mathbf{\Lambda}^{-1}_{\gamma}
\end{eqnarray}
in $\mathbb{L}_p(\mathbb{R}^+)$ space, where $c^{-s}=|c|^{-s}e^{-is\arg\,c}$ and the function $g^s_{-c\gamma,\gamma}$ is defined in \eqref{e1.11} and the last summand in \eqref{e72}
 \begin{eqnarray}\label{e73}
\mathbf{T}:=s\,\gamma\,c^{-s}\mathbf{K}^1_cW_{g^{s-1}_{-c\gamma,\gamma}}\mathbf{\Lambda}^{-1}_{\gamma}
     \;:\;\mathbb{L}_p(\mathbb{R}^+)\to\mathbb{L}_p(\mathbb{R}^+),
 \end{eqnarray}
is a compact operator.

If conditions of Corollary \ref{c3.2} hold, the Mellin convolution operator between Bessel potential spaces \eqref{e71} is lifted to the equivalent operator
\begin{eqnarray}\label{e74}
&&\hskip-15mm\mathbf{\Lambda}^s_{-\gamma}\mathbf{K}^2_c\mathbf{\Lambda}^{-s}_{\gamma}
     =c^{-s}W_{g^s_{-\gamma,-\gamma_0}}\left[\mathbf{K}^2_c-sc^{-1}\mathbf{K}^1_c\right]
     W_{g^s_{-c\gamma_0,\gamma}}\nonumber\\[2mm]
&&+s\,\gamma\,c^{-s}W_{g^s_{-\gamma,-\gamma_0}}\mathbf{K}^1_cW_{g^{s-1}_{-c\gamma_0,\gamma}}
     \mathbf{\Lambda}^{-1}_{\gamma}\nonumber\\[2mm]
&&\hskip-15mm=c^{-s}\left[\mathbf{K}^2_c-sc^{-1}\mathbf{K}^1_c\right]W_{g^s_{-\gamma,-\gamma_0}
    g^s_{-c\,\gamma_0,\gamma}}+\mathbf{T}_0
\end{eqnarray}
in $\mathbb{L}_p(\mathbb{R}^+)$ space; $\mathbf{T}_0\;:\;\mathbb{L}_p(\mathbb{R}^+)\to \mathbb{L}_p(\mathbb{R}^+)$ is a compact operator.
\end{theorem}
{\bf Proof:} Let conditions of Theorem \ref{t1.6} hold (that means ${\rm Im}\,\gamma>0$ and
${\rm Im}\,c\,\gamma<0$. Then
 \[
\frac1{(t-c)^2}=\lim_{\varepsilon\to0}\frac1{2\varepsilon i}\left[\frac1{t-c-\varepsilon i}
     -\frac1{t-c+\varepsilon i}\right]
 \]
and we have
\begin{eqnarray*}
\mathbf{\Lambda}^s_{-\gamma}\mathbf{K}^2_c\mathbf{\Lambda}^{-s}_{\gamma}&\hskip-3mm=&\hskip-3mm
     \lim_{\varepsilon\to0}\frac1{2\varepsilon i}\mathbf{\Lambda}^s_{-\gamma} \left[\mathbf{K}^1_{c +\varepsilon i}-\mathbf{K}^1_{c-\varepsilon i}\right] \mathbf{\Lambda}^{-s}_{\gamma}\\
&&\hskip-15mm=\lim_{\varepsilon\to0}\frac1{2\varepsilon i}
     \left[(c+\varepsilon i)^{-s}\mathbf{K}^1_{c +\varepsilon i}\mathbf{\Lambda}^s_{-(c+\varepsilon
     i)\gamma} - (c-\varepsilon i)^{-s}\mathbf{K}^1_{c-\varepsilon i}\mathbf{\Lambda}^s_{-(c-\varepsilon i)\gamma}\right]\mathbf{\Lambda}^{-s}_{\gamma}\\
&&\hskip-15mm=\lim_{\varepsilon\to0}\left\{\frac{(c+\varepsilon i)^{-s}-(c-\varepsilon
     i)^{-s}}{2\varepsilon i}\mathbf{K}^1_{c +\varepsilon i}\mathbf{\Lambda}^s_{-(c+\varepsilon
     i)\gamma}\right.\\
&&\hskip-10mm \left.- (c-\varepsilon i)^{-s}\frac1{2\varepsilon i}\left[
     \mathbf{K}^1_{c+\varepsilon i}-\mathbf{K}^1_{c-\varepsilon i}\right]\mathbf{\Lambda}^s_{-(c-\varepsilon i)\gamma}\right.\\
&&\hskip-10mm \left.- (c-\varepsilon i)^{-s}\mathbf{K}^1_{c-\varepsilon i}\frac1{2\varepsilon i}\left[
     \mathbf{\Lambda}^s_{-(c+\varepsilon i)\gamma}-\mathbf{\Lambda}^s_{-(c-\varepsilon i)\gamma}\right]\right\}\mathbf{\Lambda}^{-s}_{\gamma}\\
&&\hskip-15mm=-s\,c^{-s-1}\mathbf{K}^1_c\mathbf{\Lambda}^s_{-c\,\gamma}\mathbf{\Lambda}^{-s}_{\gamma}
     +c^{-s}\mathbf{K}^2_c\mathbf{\Lambda}^s_{-c\,\gamma}\mathbf{\Lambda}^{-s}_{\gamma}\\
&&\hskip-10mm + c^{-s}\mathbf{K}^1_c\lim_{\varepsilon\to0}
     \cF^{-1}\frac{(\xi-c\,\gamma-\varepsilon\gamma i)^s-(\xi-c\,\gamma+\varepsilon\gamma i)^s}{2\varepsilon i}\cF\mathbf{\Lambda}^{-s}_{\gamma}\\
&&\hskip-15mm=c^{-s}\left[\mathbf{K}^2_c-sc^{-1}\mathbf{K}^1_c\right]W_{g^s_{-c\gamma,\gamma}}
     +s\,\gamma\,c^{-s}\mathbf{K}^1_c\mathbf{\Lambda}^{s-1}_{-c\,\gamma}\mathbf{\Lambda}^{-s}_{\gamma}\\
&&\hskip-15mm=c^{-s}\left[\mathbf{K}^2_c-sc^{-1}\mathbf{K}^1_c\right]W_{g^s_{-c\gamma,\gamma}}
     +s\,\gamma\,c^{-s}\mathbf{K}^1_cW_{g^{s-1}_{-c\gamma,\gamma}}\mathbf{\Lambda}^{-1}_{\gamma}.
 \end{eqnarray*}
Formula \eqref{e72} is proved. To check that the operator $\mathbf{T}$ in \eqref{e73} is compact, let us rewrite it as follows
\begin{eqnarray}\label{e75}
\mathbf{T}=s\,\gamma\,c^{-s}\mathbf{K}^1_cW_{g^{s-1}_{-c\gamma,\gamma}}\mathbf{\Lambda}^{-1}_{\gamma}
    =s\,\gamma\,c^{-s}(1-h)\mathbf{K}^1_cW_{g^{s-1}_{-c\gamma,\gamma}}\mathbf{\Lambda}^{-1}_{\gamma}
    \nonumber\\
    +s\,\gamma\,c^{-s}h\mathbf{K}^1_cW_{g^{s-1}_{-c\gamma,\gamma}}\mathbf{\Lambda}^{-1}_{\gamma},
\end{eqnarray}
where $h\in C^\infty(\mathbb{R})^+$ is a smooth function, equal to 1 in the neighborhood of $0$ and has a compact support. Since $1-h(t)$ vanishes in the neighbourhood of $0$, the operator
$h\mathbf{K}^1_c$ has a smooth kernel and is compact in $\mathbb{L}_p(\mathbb{R}^+)$. The second summand in \eqref{e75} is compact because we can drug $h$ through Mellin $\mathbf{K}^1_c$ and Fourier $W_{g^{s-1}_{-c\gamma,\gamma}}$ convolution operators modulo compact
 \[
s\,\gamma\,c^{-s}h\mathbf{K}^1_cW_{g^{s-1}_{-c\gamma,\gamma}}\mathbf{\Lambda}^{-1}_{\gamma}
      =s\,\gamma\,c^{-s}\mathbf{K}^1_cW_{g^{s-1}_{-c\gamma,\gamma}}h\mathbf{\Lambda}^{-1}_{\gamma}+\mathbf{T}_1,
 \]
where $\mathbf{T}_1$ is compact in $\mathbb{L}_p(\mathbb{R}^+)$  (see Proposition \ref{p3.9} and \cite[Lemma 7.4]{Du79}, \cite[Lemma 1.2]{Du87}) and note, that due to the embedding theorem of Sobolev the operator $h\mathbf{\Lambda}^{-1}_{\gamma}$ is also compact in
$\mathbb{L}_p(\mathbb{R}^+)$, bacause ${\rm supp}\,h$ is compact.

And the final remark: formula \eqref{e74} is derived from \eqref{e72} as in Theorem \ref{t4.1}.
\QED
 %
\begin{remark}\label{r4.4}
The case of operators $\mathbf{K}^n_c$, $n=3,4,\ldots$, can be treated similarly as in Corollary \ref{t4.3}: with the help of perturbation the operator $\mathbf{K}^n_c$ can be represented in the form
\begin{eqnarray} \label{e103}
&&\mathbf{K}^n_c\vf=\lim_{\ve\to0}\mathbf{K}_{c_{1,\varepsilon},\dots,c_{n,\varepsilon}}\vf,\qquad
    \forall\,\vf\in\widetilde{\mathbb{H}}{}^r_p(\mathbb{R}^+)\nonumber\\
&&\mathbf{K}_{c_{1,\varepsilon},\dots,c_{n,\varepsilon}}\vf(t):=\int_0^\infty
    \cK_{c_{1,\varepsilon},\dots,c_{n,\varepsilon}}\left(\frac t\tau\right)\varphi(\tau)\frac{d\tau}\tau=\sum_{j=1}^nd_j(\varepsilon)
    \mathbf{K}^1_{c_{j,\varepsilon}}\vf(t),\nonumber\\
&&\cK_{c_{1,\varepsilon},\dots,c_{m,\varepsilon}}(t):=\frac1{(t-c_{1,\varepsilon})
    \cdots(t-c_{n,\varepsilon})}=\sum_{j=1}^n\frac{d_j(\varepsilon)}{t-c_{j,\varepsilon}},\\
&&c_{j,\varepsilon}\!=\!c(1\!+\!\ve e^{i\omega_j}),\;\; \omega_j\!\in\!(-\pi,\pi), \;\; \arg
    c_{j,\varepsilon},\ \arg c_{j,\varepsilon}\,\gamma_j\!\not=\!0,\;\; j\!=\!1,\dots,m. \nonumber
\end{eqnarray}
The points $\omega_1,\dots,\omega_n\in(-\pi,\pi]$ are pairwise different, i.e.,  $\omega_j\not=\omega_k$ for $j\not=k$ (we remind that $\arg c\not=0$ because $n=3,4,\ldots$). By equating the numerators in the formula \eqref{e103} we find the coefficients $d_1(\varepsilon),\ldots,d_{n-1}(\varepsilon)$.

Since the operators $\mathbf{K}^3_c,\mathbf{K}^4_c,\ldots$ encounter in applications rather rarely, we have confined ourselves with the exact formulae only for the operators $\mathbf{K}^1_c$ and $\mathbf{K}^2_c$.
\end{remark}

\section{Mellin Convolution Operators in Bessel Potential Spaces. Fredholm Properties}
\label{Sec6}

Let us write the symbol of a model operator
\begin{equation}\label{e6.1}
\mathbf{A}:=d_0I+W_{a_0}+\sum_{j=1}^nW_{a_j}\mathbf{K}^1_{c_j}W_{b_j},
\end{equation}
acting in the Bessel potential spaces $\wt{\mathbb{H}}^s_p(\mathbb{R}^+)\to\mathbb{H}^s_p( \mathbb{R}^+)$, compiled of the identity $I$, of Fourier $W_{a_0},\ldots,W_{a_n}$, $W_{b_1},\ldots,W_{b_n}$ and Mellin $\mathbf{K}^1_{c_1},\ldots,\mathbf{K}^1_{c_n}$ convolution operators.

We assume that $a_0,\ldots,a_n,b_1,\ldots,b_n\in C\mathfrak{M}_p(\overline{\mathbb{R}} \setminus\{0\})$, $c_1,\ldots,c_n\in \mathbb{C}$ and, if $s\leqslant\displaystyle\frac1p-1$ or $s\geqslant\displaystyle\frac1p$, the functions $a_1(\xi),\ldots,a_n(\xi)$ have bounded analytic extensions in the lower half plane ${\rm Im}\,\xi<0$, while the functions $b_1(\xi),\ldots,b_n(\xi)$ have bounded analytic extensions in the upper half plane ${\rm Im}\,\xi>0$ to ensure the proper mapping properties of the operator $\mathbf{A}:\wt{\mathbb{H}}^s_p(\mathbb{R}^+)\to \mathbb{H}^s_p(\mathbb{R}^+)$. For $\displaystyle\frac1p-1<s<\displaystyle\frac1p$ such constraints are not necessary.

Now we describe the symbol $\mathcal{A}^s_p(\omega)$ of the operator $\mathbf{A}$. For this we lift the operator $\mathbf{A}:\wt{\mathbb{H}}^s_p(\mathbb{R}^+)\to \mathbb{H}^s_p(\mathbb{R}^+)$ to the $\mathbb{L}_p$-setting and apply equality \eqref{e28}:
\begin{eqnarray}\label{e6.2a}
\mathbf{\Lambda}_{-\gamma}^s\mathbf{A}\mathbf{\Lambda}^{-s}_\gamma&:&
\mathbb{L}_p(\mathbb{R}^+)\to \mathbb{L}_p(\mathbb{R}^+),\\
\mathbf{\Lambda}_{-\gamma}^s\mathbf{A}\mathbf{\Lambda}^{-s}_\gamma
     &\hskip-3mm=&\hskip-3mm d_0\mathbf{\Lambda}_{-\gamma}^s\mathbf{\Lambda}^{-s}_\gamma
     +\mathbf{\Lambda}_{-\gamma}^sW_{a_0}\mathbf{\Lambda}^{-s}_\gamma
     +\sum_{j=1}^nW_{a_j}\mathbf{\Lambda}_{-\gamma}^s\mathbf{K}^1_{c_j}
     \mathbf{\Lambda}^{-s}_\gamma W_{b_j}\nonumber\\
\label{e6.2b}
&\hskip-3mm=&\hskip-3mm d_0W_{\left(\frac{\xi-\gamma}{\xi+\gamma}\right)^s}
     +W_{a_0(\xi)\left(\frac{\xi-\gamma}{\xi+\gamma}\right)^s}
     +\sum_{j=1}^nW_{a_j}\mathbf{K}^1_{c_j}W_{\left(\frac{\xi-c\gamma}{\xi+\gamma}
     \right)^s}W_{b_j}
\end{eqnarray}
(see Theorem \ref{t1.1}, diagram \eqref{e1.25}) if conditions of Theorem \ref{t1.6} hold (see \eqref{e67}) and to the operator
\begin{eqnarray}\label{e6.2c}
\hskip-7mm\mathbf{\Lambda}_{-\gamma}^s\mathbf{A}\mathbf{\Lambda}^{-s}_\gamma
     = d_0W_{\left(\frac{\xi-\gamma}{\xi+\gamma}\right)^s}
     +W_{a_0(\xi)\left(\frac{\xi-\gamma}{\xi+\gamma}\right)^s}
     +\sum_{j=1}^nW_{a_j}\mathbf{K}^1_{c_j}W_{\left(\frac{\xi-\gamma}{\xi-\gamma_0}
     \right)^s\left(\frac{\xi-c\gamma_0}{\xi+\gamma}\right)^s}W_{b_j}+\mathbf{T},
\end{eqnarray}
where $\mathbf{T}\;:\;\mathbb{L}_p(\mathbb{R}^+) \to\mathbb{L}_p(\mathbb{R}^+)$ is a compact operator, if conditions of Corollary \ref{c3.2} hold (see \eqref{e68}).

The symbol of the lifted operator  \eqref{e6.2a}-- \eqref{e6.2c} in $\mathbb{L}_p(\mathbb{R}^+)$-space we declare the symbol of the operator $\mathbf{A}$ in the Bessel potential space. This symbol, written according formulae \eqref{e2.7} and \eqref{e2.8}, has the form:
\begin{equation}\label{e6.2}
\mathcal{A}^s_p(\omega):=d_0\mathcal{I}^s_p(\omega)+\mathcal{W}^s_{a_0,p}(\omega)
     +\sum_{j=1}^n\mathcal{W}^0_{a_j,p}(\omega)\mathcal{K}^{1,s}_{c_j,p}(\omega)
     \mathcal{W}^0_{b_j,p}(\omega),
\end{equation}
where $\mathcal{I}^s_p(\omega)$, $\mathcal{W}^s_{a_0,p}(\omega)$, $\mathcal{W}^0_{a_j,p}(\omega)$, $\mathcal{K}^{1,s}_{c_j,p}(\omega)$ and $\mathcal{W}^0_{b_j,p}(\omega)$ are the symbols of the operators $W_{\left(\frac{\xi-\gamma}{\xi-\gamma}\right)^s}$ in $\mathbb{L}_p$ (of $I$ in $\mathbb{H}^s_p$), of $W_{a_0(\xi)\left(\frac{\xi-\gamma}{\xi-\gamma}\right)^s}$ in $\mathbb{L}_p$ (of $W_{a_0}$ in $\mathbb{H}^s_p$),
of $W_{a_j}$ in $\mathbb{L}_p$ (and in $\mathbb{H}^s_p$) of $\mathbf{K}^1_{c_j}W_{\left(\frac{\xi-c\gamma}{\xi-\gamma}\right)^s}$ in $\mathbb{L}_p$ (of $\mathbf{K}^1_{c_j}$ in $\mathbb{H}^s_p$),  of $W_{b_j}$ in $\mathbb{L}_p$ (and in $\mathbb{H}^s_p$).

Now it suffices to expose the symbols $\mathcal{I}^s_p(\omega)$, $\mathcal{W}^s_{a_0,p}(\omega)$, $\mathcal{W}^0_{a_j,p}(\omega)$ and $\mathcal{K}^{1,s}_{c_j,p}(\omega)$ of the operators $I$, $W_{a_0}$, $W_{a_j} (j=1,2,\ldots,n)$ and $\mathbf{K}^1_c$ separately (the symbol $\mathcal{W}^0_{b_j,p}(\omega)$ of $W_{b_j} (j=1,2,\ldots,n)$ is written analogously):
\begin{subequations}
 \begin{align}\label{e6.3a}
\mathcal{I}^s_p(\omega)&:=\begin{cases}
    g^s_{-\gamma,\gamma,p}(\infty,\xi), & \omega=(\xi,\infty)\in\overline{\Gamma}_1,
    \\[1ex]
\left(\displaystyle\frac{\eta-\gamma}{\eta+\gamma}\right)^{\mp s}, &
     \omega=(+\infty,\eta)\in\Gamma^\pm_2, \\[1ex]
     e^{\pi si}, &\omega=(\xi,0)\in\overline{\Gamma}_3,\end{cases}\\[2ex]
\label{e6.3c}
\mathcal{W}^s_{a,p}(\omega)&:=\begin{cases}
        a^s_p(\infty,\xi), & \omega=(\xi,\infty)\in\overline{\Gamma}_1, \\
        a(\mp\eta)\left(\displaystyle\frac{\eta-\gamma}{\eta+\gamma}\right)^{\mp s}, & \omega=(+\infty,\eta)\in\Gamma^\pm_2, \\
        e^{\pi si}a_p(0,\xi), &
        \omega=(\xi,0)\in\overline{\Gamma}_3,\end{cases}\\[2ex]
\label{e6.3b}
\mathcal{W}^0_{a,p}(\omega)&:=\begin{cases}
        a_p(\infty,\xi), & \omega=(\xi,\infty)\in\overline{\Gamma}_1, \\
        a(\mp\eta), & \omega=(+\infty,\eta)\in\Gamma^\pm_2, \\
        a_p(0,\xi), &
        \omega=(\xi,0)\in\overline{\Gamma}_3,\end{cases}\\[2ex]
 \label{e6.3d}
\mathcal{K}^{1,s}_{c,p}(\omega)&:=\begin{cases}
    \displaystyle\frac{c^{-s}(-c)^{\frac1p-i\xi-1}}{\sin\pi(\frac1p-i\xi)},
    &\omega=(\xi,\infty)\in\overline{\Gamma}_1,,\\[1ex]
    0, &\omega=(\pm\infty,\eta)\in\Gamma^\pm_2, \\[1ex]
    \displaystyle\frac{c^{-s}(-c)^{\frac1p+s-i\xi-1}}{\sin\pi(\frac1p-i\xi)},
    &\omega=(\xi,0)\in\overline{\Gamma}_3,\qquad \text{for} \quad \arg\,c\not=0,\end{cases}\\[1.5ex]
\label{e6.3e}
\mathcal{K}^{1,s}_{1,p}(\omega)&:=\begin{cases}
    -i\cot\pi(\frac1p-i\xi),&\omega=(\xi,\infty)\in\overline{\Gamma}_1,,\\[1ex]
    \pm1, &\omega=(\pm\infty,\eta)\in\Gamma^\pm_2, \\[1ex]
    i\cot\pi(\frac1p-i\xi),&\omega=(\xi,0)
    \in\overline{\Gamma}_3,\end{cases}\\[1.5ex]
a^s_p(\infty,\xi)&:=\frac{e^{2\pi si}a(\infty)+a(-\infty)}2+\frac{e^{2\pi si}
    a(\infty)-a(-\infty)}{2i}\cot\pi\Big(\frac1p-i\xi\Big),\nonumber
 \end{align}
 \begin{align*}
&\hskip-7mm a_p(x,\xi):=\frac{a(x+0)+a(x-0)}2+\frac{a(x+0)-a(x-0)}{2i}
     \cot\pi\Big(\frac1p-i\xi\Big),\quad x=0,\infty,\\[1.5ex]
&\hskip-7mm g^s_{-\gamma,\gamma,p}(\infty,\xi):=\frac{e^{2\pi si}+1}2
     +\frac{e^{2\pi si}-1}{2i}\cot\pi\Big(\frac1p-i\xi\Big)=e^{\pi si}\frac{\sin\pi\Big(\frac1p+s-i\xi\Big)}
     {\sin\pi\Big(\frac1p-i\xi\Big)},\\
&\hskip88mm\xi\in\mathbb{R},\quad \eta\in\mathbb{R}^+,\nonumber
 \end{align*}
where
 \[
-\pi\leqslant\arg\,c<\pi,\quad-\pi<\arg(c\,\gamma)<0,\quad
0<\arg\gamma<\pi
 \]
and $c^s=|c|^se^{is\arg\,c}$, $(-c)^\delta=e^{\pm\pi i\delta}|c|^\delta e^{i\delta\arg\,c}$ for $c,\delta\in\mathbb{C}$; the sign "+" is chosen for $\pi<\arg\,c<2\pi$ and the sign "-" is chosen for $0<\arg\,c<\pi$.

Note that, we got the equal symbol $\mathcal{K}^{1,s}_{c,p}(\omega)$ of the operator $\mathbf{K}^1_{c_j}$ in the cases \eqref{e6.2b} and \eqref{e6.2c} because the functions
 \[
g^s_{-\gamma,-\gamma_0}(\xi)g^s_{-c\gamma_0,\gamma}(\xi):=\left(\frac{\xi-\gamma}{\xi
     -\gamma_0}\right)^s\left(\frac{\xi-c\gamma_0}{\xi+\gamma}\right)^s\quad{\rm and}\quad g^s_{-c\gamma,\gamma}(\xi):=\left(\frac{\xi-c\gamma}{\xi+\gamma}
     \right)^s
 \]
have equal limits at infinity $g^s_{-c\gamma,\gamma}(\pm\infty)=g^s_{-\gamma,-\gamma_0} (\pm\infty)g^s_{-c\gamma_0, \gamma}(\pm\infty)=1$ and\linebreak $g^s_{-c\gamma,\gamma}(0) =g^s_{-\gamma,-\gamma_0}(0)g^s_{-c\gamma_0, \gamma}(0)=(-c)^s$.

If $a(-\infty)=1$ and $a(+\infty)=e^{2\pi\alpha i}$, then $a_\infty^-=0$, $a_\infty^+=2\alpha$ and the symbol $a^s_p(\infty,\xi)$ acquires the form:
 \begin{equation}\label{e6.3h}
a^s_p(\infty,\xi)=e^{\pi(s+\alpha)i}\frac{\sin\pi\left(\frac1p+s+\alpha-i\xi\right)}{
     \cos\pi\Big(\frac1p-i\xi\Big)}.
 \end{equation}
\end{subequations}

Note that, the Mellin convolution operator
 \[
\begin{array}{c}
\mathbf{K}^1_{-1}\varphi(t):=\int\limits_0^\infty\displaystyle\frac{\varphi(\tau)
     \,d\tau}{t+\tau}=\mathfrak{M}^0_{\cM_\frac1p\cK^1_{-1}},\qquad
\cM_\frac1p\cK^1_{-1}(\xi)=\displaystyle\frac{\pi d\,^{\beta-i\xi-1}}{\sin\pi(\beta-i\xi)}
 \end{array}
 \]
(see \eqref{e3.11d}), which we encounter in applications, has a rather simple symbol in the Bessel potential space $\mathbb{H}^s_p(\mathbb{R}^+)$ (see \eqref{e6.3c}, where $c=-1$):
 \begin{eqnarray*}
\hskip-5mm\mathcal{K}^{1,s}_{-1,p}(\omega):=\begin{cases}
     \displaystyle\frac{e^{\pi si}}{\sin\pi(\beta-i\xi)},&\omega=(\xi,\infty) \in\overline{\Gamma_1} \cup\overline{\Gamma_3},\\
     0, &\omega=(\pm\infty,\eta)\in\Gamma^\pm_2,. \end{cases}
 \end{eqnarray*}
 %
\begin{theorem}\label{t6.1}
Let $1<p<\infty$, $s\in\mathbb{R}$. The operator
 \begin{equation}\label{e6.4}
\mathbf{A}:\widetilde{\mathbb{H}}{}^s_p(\mathbb{R}^+)\longrightarrow
     \mathbb{H}^s_p (\mathbb{R}^+),
 \end{equation}
defined in \eqref{e6.1}, is Fredholm if and only if its symbol $\mathcal{A}^s_p(\omega)$ defined in \eqref{e6.2} and \eqref{e6.3a}--\eqref{e6.3h}, is elliptic.

If $\mathbf{A}$ \, is Fredholm, the index of the operator has the value
 \[
    \Ind\mathbf{A}=-\ind\det\mathcal{A}^s_p.
 \]
\end{theorem}
{\bf Proof:}
Let $c_j, d_j\in\mathbb{C}$, $-\pi\leqslant\arg\,c_j<\pi$
$\arg\,c_j\not=0$, for $j=1,\ldots,n$. Lifting the operator $\mathbf{A}$ to the $\mathbb{L}_p(\mathbb{R}^+)$ space we get
\begin{eqnarray}\label{eqn51a}
\mathbf{\Lambda}^s_{-\gamma}\mathbf{A}\mathbf{\Lambda}^{-s}_{\gamma}
     =d_0\mathbf{\Lambda}^s_{-\gamma}\mathbf{\Lambda}^{-s}_{\gamma}
     +\mathbf{\Lambda}^s_{-\gamma}W_{a_0}\mathbf{\Lambda}^{-s}_{\gamma}
+\sum_{j=1}^nW_{a_j}\mathbf{\Lambda}^s_{-\gamma}\mathbf{K}^1_{c_j}
     \mathbf{\Lambda}^{-s}_{\gamma}W_{b_j},
\end{eqnarray}
where $c^{-s}=|c|^{-s}e^{-is\arg\,c}$ and $\gamma_j$ are such that $0<\arg
\gamma_j<\pi$, $-\pi<\arg(c_j\,\gamma_j)<0$ for $j=m+1,\ldots,n$.

To derive \eqref{eqn51a} we have applied the following property of convolution operators $\mathbf{\Lambda}^s_{-\gamma}W_{a_j}=W_{a_j}\mathbf{\Lambda}^s_{-\gamma}$
and $W_{b_j}\mathbf{\Lambda}^s_\gamma=\mathbf{\Lambda}^s_\gamma W_{b_j}$, $\mathbf{\Lambda}^{\mp s}_{\pm\gamma}=W_{\lambda^{\mp s}_{\pm\gamma}}$,
which are based on the analytic extension properties of the symbols $\lambda^s_{-\gamma}, a_1(\xi),\ldots,a_n(\xi)$ in the lower half plane ${\rm Im}\,\xi<0$ and of symbols $\lambda^{-s}_\gamma,b_1(\xi),\ldots,b_n(\xi)$ in the upper half plane ${\rm Im}\,\xi>0$ (see \eqref{e1.24}).

The model operators $I$, $\mathbf{K}^1_c$ and $W_{a}$ Lifted to the space $\mathbb{L}_p(\mathbb{R}^+)$ acquire the form
 \begin{eqnarray}\label{e6.5}
&&\hskip0mm\mathbf{\Lambda}^s_{\gamma}I\mathbf{\Lambda}^{-s}_{\gamma}=W_{g^s_{-\gamma,\gamma}},
     \qquad \mathbf{\Lambda}^s_{\gamma}W_{a_k}\mathbf{\Lambda}^{-s}_{\gamma}
     =W_{a_kg^s_{-\gamma,\gamma}},\nonumber\\[2mm]
&&\hskip-10mm\mathbf{\Lambda}^s_{\gamma}\mathbf{K}^1_c\mathbf{\Lambda}^{-s}_{\gamma}
     =\left\{\begin{array}{ll}c^{-s}\mathbf{K}^1_c W_{g^s_{-c\,\gamma,\gamma}}\! &\text{for}\; -\pi<\arg(c\,\gamma)<0,\\[2mm]
     c^{-s}\mathbf{K}^1_c W_{g^s_{-\,\gamma,-\gamma_0}g^s_{-c\,\gamma_0,\gamma}}+\mathbf{T} &\text{for}\quad 0<\arg(c\,\gamma)<\pi,\\
     &\hskip3mm -\pi<\arg(c\,\gamma_0)|<0,
     \end{array}\right.
 \end{eqnarray}
where $\mathbf{T}$ is a compact operator. Here, as above, $-\pi\leqslant\arg\,c<\pi$, $\arg\,c\not=0$, $0<\arg\,\gamma<\pi$, $0<\arg\,\gamma_0<\pi$ and either $-\pi<\arg(c\,\gamma)<0$ or, if $-\pi<\arg(c\,\gamma)<0$, then $-\pi<\arg(c\,\gamma_0)|<0$. Here $c^{-s}=|c|^{-s}e^{-is\arg\,c}$.

Therefore the operator $\mathbf{\Lambda}^s_{-\gamma}\mathbf{A}\mathbf{\Lambda}^{-s}_{ \gamma}$ in \eqref{eqn51a} is rewritten as follows:
\begin{eqnarray}\label{eqn51b}
\mathbf{\Lambda}^s_{-\gamma}\mathbf{A}\mathbf{\Lambda}^{-s}_{\gamma}
     =d_0W_{g^s_{-\gamma,\gamma}}+W_{a_0g^s_{\gamma,\gamma}}+\sum_{j=1}^mc^{-s}_j
     W_{a_j}\mathbf{K}^1_{c_j}W_{g^s_{-c_j\gamma,\gamma}}W_{b_j}\nonumber\\
+\sum_{j=m+1}^nc^{-s}_jW_{a_j}\mathbf{K}^1_{c_j}W_{g^s_{-\gamma,-\gamma_j}
     g^s_{-c_j\gamma_j,\gamma}}W_{b_j}+\mathbf{T}\;:\;\mathbb{L}_p(\mathbb{R}^+)
     \longrightarrow\mathbb{L}_p(\mathbb{R}^+),
\end{eqnarray}
where $\mathbf{T}$ is a compact operator and we ignore it when writing the symbol of $\mathbf{A}$.

We declare the symbol of the lifted operator $\mathbf{\Lambda}^s_{-\gamma}\mathbf{A}\mathbf{\Lambda}^{-s}_{\gamma}$ (see \eqref{eqn51b}) in the Lebes\-gue space $\mathbb{L}_p(\mathbb{R}^+)$ as the symbol of the initial operator $\mathbf{A}\;:\;\widetilde{\mathbb{H}}^s_p(\mathbb{R}^+) \rightarrow\mathbb{H}^s_p(\mathbb{R}^+)$ in \eqref{e6.1}.

The function $g^s_{-\gamma,\gamma}\in C(\mathbb{R})$ is continuous on $\mathbb{R}$, but has different limits at the infinity
\begin{subequations}
 \begin{eqnarray}\label{e6.6a}
 g^s_{-\gamma,\gamma}(-\infty)=1, \quad g^s_{-\gamma,\gamma}(+\infty)=e^{2\pi si},\qquad
       g^s_{-\gamma,\gamma}(0)=e^{\pi si},
 \end{eqnarray}
while the functions $g^s_{-\gamma,-\gamma_0},\ g^s_{-c\gamma,\gamma}, g^s_{-c\gamma_0,\gamma}\in C(\mathbb{R})$ are continuous on $\mathbb{R}$ including the infinity
 \begin{eqnarray}\label{e6.6b}
 \begin{array}{c}
g^s_{-c\,\gamma,\gamma}(\pm\infty)=g^s_{-\gamma,-\gamma_0}(\pm\infty)
    =g^s_{-c\,\gamma_0,\gamma}(\pm\infty)=1,\\[2mm]
g^s_{-\gamma,-\gamma_0}(0)g^s_{-c\,\gamma_0,\gamma}(0)
     =\left(\displaystyle\frac{-\gamma}{-\gamma_0}\right)^s
     \left(\displaystyle\frac{-c\gamma_0}{\gamma}\right)^s
     =\left(-c\right)^s,\\[2mm]
g^s_{-c\,\gamma,\gamma}(0)=\left(-c\right)^s\quad\text{if}\quad
     -\pi\leqslant\arg\,c<\pi, \quad \arg\,c\not=0.
 \end{array}
 \end{eqnarray}
\end{subequations}

In the Lebesgue space $\mathbb{L}_p(\mathbb{R}^+)$, the symbols of the first two operators in \eqref{eqn51b}, are written according the formulae \eqref{e2.7}--\eqref{e2.8} by taking into account the equalities \eqref{e6.6a} and \eqref{e6.6a}. The symbols of these operators have, respectively, the form \eqref{e6.3a} and \eqref{e6.3c}.

The symbols of operators $W_{a_1},\ldots,W_{a_n}$ and $W_{b_1},\ldots,W_{b_n}$ are written with the help of the formulae \eqref{e2.7}--\eqref{e2.8} and have the
form \eqref{e6.3b}.

The lifted Mellin convolution operators
 \[
\mathbf{\Lambda}^s_{\gamma}\mathbf{K}^1_{c_j}\mathbf{\Lambda}^{-s}_{\gamma}
     \;:\;\mathbb{L}_p(\mathbb{R}^+)\longrightarrow,\mathbb{L}_p(\mathbb{R}^+)
 \]
 are of mixed type and comprise both the Mellin convolution operators $\mathbf{K}^1_{c_j}=\mathfrak{M}^0_{\mathcal{K}^1_{c_j,p}(\xi)}$, where the symbol $\mathcal{K}^1_{c_j,p}(\xi):=\mathcal{M}_{1/p}\mathcal{K}^1_{c_j}(\xi)$ is defined in \eqref{e3.11d} and \eqref{e3.11e}, and the Fourier convolution operators $W_{g^s_{-c_j\,\gamma_0,\gamma}}$ and $W_{g^s_{-\,\gamma,-\gamma_0} g^s_{-c_j\,\gamma_0, \gamma}}$. The symbol of the operators $\mathbf{\Lambda}^s_{\gamma}\mathbf{K}^1_{c_j} \mathbf{\Lambda}^{-s}_{\gamma}$ from \eqref{e6.5} in the Lebesgue space $\mathbb{L}_p(\mathbb{R}^+)$ is found according formulae \eqref{e2.7}--\eqref{e2.8}, has the form \eqref{e6.3d} and is declared the symbol of $\mathbf{K}^1_{c_j}\;:\;\widetilde{\mathbb{H}}^s_p(\mathbb{R}^+) \rightarrow\mathbb{H}^s_p(\mathbb{R}^+)$. The symbols of Fourier convolution factors $W_{g^s_{-c_j\,\gamma_0,\gamma}}$ and $W_{g^s_{-\,\gamma,-\gamma_0}g^s_{-c_j\,\gamma_0, \gamma}}$, which contribute the symbol of $\mathbf{K}^1_{c_j}=\mathfrak{M}^0_{ \mathcal{K}^1_{c_j,p}}$ are written again according formulae \eqref{e2.7}--\eqref{e2.8} by taking into account the equalities \eqref{e6.6a} and \eqref{e6.6b}.

To the lifted operator applies Theorem \ref{t4.2} and gives the result formulated in
Theorem \ref{t6.1}.
\QED
 %
\begin{corollary}\label{c6.2}
Let $1<p<\infty$, $s\in\mathbb{R}$. The operator
 \begin{equation*}
\mathbf{A}:\widetilde{\mathbb{H}}{}^s_p(\mathbb{R}^+)\longrightarrow
     \mathbb{H}^s_p (\mathbb{R}^+),
 \end{equation*}
defined in \eqref{e2.15}, is locally invertible at $0\in\mathbb{R}^+$ if and only if its symbol $\mathcal{A}^s_p(\omega)$, defined in \eqref{e6.2} and \eqref{e6.3a}--\eqref{e6.3h}, is elliptic on $\Gamma_1$, i.e.
 \[
\inf_{\omega\in\Gamma_1}\left|\det\,\mathcal{A}^s_p(\omega)\right|
     =\inf_{\xi\in\mathbb{R}}\left|\det\,\mathcal{A}^s_p(\xi,\infty)\right|>0.
 \]
\end{corollary}
{\bf Proof:}

\QED

For the definition of the Sobolev--Slobodeckij (Besov) spaces $\mathbb{W}^s_p(\Omega) =\mathbb{B}^s_{p,p}(\Omega)$, $\widetilde{\mathbb{W}}^s_p(\Omega)=\widetilde{ \mathbb{B}}{}^s_{p,p}(\Omega)$ for an arbitrary domain $\Omega\subset\mathbb{R}^n$, including the half axes $\mathbb{R}^+$, we refer to the monograph \cite{Tr95}.
 %
\begin{corollary}\label{c6.3}
Let $1<p<\infty$, $s\in\mathbb{R}$. If the operator $\mathbf{A}:\widetilde{\mathbb{H}}{}^s_p(\mathbb{R}^+) \longrightarrow\mathbb{H}^s_p (\mathbb{R}^+)$, defined in \eqref{e2.15}, is Fredholm $($is invertible$)$ for all $a\in(s_0,s_1)$ and $p\in(p_0,p_1)$, where $-\infty<s_0<s_1 <\infty$,  $1<p_o<p_1<\infty$, then
\begin{equation}\label{e6.7}
    \mathbf{A}:\widetilde{\mathbb{W}}{}^s_p (\mathbb{R}^+) \longrightarrow \mathbb{W}^s_p (\mathbb{R}^+),\;\;\; s\in(s_0,s_1), \;\; p\in(p_0,p_1)
\end{equation}
is Fredholm and has the equal index
\begin{equation}\label{e6.8}
    \Ind\mathbf{A}=-\ind\det\mathcal{A}^s_p.
\end{equation}
$($is invertible, respectively$)$ in the Sobolev--Slobodeckij $($Besov$)$ spaces $\mathbb{W}^s_p =\mathbb{B}^s_{p,p}$.
\end{corollary}

{\bf Proof:}
First of all recall that the Sobolev--Slobodeckij (Besov) spaces $\mathbb{W}^s_p=\mathbb{B}^s_{p,p}$ emerge as the result of interpolation with the real interpolation method between Bessel potential spaces

\begin{equation}\label{e6.9}
\begin{aligned}
\big(\bH_{p_0}^{s_0}(\Omega),\bH_{p_1}^{s_1}(\Omega)\big)_{\theta,p} & =\bW_p^s(\Omega),\;\;
    s:=s_0(1-\theta)+s_1\theta, \\
\big(\widetilde{\bH}{}_{p_0}^{s_0}(\Omega),\widetilde{\bH}{}_{p_1}^{s_1}(\Omega)\big)_{\theta,p}
    & =\widetilde{\bW}{}_p^s(\Omega),\;\;p:=\frac1{p_0}\,(1-\theta)+\frac1{p_1}\,\theta, \;\; 0<\theta<1.
\end{aligned}
\end{equation}

If $\mathbf{A}:\widetilde{\mathbb{H}}{}^s_p (\mathbb{R}^+) \longrightarrow\mathbb{H}^s_p (\mathbb{R}^+)$ is Fredholm (or is invertible) for all $s\in(s_0,s_1)$ and $p\in(p_0,p_1)$, it has a regularizer $\mathbf{R}$ (has the inverse $\mathbf{A}^{-1}=\mathbf{R}$, respectively), which is bounded in the setting
 \[
\mathbf{R}:\mathbb{W}^s_p (\mathbb{R}^+) \longrightarrow\widetilde{\mathbb{W}}{}^s_p (\mathbb{R}^+)
 \]
due to the interpolation \eqref{e6.9} and
 \[
\mathbf{R}\mathbf{A}=I+\mathbf{T}_1, \quad \mathbf{A}\mathbf{R}=I+\mathbf{T}_2,
 \]
where $\mathbf{T}_1$ and $\mathbf{T}_2$ are compact in $\widetilde{\mathbb{H}}{}^s_p (\mathbb{R}^+)$ and in $\mathbb{H}^s_p(\mathbb{R}^+)$, respectively ($\mathbf{T}_1=\mathbf{T}_2=0$ if $\mathbf{A}$ is invertible).

Due to the Krasnoselskij interpolation theorem (see \cite{Tr95}), $\mathbf{T}_1$ and $\mathbf{T}_2$ are compact in $\widetilde{\mathbb{W}}{}^s_p(\mathbb{R}^+)$ and in $\mathbb{W}^s_p (\mathbb{R}^+)$, respectively for all $s\in(s_0,s_1)$ and $p\in(p_0,p_1)$ and, therefore, $\mathbf{A}$ in \eqref{e6.7} is Fredholm (is invertible, respectively).

The index formulae \eqref{e6.8} follows from the embedding properties of the Sobolev--Slobodeckij and Bessel potential spaces by standard well-known arguments.                       \QED


\begin{thebibliography}{99}
\addcontentsline{toc}{section}{References}

\bibitem{BCC12a} \textsc{A.-S. Bonnet-Ben Dhia, L. Chesnel, P, Ciarlet, Jr.} $T$-coercivity for scalar interface problems between dielectrics and metamaterials. \emph{ESAIM Math. Model. Numer. Anal.}  \textbf{46}  (2012), No. 6, 1363--1387.

\bibitem{BCC12b} \textsc{A.-S. Bonnet-Ben Dhia, L. Chesnel, X. Claeys}, Radiation condition for a non-smooth interface between a dielectric and a metamaterial. \emph{Math. Models Methods Appl. Sci.}  \textbf{23}  (2013), No. 9, 1629--1662; \verb"http://hal.inria.fr/hal-00651008/"

\bibitem{BDKT13} \textsc{T. Buchukuri, R. Duduchava, D. Kapanadze, M. Tsaava}, The boundary value problems for the Helmholtz equation in arbitrary $2D$-sectors. \emph{Manuscript}.

\bibitem{Cos83} \textsc{M. Costabel}, Boundary integral operators on curved polygons. \emph{Ann. Mat. Pura Appl. (4)} \textbf{133} (1983), 305--326.

\bibitem{CS84} \textsc{M. Costabel, E. Stephan}, The method of Mellin transformation for boundary integral equations on curves with corners. In: A. Gerasoulis, R. Vichnevetsky (Eds.), {\em Numerical Solutions of Singular Integral Equations, IMACS, New Brunswick,} 1984, pp. 95--102.

\bibitem{Co69} \textsc{H. O. Cordes}, Pseudo-differential operators on a half-line. \emph{J. Math. Mech.}  \textbf{18}  (1968/69), 893--908.

\bibitem{DiDu15} \textsc{V. D. Didenko, R. Duduchava}, Mellin convolution operators in the Bessel potential spaces. Submitted for a publication.\\ Preprint: http://arxiv.org/pdf/1502.02756.pdf

\bibitem{DiSi08} \textsc{V. D. Didenko, B. Silbermann}, Approximation of additive convolution-like operators. Real $C^{*}$-algebra approach. \emph{Frontiers in Mathematics. Birkh\"{a}user Verlag, Basel}, 2008.

\bibitem{DiRS95} \textsc{V. D. Didenko, S. Roch, B. Silbermann}, Approximation methods for singular integral equations with conjugation on curves with corners. \emph{SIAM J. Numer. Anal.}  \textbf{32}  (1995), No. 6, 1910--1939.

\bibitem{DiV96} \textsc{V. Didenko, E. Venturino}, Approximate solutions of some Mellin equations with conjugation. \emph{Integral Equations Operator Theory}  \textbf{25}  (1996), No. 2, 163--181.

\bibitem{Du74} \textsc{R. Duduchava}, Convolution integral operators with discontinuous coefficients. (Russian) \emph{Dokl. Akad. Nauk SSSR}  \textbf{218}  (1974), 264--267; translation in {\em Sov. Math. Doklady} {\bf 15} (1974), 1302--1306. 

\bibitem{Du75a} \textsc{R. V. Duduchava}, On Wiener--Hopf integral operators. \emph{Math. Nachr.} \textbf{65} (1975), No. 1, 59--82.

\bibitem{Du75b} \textsc{R. V. Duduchava}, Convolution integral operators with discontinuous symbols. (Russian) \emph{Tr. Tbilis. Mat. Inst. Razmadze} \textbf{50} (1975), 34--41.

\bibitem{Du77} \textsc{R. V. Duduchava}, Integral operators of convolution type with disconnected coefficients. \emph{Math. Nachr.} \textbf{79} (1977), 75--98.

\bibitem{Du78} \textsc{R. V. Duduchava}, Integral equations of convolution type with discontinuous coefficients. (Russian) \emph{Soobshch. Akad. Nauk Gruzin. SSR}  \textbf{92}  (1978), No. 2, 281--284.

\bibitem{Du79} \textsc{R. Duduchava}, Integral equations in convolution with discontinuous presymbols, singular integral equations with fixed singularities, and their applications to some problems of mechanics. \emph{BSB B. G. Teubner Verlagsgesellschaft, Leipzig}, 1979.

\bibitem{Du82} \textsc{R. Duduchava}, An application of singular integral equations to some problems of elasticity. \emph{Integral Equations Operator Theory}  \textbf{5}  (1982), No. 4, 475--489.

\bibitem{Du84a} \textsc{R. Duduchava}, On multidimensional singular integral operators. I. The half-space case. \emph{J. Operator Theory}  \textbf{11}  (1984), Nno. 1, 41--76; II. The case of compact manifolds. \emph{J. Operator Theory}  \textbf{11}  (1984), No. 2, 199--214.

\bibitem{Du84b} \textsc{R. Duduchava}, On general singular integral operators of the plane theory of elasticity. \emph{Rend. Sem. Mat. Univ. Politec. Torino}  \textbf{42}  (1984), No. 3, 15--41.

\bibitem{Du86} \textsc{R. Duduchava}, General singular integral equations and fundamental problems of the plane theory of elasticity. (Russian) \emph{Trudy Tbiliss. Mat. Inst. Razmadze Akad. Nauk Gruzin. SSR} \textbf{82} (1986), 45--89.

\bibitem{Du87} \textsc{R. Duduchava}, On algebras generated by convolutions and discontinuous functions. \emph{Special issue: Wiener--Hopf problems and applications (Oberwolfach, 1986). Integral Equations Operator Theory}  \textbf{10}  (1987), No. 4, 505--530.

\bibitem{Du13} \textsc{R. Duduchava}, Mellin convolution operators in Bessel potential spaces with admissible meromorphic kernels, {\em Memoirs on Differential Equations and Mathematical Physics} {\bf 60},  135-177, 2013.

\bibitem{DL85} \textsc{R. Duduchava, T. Latsabidze}, The index of singular integral equations with complex-conjugate functions on piecewise-smooth lines. (Russian) \emph{Trudy Tbiliss. Mat. Inst. Razmadze Akad. Nauk Gruzin. SSR} \textbf{76} (1985), 40--59.

\bibitem{DLS95} \textsc{R. Duduchava, T. Latsabidze, A.  Saginashvili}, Singular integral operators with the complex conjugation on curves with cusps. \textsc{Integral Equations Operator Theory}  \textbf{22}  (1995), No. 1, 1--36.

\bibitem{DS93} \textsc{R. Duduchava, F.-O. Speck}, Pseudodifferential operators on compact manifolds with Lipschitz boundary. \emph{Math. Nachr.} \textbf{160} (1993), 149--191.

\bibitem{DTT14} \textsc{R. Duduchava, M. Tsaava, T. Tsutsunava},
 Mixed boundary value problem on hypersurfaces. {\em International Journal of Differential Equations}, Hindawi Publishing Corporation, Volume 2014, Article ID 245350, 8 pages.

\bibitem{Es81} \textsc{G. I. Eskin}, Boundary value problems for elliptic pseudodifferential equations. Translated from the Russian by S. Smith. \emph{Translations of Mathematical Monographs}, 52. \emph{American Mathematical Society, Providence, R.I.}, 1981.

\bibitem{GF74} \textsc{I. Z. Gochberg, I. A. Feldman}, Faltungsgleichungen und Projektionsverfahren zu ihrer Losung. (German) \"{U}bersetzung aus dem Russischen von Reinhard Lehmann und J\"{u}rgen Leiterer. \emph{Mathematische Reihe}, Band 49. \emph{Birkh\"{a}user Verlag, Basel--Stuttgart}, 1974.

\bibitem{GK79} \textsc{I. Gohberg, N. Krupnik}, One-dimensional linear singular integral equations. I, II. \emph{Operator Theory: Advances and Applications} 53--54. \emph{Birkh\"{a}user Verlag, Basel etc.}, 1992.

\bibitem{GR94} \textsc{I. S. Gradshteyn, I. M. Ryzhik}, Table of integrals, series, and products. \emph{Academic Press, Inc. Boston, MA}, 1994.

\bibitem{GB10} \textsc{D. K. Gramotnev, S. I. Bozhevolnyi}, Plasmonics beyond the diffraction limit. \emph{Nature Photonics} \textbf{4} (2010), No. 2, 83--91.

\bibitem{Kh57} \textsc{B. V. Khvedelidze}, Linear discontinuous boundary problems in the theory of functions, singular integral equations and some of their applications. (Russian)\emph{ Akad. Nauk Gruzin. SSR. Trudy Tbiliss. Mat. Inst. Razmadze}  \textbf{23} (1956), 3--158.

\bibitem{Kr60} \textsc{M. Krasnosel'skij},On a theorem of M. Riesz, \emph{ Sov. Math. Dokl.} {\bf 1}, 1960, 229-231; translation from \emph{ Dokl. Akad. Nauk SSSR} {\bf 131}, 1960,  246-248.

\bibitem{Ku04} \textsc{P. Kuchment}, Quantum graphs. I. Some basic structures. \emph{Waves Random Media}  \textbf{14}  (2004), No. 1, S107--S128.

\bibitem{Ku05} \textsc{P. Kuchment}, Quantum graphs. II. Some spectral properties of quantum and combinatorial graphs. \emph{J. Phys. A}  \textbf{38}  (2005), No. 22, 4887--4900.

\bibitem{RR12} \textsc{V. S. Rabinovich, S. Roch}, Pseudodifferential operators on periodic graphs. \emph{Integral Equations Operator Theory}  \textbf{72}  (2012), No. 2, 197--217.

\bibitem{Sc85} \textsc{R. Schneider}, Integral equations with piecewise continuous coefficients in $L^p$-spaces with weight. \emph{J. Integral Equations}  \textbf{9}  (1985), No. 2, 135--152.

\bibitem{Se66} \textsc{R.T. Seeley}, Singular integrals and boundary value problems {\em American Journal of Mathematics} {\bf 88}, 1966, 781-809.

\bibitem{Si65} \textsc{I. B. Simonenko}, A new general method of investigating linear operator equations of singular integral equation type. I. (Russian) \emph{Izv. Akad. Nauk SSSR Ser. Mat.}  \textbf{29}  (1965), 567--586.

\bibitem{Th85} \textsc{G. Thelen}, Zur Fredholmtheorie singul\"arer Integrodifferentialoperatoren auf der Halbachse. \emph{Dissertation Dr. rer. nat. Darmstadt}, 1985.

\bibitem{Tr95} \textsc{H. Triebel}, Interpolation theory, function spaces, differential operators. \emph{Johann Ambrosius Barth, Heidelberg}, 1995.
\end{thebibliography}
\end{document}